\documentclass[dvipsnames,11pt,reqno,twoside,final]{amsart}
\usepackage[a4paper,left=3cm,right=3cm,top=3cm,bottom=3cm]{geometry}
\usepackage[T1]{fontenc}
\usepackage[utf8]{inputenc}
\usepackage{xcolor}
\usepackage{amssymb,mathtools}  
\usepackage{dsfont}
\usepackage{enumitem}
\usepackage{subcaption}
\usepackage{centernot}
\usepackage{booktabs}
\usepackage{enumerate}
\usepackage{mathrsfs}
\usepackage{mathtools}
\mathtoolsset{showonlyrefs}

\usepackage{preamble/mhequ}

\usepackage{amsthm}

\newcounter{counterEnvDefault}
\numberwithin{counterEnvDefault}{section}

\theoremstyle{plain}

\newtheorem{lemma}[counterEnvDefault]{Lemma}
\newtheorem*{lemma*}{Lemma}

\newtheorem{theorem}[counterEnvDefault]{Theorem}
\newtheorem*{theorem*}{Theorem}

\newtheorem{proposition}[counterEnvDefault]{Proposition}
\newtheorem*{proposition*}{Proposition}

\newtheorem{corollary}[counterEnvDefault]{Corollary}
\newtheorem*{corollary*}{Corollary}

\newtheorem{assumption}[counterEnvDefault]{Assumption}
\newtheorem*{assumption*}{Assumption}

\theoremstyle{definition}

\newtheorem*{exercise*}{Exercise}

\newtheorem{definition}[counterEnvDefault]{Definition}
\newtheorem*{definition*}{Definition}

\newtheorem*{notation*}{Definition}

\newtheorem{remark}[counterEnvDefault]{Remark}
\newtheorem*{remark*}{Remark}

\newtheorem{claim}[counterEnvDefault]{Claim}
\newtheorem*{claim*}{Claim}

\newtheorem*{assertion*}{Assertion}

\usepackage{microtype}
\renewcommand\epsilon\varepsilon

\usepackage{hyperref}
\definecolor{colorlinks}{RGB}{0, 24, 168}
\definecolor{colorcites}{RGB}{0, 138, 118}
\hypersetup{
    colorlinks=true,
    linkcolor=colorlinks,
    citecolor=colorcites,
    urlcolor=colorlinks,
    pdfborder={0 0 0}
}

\usepackage{xargs}
\usepackage[colorinlistoftodos,prependcaption,textsize=tiny]{todonotes}

\newcommandx\work[2][1=]{\todo[linecolor=RoyalBlue,backgroundcolor=RoyalBlue!25,bordercolor=RoyalBlue,#1]{\textsc{todo} #2}}
\newcommandx\comment[2][1=]{\todo[linecolor=OliveGreen,backgroundcolor=OliveGreen!25,bordercolor=OliveGreen,#1]{\textsc{comment} #2}}
\newcommandx\mistake[2][1=]{\todo[linecolor=red,backgroundcolor=red!25,bordercolor=red,#1]{\textsc{mistake} #2}}
\newcommandx\improve[2][1=]{\todo[linecolor=orange,backgroundcolor=orange!25,bordercolor=orange,#1]{\textsc{improve} #2}}
\newcommandx\change[2][1=]{\todo[linecolor=yellow,backgroundcolor=yellow!25,bordercolor=yellow,#1]{\textsc{change} #2}}
\newcommandx\mem[2][1=]{\todo[linecolor=orange,backgroundcolor=orange!25,bordercolor=orange,#1]{\textsc{mem} #2}}
\newcommandx\status[2][1=]{\todo[linecolor=Blue,backgroundcolor=Blue!25,bordercolor=Blue,#1]{\textsc{Status} #2}}

\newtheorem*{Notationsandconventions}{Notation and conventions}
\newtheorem*{theoremnonumber}{Theorem}

\newcommand{\nidiam}{\mathsf{NonIntDiam}}
\newcommand{\N}{\ensuremath{\mathbb{N}}}
\newcommand{\R}{\ensuremath{\mathbb{R}}}

\newcommand{\Z}{\ensuremath{\mathbb{Z}}}

\newcommand{\PP}{\ensuremath{\mathbb{P}}}

\newcommand{\e}{\ensuremath{\mathrm{e}}}
\newcommand{\om}{\omega}

\newcommand{\eps}{\varepsilon}
\newcommand{\E}{\mathbb{E}}
\newcommand{\dif}{\mathrm{d}}

\newcommand{\1}{\mathds{1}}

\newcommand{\oz}{Ornstein-Zernike \;}

\newcommand{\Vol}{\text{Vol}}

\newcommand{\norme}[1]{\left\Vert #1\right\Vert_2}
\newcommand{\goes}[3]{\xrightarrow[#2]{#1} {#3}}

\newcommand{\EXT}{\mathsf{EXT}}
\newcommand{\Ext}{\mathsf{Ext}}

\newcommand{\Diam}{\mathsf{Diam}}

\newcommand{\Gap}{\mathsf{Gap}}
\newcommand{\BigDiam}{\mathsf{BigDiam}}

\newcommand{\bw}{\mathsf{BW}}

\newcommand{\nonint}{\mathsf{NI}}
\newcommand{\con}{\mathsf{Con}}

\newcommand{\hit}{\mathsf{Hit}}

\newcommand{\RW}{\mathsf{RW}}

\def\restriction#1#2{\mathchoice
              {\setbox1\hbox{${\displaystyle #1}_{\scriptstyle #2}$}
              \restrictionaux{#1}{#2}}
              {\setbox1\hbox{${\textstyle #1}_{\scriptstyle #2}$}
              \restrictionaux{#1}{#2}}
              {\setbox1\hbox{${\scriptstyle #1}_{\scriptscriptstyle #2}$}
              \restrictionaux{#1}{#2}}
              {\setbox1\hbox{${\scriptscriptstyle #1}_{\scriptscriptstyle #2}$}
              \restrictionaux{#1}{#2}}}
\def\restrictionaux#1#2{{#1\,\smash{\vrule height .8\ht1 depth .85\dp1}}_{\,#2}}

\mathtoolsset{showonlyrefs}
\usepackage{enumerate}

\makeatletter
\@namedef{subjclassname@2020}{\textup{2020} Mathematics Subject Classification}
\makeatother
\microtypesetup{expansion=false}

\title[Scaling limit of subcritical FK interfaces]{Entropic repulsion and scaling limit for a finite number of non-intersecting subcritical FK interfaces}

\author{Lucas D'Alimonte}
\address{University of Fribourg, Switzerland}
\email{lucas.dalimonte@unifr.ch}

\date{April 16, 2024}

\begin{document}

\maketitle

\begin{abstract}
This article is devoted to the study of a finite system of long clusters of subcritical 2-dimensional FK-percolation with $q\geq 1$, conditioned on mutual avoidance. We show that the diffusive scaling limit of such a system is given by a system of Brownian bridges conditioned not to intersect: the so-called \emph{Brownian watermelon}. Moreover, we give an estimate of the probability that two sets of $r$ points at distance $n$ of each other are connected by distinct clusters. As a byproduct, we obtain the asymptotics of the probability of the occurrence of a large finite cluster in a supercritical random-cluster model.
\end{abstract}

\setcounter{tocdepth}{2}

\makeatletter
\def\l@section{\@tocline{1}{0pt}{1em}{}{\bfseries}}
\def\l@subsection{\@tocline{2}{0pt}{4em}{}{}}
\makeatother

\tableofcontents

\section{Introduction}

Rigorous understanding of the behaviour of \textit{interfaces} in statistical mechanics models has been the focus of intensive study for more than 50 years, especially in the case of the Ising model. The first rigorous results were perturbative and made use of the Pirogov--Sinaï theory to show that a low temperature two-dimensional Ising interface converges, after an appropriate diffusive scaling, towards a Brownian bridge \cite{Gallavotti, higuchi}. However, these works are restricted to the very low temperature regime, even if the belief was that the result should hold for any subcritical temperature.

In the beginning of the XXI$^\text{st}$ century, the development and the understanding of the rigorous \textit{Ornstein--Zernike theory}, first in Bernoulli percolation and later on in the context of more dependent models such as the Ising and Potts models  \cite{ornsteinzernikebernoulli2002, ozisingcampanino,ozrandomclustercampanino}, provided a new powerful tool for a detailed study of the subcritical phase of these percolation or spin models. The structural output of this theory is the probabilistic description of long clusters (or equivalently of long interfaces as we shall see below) in terms of one-dimensional ``irreducible pieces'' behaving almost independently  (for a precise statement, we refer to Theorem~\ref{oz theorem}). In particular, the diffusive scaling limit of interfaces at any subcritical temperature could be obtained in the case of the Ising model as a quite simple byproduct of this robust theory in the work of Greenberg and Ioffe \cite{greenbergioffe} (see \cite{kovchegov} for the simpler case of Bernoulli percolation). 

Later on, this technique has been found to be efficient for studying interfaces interacting with their environment. Indeed, the above mentioned works deal with unconstrained (also called \emph{free}) interfaces, but recent works have been extending the study of these interfaces to broader settings in which non-trivial interactions with the environment are added.  Let us cite \cite{Ott_2018} for the case of a defect line in the Potts model, and --- much more related to this work --- \cite{ottvelenikwachtelioffe} for the treatment of a Potts interface above a boundary wall. These examples of interfaces interacting with their close environment have turned out to be more delicate to handle and in certain conditions have been shown to exhibit highly non-trivial behaviours such as \textit{wetting transitions}, which have been studied in \cite{IV18}. 

Of the same nature is the study of a system of multiple interacting interfaces, which is the focus of the present work. Indeed, this paper determines the scaling behaviour of a finite number of long clusters of subcritical Fortuin--Kasteleyn (FK) percolation, conditioned not to intersect; subcritical percolation clusters mimic interfaces in the low temperature regime. 

An interesting feature of this setting is that when conditioned not to intersect, the interaction between the clusters can turn out to be \emph{attractive}, \emph{a priori} allowing the existence of a \emph{pinning} transition --- a regime where this attraction is so strong that the clusters actually remain at a bounded distance from each other. 
We rule out the existence of such a transition. In the fashion of \cite{ottvelenikwachtelioffe}, we show that the behaviour of this system obeys an \emph{entropic repulsion} phenomenon: the entropy caused by the large number of possible clusters wins over the energetic reward obtained by keeping them  close together, all the way up to the critical point. Such a phenomenon has been previously identified in a variety of settings, for instance in the three-dimensional semi-infinite Ising model at low temperatures \cite{frohlichpfisterSemiInfiniteIsing}, the 2+1-dimensional SOS model above a hard wall \cite{caputomartinellitoninelliEntropicRepulsionSOS}, a 1+1-dimensional interface above an attractive field in presence of a magnetic field \cite{Ve2004} or a supercritical Potts interface above a wall \cite{ottvelenikwachtelioffe}, to mention but a few works studying this phenomenon. 

In this work, entropic repulsion of the FK clusters at any subcritical temperature is established in Proposition~\ref{proposition global entropic repulsion}, which is probably the most important output of this work. As a byproduct, we derive two results regarding the global behaviour of such a system of conditioned clusters. The first one is the diffusive scaling limit of such a system, which is shown to be a system of Brownian bridges conditioned not to intersect: the so-called \emph{Brownian watermelon}. Moreover, we observe that the entropic repulsion phenomenon also allows the computation --- up to a multiplicative constant --- of the probability of the existence of such a system of clusters. Finally, as a byproduct of the latter observation and using the duality property of the planar random-cluster model, we also obtain the asymptotic behaviour of the probability of the occurrence of a large finite connection in the supercritical random-cluster model.  

We mention that this work deals with the analysis of the asymptotic behaviour of the probability of \emph{pairwise} connections in the subcritical random-cluster model, while results about triple connections of points in a subcritical percolation cluster by disjoint paths in dimension $d \geq 2$ already appeared in~\cite{CampaninoGianfelicetripleconnections}. This result is however restricted to Bernoulli percolation (the case $q=1$ in the random-cluster model). Finally, we close this bibliographical introduction by mentioning that perturbative results have been obtained in~\cite{BragaProcacciSanchisPerturbativeOZsupercriticalBernoulliperco, CampaninoGianfeliceSupercriticalOZforperturbativeRC} establishing the Ornstein-Zernike behaviour of truncated point-to-point correlation functions for supercritical random-cluster models in dimensions $d \geq 3$ in the regime in which $p$ is very close to 1. In contrast, our estimate in dimension 2 shows an \emph{anomalous} prefactor (an exponent 2 rather than 1/2) to the exponential decay of the truncated correlation function. This is a special instance of the planarity of the model, leading to the entropic repulsion phenomenon, first observed in~\cite{louidor} in the case of Bernoulli percolation.

The method is in spirit close to that of \cite{ottvelenikwachtelioffe}, but with considerable additional difficulties. These are essentially due to the fact that the interaction is not only between a random object and a deterministic one, but between several random objects, forcing one to control their \emph{joint} behaviour. The proofs make heavy use of the Ornstein--Zernike theory for subcritical random-cluster models, developed in \cite{ozrandomclustercampanino}.

\subsection{Definitions of the random-cluster model and the Brownian watermelon}
\subsubsection{The random-cluster model}
The statistical mechanics model of interest is the so-called \textit{random-cluster model} (also known as \textit{FK-percolation}). We first recall its definition and a few basic properties (we refer to~\cite{duminilcopin2017lectures} for a complete exposition).
The random-cluster model on $\Z^2$ is a model of random subgraphs of $\Z^2$. Its law is described by two parameters, $p \in [0,1]$ and $q > 0$. 

Let $G = \left( V(G), E(G) \right)$ be a finite subgraph of $\Z^2$. We denote its \textit{inner boundary} (resp. \emph{outer boundary}) by 
\begin{align*}
\partial G &= \left\lbrace x \in V(G), \exists y \notin V(G), \lbrace x,y \rbrace \in E(\Z^2) \right\rbrace \text{ and}
\\ 
\partial_{\mathsf{ext}} G &= \left\{ x \notin V(G), \exists y \in V(G), \lbrace x,y \rbrace \in E(\Z^2) \right\}, \text{ respectively.}
\end{align*}
A \textit{percolation configuration} on $G$ is an element $\om \in \left\lbrace 0,1 \right\rbrace^{E(G)}$. We say that an edge $e \in G$ is \textit{open} if $\om(e) = 1$ and \textit{closed} otherwise. Two vertices $x,y \in \Z^2$ are said to be connected if there exists a path of nearest neighbour vertices $x=x_0, x_1, \dots, x_n = y$ such that the edges $\lbrace x_i, x_{i+1} \rbrace$ are open for every $0 \leq i \leq n-1$. In this case, we say that the event $\lbrace x \leftrightarrow y \rbrace$ occurs.  A \textit{vertex cluster} of $\om$ is a maximal connected component of the set of vertices (it can be an isolated vertex). Given a percolation configuration $\om$, we denote by $o(\om)$ its number of open edges, and by $k(\om)$ its number of vertex clusters.

A \textit{boundary condition} on $G$ is a partition $\eta = P_1 \cup \dots \cup P_k$ of $\partial G$. From a configuration $\om \in \left\lbrace 0,1 \right\rbrace^{E(G)}$, we create a configuration $\om^\eta$ by identifying the vertices that belong to the same $P_i$ of $\eta$. Two particular boundary conditions, that we shall call the \textit{free boundary} condition (resp. wired boundary condition), consist in the partition made of singletons (resp. of the whole set $\partial G$). We shall write $\eta= 0$ (resp. $\eta=1$) for this specific boundary condition.
\begin{definition}
Let $G =  \left( V(G), E(G) \right)$ be a finite subgraph of $\Z^2$, and $\eta$ be a boundary condition on $G$. Let $p \in [0,1]$ and $q > 0$. The random-cluster measure on $G$ with boundary condition $\eta$ is the following probability measure on percolation configurations on $G$:
\begin{equation}
\phi^\eta_{p,q,G}\left( \om \right) = \frac{1}{Z^\eta_{p,q,G}} \left(\frac{p}{1-p}\right)^{o(\om)}q^{k(\om^\eta)},
\end{equation}
where $Z^\eta_{p,q,G}>0$ is the normalisation constant ensuring that $\phi^\eta_{p,q,G}$ is indeed a probability measure. We shall refer to $Z^\eta_{p,q,G}$ as the \textit{partition function} of the model.
\end{definition}

It is classical that for $\eta = 0$ and $\eta = 1$, weak limits of the measures $\phi^\eta_{p,q,G_n}$ can be taken over any exhaustion $\left(G_n\right)_{n \in \N}$ of $\Z^2$, and that the limit measure does not depend of the choice of the exhaustion. Below, we will write $\phi_{p,q}^\eta$ for those two measures on $\Z^2$.

A very fundamental feature of this model is that it undergoes a \textit{phase transition}. Namely for any $q \geq 1$, there exists a critical parameter $p_c = p_c(q) \in (0,1)$ such that:

\begin{itemize}
    \item $\forall p < p_c(q), \phi^1_{p,q}\left( 0 \leftrightarrow \infty \right) = 0$;
    \item $\forall p > p_c(q), \phi^0_{p,q} \left( 0 \leftrightarrow \infty \right) > 0$,
\end{itemize}
where $\lbrace 0 \leftrightarrow \infty \rbrace$ is the event that the cluster of 0 is infinite.

We are going to be interested in the first case --- called \textit{the subcritical regime}. In this case it is well known that the choice of boundary conditions does not affect the infinite volume measure. We thus drop $\eta$ from the notation and simply write $\phi_{p,q}$ for the unique infinite volume measure when $p<p_c(q)$. Another important feature of the subcritical random-cluster model is the existence and the positivity of the following limit:
\begin{equation}\label{Def tau}
    \tau_{p,q} := \lim_{n \rightarrow \infty} -\frac{1}{n}\log\left[\phi_{p,q}\left(0 \leftrightarrow (n, 0) \right)\right].
\end{equation}
We call this quantity the \textit{inverse correlation length} in the direction $\Vec{e_1}$. Moreover, standard subadditivity arguments yield that 
\begin{equation}\label{eq:tau}
    \forall  x \in \Z^2, \phi_{p,q} \left[ x \leftrightarrow x + (n,0)  \right] \leq \e^{-\tau_{p,q} n}.
\end{equation}
Since $p,q$ will be fixed through this work, we shall simply write $\tau>0$ instead of $\tau_{p,q}$.

\subsubsection{The Brownian watermelon}\label{subsection bw}
The Brownian watermelon is a stochastic process that arises in various areas of probability theory, like random matrix theory~\cite{baik2007}, integrable probability~\cite{johansonn2004}, but also more recently in the study of the KPZ universality class~\cite{hammond2016brownian}.

We give a brief definition of this object, and we refer to~\cite{oconnellyor},~\cite{GRABINER1999177} and~\cite{conditionallimittheoremsfororderedrandomwalks} for the full construction and details. Let $r \geq 1$ be an integer. We define the Weyl chamber of order $r$:
\begin{equation}
    W = \left\lbrace (x_1, \dots, x_r) \in \R^r, x_1 < \dots< x_r \right\rbrace.
\end{equation}
We shall also introduce the functional Weyl chamber in the interval $[s,t]$ for $0\leq s<t$ (for any set $A \subset \R$, the set $\mathcal{C}(A, \R^r)$ denotes the space of continuous functions from $A$ to $\R^r$):
\begin{equation}
    \mathcal{W}_{[s,t]} = \lbrace f \in \mathcal{C}([s,t], \R^r), \forall s \leq \ell \leq t, f(\ell) \in W  \rbrace.
\end{equation}
Moreover let $\Delta$ denote the Vandermonde function, defined for any $(x_1, \dots, x_r) \in \R^r$ by:
\begin{equation}
    \Delta(x_1, \dots, x_r) = \prod_{1 \leq i < j \leq r}(x_j-x_i). 
\end{equation}
\begin{definition}[Brownian watermelon]
The Brownian watermelon with $r$ bridges is the continuous process $\big(\bw^{(r)}_t\big)_{0\leq t \leq 1}$ obtained by conditioning $r$ independent standard Brownian bridges not to intersect in $(0,1)$. It is a random object of $\mathcal{C}([0,1], \R^r)$.
\end{definition}
\begin{remark}
Since the non-intersection event has null probability for $r$ random bridges as soon as $r \geq 2$, the latter conditioning is rigorously done by means of a Doob $h$-transformation by the harmonic function $\Delta$. We refer to~\cite{oconnellyor} and~\cite{Invarianceprinciplesforrandomwalksincones} for the details of the construction (and the fact that $\Delta$ is harmonic for a system of $r$ standard bridges). Moreover, it can be shown, by means of the Karlin--McGregor formula, that for any $0<t<1$
\begin{equation}\label{equation marginale BW}
    \PP\left[ \bw^{(r)}_t \in \mathrm{d}z\right] \propto \frac{1}{\left(t(1-t)\right)^{r^2/2}}\Delta^2(z)\e^{-\frac{|z|^2}{2t(1-t)}}\1_{z \in W}\dif z.
\end{equation}
\end{remark}
\begin{remark}
Alternatively, the Brownian watermelon can be built \emph{via} the following method: consider a system $(B^\eps_t)_{0 \leq t \leq 1}$ of $r$ independent standard Brownian bridges started from 0, $\eps, \dots, (r-1)\eps$ respectively. Then under the conditioning on the event $\left\lbrace B^\eps_t \in \mathcal{W}_{[0,1]} \right\rbrace$ (this happens with positive probability), the following weak limit exists in $\mathcal{C}([0,1], \R^r)$ when $\eps \rightarrow 0$ and is called the Brownian watermelon:
\begin{equation}
   ( B^\eps_t )_{0\leq t \leq 1} \goes{(d)}{\eps \rightarrow 0}  {(\bw^{(r)}_t)_{0\leq t\leq 1}}.
\end{equation}
For more information on this construction, see~\cite{oconnellyor} and~\cite{hammond2016brownian}.
\end{remark}

\begin{Notationsandconventions} If $a_n$ and $b_n$ are two sequences of real numbers, we shall write $a_n \sim b_n$ when $\frac{a_n}{b_n} \goes{}{n \rightarrow \infty}{1}$. We shall also write $a_n = o(b_n)$ when $\frac{a_n}{b_n} \goes{}{n \rightarrow \infty}{0}$ and $a_n = O(b_n)$ when there exists a constant $C> 0$ such that $\vert a_n \vert \leq C\vert b_n \vert$ for all $n \geq 0$. Moreover, we shall write $a_n \asymp b_n$ whenever $a_n = O(b_n)$ and $b_n = O(a_n)$. Finally, the generic notation $c,C>0$ will denote constants depending only on $p$ and $q$, that may change from line to line during computations. We denote by $\Vert \cdot \Vert$ the Euclidian norm on $\R^d$. 
\end{Notationsandconventions}

\subsection{Exposition of the results}

In this paper, we study the scaling limit of a system of subcritical clusters conditioned on a connection and a non-intersection event. We first start by defining these percolation events.

\begin{definition}[Connection event, Non-intersection event]
Let $x, y \in W\cap\Z^r$ and $n\geq 0$.
Then we define the multiple connection event $\con^n_{x,y}$ by
\begin{equation}
    \con^n_{x,y} = \left\lbrace \forall 1 \leq i \leq r, (0,x_i) \leftrightarrow (n,y_i) \right\rbrace.
\end{equation}
The non-intersection event will be defined by
\begin{equation}
    \nonint_{x} = \left\lbrace \forall 1 \leq i < j \leq r, ~\mathcal{C}_{(0,x_i)} \cap \mathcal{C}_{(0, x_j)} = \emptyset \right\rbrace,
\end{equation}
where $\mathcal{C}_u$ denotes the cluster of the vertex $u \in \mathbb Z^2$.
In the rest of this work, as $n, x, y$ will be fixed, we shall abbreviate $\con^n_{x,y}$ by $\con$ and $\nonint_x$ by $\nonint$. Moreover, we will also abbreviate $\mathcal{C}_{(0,x_i)}$ by $\mathcal{C}_i$. 
\end{definition}

Our main result consists in the estimation of the probability that $\lbrace \con, \nonint \rbrace := \lbrace \con \cap \nonint \rbrace  $ occurs in a subcritical random-cluster measure. 
\begin{theorem}\label{theoreme estimation}
Let $q\geq 1$, and $0<p<p_c(q)$. Let $r\geq 1$ be a fixed integer. Then, there exist two constants $C_-, C_+ > 0$ such that for any sequences $x(n), y(n)$ of elements of $W$ satisfying $\Vert x(n) \Vert, \Vert y(n) \Vert = o(\sqrt{n})$, when $n$ is sufficiently large,
\begin{equation} 
  C_-V(x(n))V(y(n))n^{-\frac{r^2}{2}}\e^{-\tau rn} \leq \phi\left[\con, \nonint \right] \leq C_+ V(x(n))V(y(n))n^{-\frac{r^2}{2}}\e^{-\tau rn},
\end{equation}
where $V$ is the function defined in Theorem~\ref{theoreme local limit srw}.
\end{theorem}
\begin{remark}
The function $V$ is not explicit. However, it is known that (see Theorem~\ref{theoreme local limit srw}) :
\begin{equation} \text{When }\min_{1\leq i \leq r-1} \lbrace |x(n)_{i+1}- x(n)_i | \rbrace \goes{}{n\rightarrow\infty}{+\infty}, \text{ then } V(x(n))\sim\Delta(x(n)) \text{ as } n\rightarrow \infty. \end{equation} 
Moreover, when $x,y$ are fixed elements of $W$, the statement simplifies as 
\begin{equation} 
  \phi\left[\con, \nonint \right] \asymp n^{-\frac{r^2}{2}}\e^{-\tau rn}.
\end{equation}
\end{remark}
An interesting corollary, which is a direct consequence of Theorem~\ref{theoreme estimation} in the case $r=2$, can be obtained using the methods of~\cite{louidor}, where the same result is proved in the case of Bernoulli percolation (corresponding to $q=1$ in the random-cluster model). Let us define the \textit{truncated inverse correlation length} in the direction $\Vec{e_1}$ by 
\begin{equation}
    \tau^\mathsf{f}_p = \lim_{n \rightarrow \infty} -\frac{1}{n}\log\phi\left[ 0 \leftrightarrow (n, 0), \left| \mathcal{C}_0 \right| < \infty\right].
\end{equation}
It is well known that on $\Z^2$, whenever $p \neq p_c(q)$, one has that $\tau^\mathsf{f}_p > 0$. Moreover, it is clear that whenever $p < p_c(q)$, $\tau^\mathsf{f}_p = \tau_p$, where $\tau_p$ has been defined in~\eqref{Def tau}.
Then, Theorem~\ref{theoreme estimation} allows to compute the prefactor in the supercritical truncated correlation function. 
\begin{corollary}
Let $q\geq 1$ and $p\in (p_c, 1)$. Let $\phi$ be the unique infinite-volume random-cluster measure on $\Z^2$. Then, 
\begin{equation}
    \phi\left[ 0 \leftrightarrow (n,0), \left| \mathcal{C}_0 \right| < \infty \right] \asymp \frac{1}{n^2}\e^{-2\tau^{\mathsf{f}}_{p^*}n},
\end{equation}
where $p^*$ stands for the dual parameter of $p$ (see~\eqref{eqdualite} for the relation linking $p$ and $p^*$).
\end{corollary}
\begin{remark}
In particular, we obtain the following equality, holding for any supercritical $p>p_c$
\begin{equation}
    \tau^\mathsf{f}_{p} = 2\tau^\mathsf{f}_{p^*} (= 2\tau_{p^*}).
\end{equation}
This is a very specific instance of duality, and such a relation is not expected to hold in higher dimensions. The result was already well known in the case of Bernoulli percolation, see for instance~\cite{correlationlengthonehalf} or~\cite[Theorem 11.24]{grimmett}.
\end{remark}

\begin{figure}
    \centering
    \includegraphics[width = .55\textwidth, page = 5]{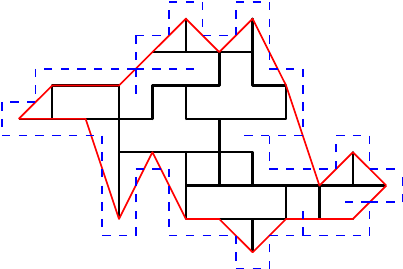}
    \caption{Depiction (in red) of the envelopes $\Gamma^+$ and $\Gamma^-$ of a percolation cluster. The blue dashed path corresponds to the more natural notion of interface that could have been considered instead. However, as explained in Remark~\ref{remarque interface}, this blue path converges in the space of continuous curves towards the Brownian watermelon as well.}
    \label{figure illustration interfaces}
\end{figure}

Our second result consists in the study of the behaviour of the $r$ clusters created by conditioning on $\lbrace \con, \nonint \rbrace$. It will be formulated in terms of the \textit{envelopes} of a cluster.
\begin{definition}[Upper and lower envelopes of a cluster]\label{definition interfaces}
Let $\omega \in \con$. Then for any $0\leq k \leq n$ and $1 \leq i\leq r$ we define (see Figure~\ref{figure illustration interfaces})
\begin{equation}
    \Gamma^+_i(k) = \max\left\lbrace \ell \in \Z, (k,\ell) \in \mathcal{C}_i \right\rbrace \text{ and } \Gamma^-_i(k) = \min\left\lbrace \ell\in \Z, (k,\ell) \in \mathcal{C}_i \right\rbrace.
\end{equation}
It is clear that $\Gamma_i^\pm$ are well defined, since all clusters are almost surely finite in the subcritical regime, and the sets above are not empty due to $\con$. We will see these quantities as functions from $[0,n]$ to $\R$ by considering the piecewise affine functions $\Gamma^\pm_i(t)$ that coincide with $\Gamma^\pm_i$ on the integers $t=k$. 
\end{definition}
Our second result is the following:
\begin{theorem}\label{Theoreme main}
Fix $x,y \in W\cap\Z^r$ and $p\in(0,p_c(q))$. Then under the family of measures $\phi_{p,q}\left[ ~\cdot \vert \con,\nonint\right]$ (we recall that $\con, \nonint$ depend on $n$), there exists $\sigma > 0$ such that:
\begin{equation}\label{equation convergence main}
    \left( \frac{1}{\sqrt{n}}\left(\Gamma^+_1(nt), \dots, \Gamma^+_r(nt) \right)\right)_{0\leq t \leq 1} \goes{(d)}{n\rightarrow \infty}{ (\sigma\bw^{(r)}_t)_{0\leq t \leq 1}},
\end{equation}
where $\bw^{(r)}$ is the Brownian watermelon with $r$ bridges, and where the convergence holds in the space $\mathcal{C}\left([0,1],\R^r\right)$ endowed with the topology of uniform convergence. Moreover, almost surely, for all $1\leq i\leq r$,
\begin{equation}\label{equation shrink interface thm convergence}
\frac{1}{\sqrt{n}}\left\Vert \Gamma^+_i - \Gamma^-_i \right\Vert_{\infty} \goes{}{n \rightarrow \infty}{0}
\end{equation}
\end{theorem}
\begin{remark}
A consequence of~\eqref{equation shrink interface thm convergence} is that in the setting of Theorem~\ref{Theoreme main}, the clusters remain of width $o(\sqrt{n})$. Actually, we prove that almost surely, $\Vert \Gamma^+_i - \Gamma^-_i \Vert_{\infty} = O(\log n)$. In particular, the choice of the upper interfaces $\Gamma_i^+$ in~\eqref{equation convergence main} is arbitrary and can be replaced by any assignment of $\pm$ for the choice of interfaces to converge.
\end{remark}
\begin{remark}
The result is stated for fixed $x,y \in W\cap \Z^r$. However, the careful reader may check that our method allows to treat the case where $x$ and $y$ depend on $n$. Indeed, as soon as $x(n), y(n)$ are two sequences of $W\cap\Z^r$ satisfying
\begin{equation}
    \Vert x(n) \Vert = o(\sqrt{n}) ~~~\text{and}~~~\Vert y(n)\Vert = o(\sqrt{n}),
\end{equation}
our methods may apply and yield the same scaling limit.
\end{remark}
\begin{remark}\label{remarque interface}
For the reader familiar with statistical mechanics, it might seem strange that our result is formulated in terms of these envelopes and not in terms of the upper and lower interfaces running along the boundary of the clusters $\mathcal{C}_i$. However, it may be shown that the interfaces also converge to the paths of $\bw^{(r)}_t$ (as paths in $[0,1] \times \R^r$). We then chose to work with $\Gamma^\pm$ since we can use the space of continuous functions from $[0,n]$ to $\R^r$ for studying convergence questions, which is easier to treat than the space of continuous curves which would be needed when considering those interfaces. 
\end{remark}

\subsection{Background on the random-cluster model}

We first recall some basic properties of the random-cluster model (once again we refer to~\cite{duminilcopin2017lectures} for a complete exposition). These properties are valid for any choice of parameters $p$ and $q$.
\medskip

\noindent{\bf Positive association.}
The space $\left\lbrace 0,1 \right\rbrace^{E(\Z^2)}$ can be equipped with a partial order: we say that $\om_1 \leq \om_2$ if for any $e \in E(\Z^2)$, $\om_1(e) \leq \om_2(e)$. An event $\mathcal{A}$ will be called \textit{increasing} if for any $\om_1 \leq \om_2$, $\om_1 \in \mathcal{A} \Rightarrow \om_2 \in \mathcal{A}$. The \textit{FKG inequality} then states that for any increasing events $\mathcal{A}, \mathcal{B}$, any graph $G$ and any boundary conditions $\eta$, 
\begin{equation}\tag{FKG}\label{equation FKG}
    \phi_{G,p,q}^\eta[\mathcal{A} \cap \mathcal{B}] \geq \phi^\eta_{G,p,q}[\mathcal{A}]\phi^\eta_{G,p,q}[\mathcal{B}].
\end{equation}
This property implies in particular that for any boundary conditions $\eta_1 \leq \eta_2$ (meaning that the partition $\eta_1$ is finer than $\eta_2$), for any increasing event $\mathcal{A}$,
\begin{equation}\tag{CBC}\label{Comparaison boundary conditions}
    \phi_{G,p,q}^{\eta_1}\left[\mathcal{A}\right] \leq \phi_{G,p,q}^{\eta_2}\left[\mathcal{A}\right].
\end{equation}
This property is called the \textit{comparison of boundary conditions} and may also be stated as ``$\phi^{\eta_1}$ is stochastically dominated by $\phi^{\eta_2}$".
\medskip

\noindent{\bf Duality.}
Let $(\Z^2)^* = (\frac{1}{2},\frac{1}{2}) + \Z^2$ and consider the lattice $(\Z^2)^*$ with edges between nearest neighbours. This lattice is called the \textit{dual lattice}. It has the property that for any $e \in E(\Z^2)$, there exists a unique edge $e^* \in E((\Z^2)^*)$ that crosses $e$. To a percolation configuration $\om \in \left\lbrace 0,1 \right\rbrace^{E(\Z^2)}$ we can associate a dual configuration $\om^*$ on the dual lattice by setting $\om^*(e^*) = 1 - \om(e)$. Then we remark that --- as soon as the parameters guarantee that there exists a unique Gibbs measure --- if $\om$ is sampled according to $\phi_{p,q}$, then $\om^*$ has the distribution of $\phi_{p^*,q^*}$, where 
\begin{equation}\label{eqdualite}
    q=q^* \text{ and } \frac{pp^*}{(1-p)(1-p^*)} = q.
\end{equation}
It has been proved by V. Beffara and H. Duminil-Copin in~\cite{duminilbeffara} that $p_c(q) = p^*_c(q)$, meaning that the parameter $p_c(q)$ is \textit{self-dual}. Also observe that if $\phi_{p,q}$ is subcritical, then $\phi_{p^*,q^*}$ is supercritical and vice-versa. 
\medskip

\noindent{\bf Spatial Markov property.} Let $G$ be a subgraph of $\Z^2$, and $G' \subset G$ a subgraph of $G$. Let $\xi$ be a percolation configuration on $\Z^2$. Observe that it induces a boundary condition on $G$ --- that we name $\eta(\xi)$ --- by identifying the vertices wired together by $\xi$ outside $G$, and a boundary condition on $G'$ - that we name $\eta'(\xi)$ by the same principle. Then,
\begin{equation}\tag{SMP}\label{equation smp}
    \phi_{G,p,q}^{\eta(\xi)}\left[ ~\cdot~\vert \om(e) = \xi(e), ~\forall e \notin G' \right] = \phi^{\eta'(\xi)}_{G',p,q}[\cdot].
\end{equation}

\noindent{\bf Finite energy property.} When $p \notin \lbrace 0,1 \rbrace$, there exists a constant $\eps>0$ depending only on $p$ and $q$ such that for any finite graph $G$, any finite $F \subset E(G)$, any boundary condition $\eta$, and any percolation configuration $\om_0$,
\begin{equation}
    \eps^{|F|} \leq  \phi_{G,p,q}^\eta\left[\om(e) = \om_0(e),~ \forall e \in F \right] \leq (1-\eps)^{|F|}.
\end{equation}
\medskip

\noindent{\bf Weak ratio mixing.}
In the subcritical regime, the random-cluster measure also enjoys the following \textit{weak ratio mixing property}. For two finite connected sets of edges $E_1$ and $E_2$, define their distance $d(E_1, E_2)$ as the Euclidean distance between the set of their respective endpoints. Then, for any graph $G$, any boundary condition $\eta$, any $q\geq 1$ and any $p<p_c(q)$, there exists a constant $c>0$ such that for any events $\mathcal{A}$ and $\mathcal{B}$ depending on edges of $E_1$ and $E_2$ respectively,
\begin{equation}\label{weak ratio mixing}\tag{MIX}
    \left| 1- \frac{\phi_{G,p,q}^\eta[\mathcal{A}\cap\mathcal{B}]}{\phi_{G,p,q}^\eta[\mathcal{A}]\phi_{G,p,q}^\eta[\mathcal{B}]} \right| < \e^{-cd(E_1,E_2)}.
\end{equation}
\medskip

\subsection{Outline of the proof}

The main idea of modern Ornstein--Zernike theory is to couple a subcritical percolation cluster conditioned on realizing a connection event $\left\lbrace x \leftrightarrow y \right\rbrace$ with a random walk started from $x$ and conditioned to reach $y$. Such a cluster is essentially a one-dimensional object. As the knowledge on conditioned random walks is very broad, in particular in terms of Local Limit Theorems and invariance principles, such a coupling allows to derive properties of the original cluster. In our setting, we would like to couple a system of $r$ percolation clusters conditioned on $\con \cap \nonint$ with a system of $r$ random walks conditioned on a hitting event and on not intersecting each other. However, such a coupling is not immediately available in this setting and we have to rely on several comparison principles to show that the behaviours of these two types of systems are close. Once this task is accomplished, we use an invariance principle for a system of non-intersecting random walks to derive Theorem~\ref{Theoreme main}. 

 Let us be a bit more precise about the method. We first show that Ornstein--Zernike theory extends to $r$ non-intersecting clusters sampled according to $\phi^{\otimes r}$ (the product of $r$ random cluster measures on $\Z^2$) and thus interacting only through the conditioning. This allows us to derive an invariance principle for this product measure. 
 
 The next step is to transmit the results obtained for the product measure to the ``true" FK-percolation measure. As crucially observed in~\cite{ottvelenikwachtelioffe}, this can be done proving an \textit{a priori} (meaning independent of the above mentioned coupling) repulsion estimate: under the conditioned random-cluster measure, the clusters naturally move far from each other and never come near each other again. This input will then allow us to use the mixing property of subcritical FK-percolation to derive Theorems~\ref{theoreme estimation} and~\ref{Theoreme main}.

\subsection{Organization of the paper}

We first focus on Theorems~\ref{theoreme estimation} and~\ref{Theoreme main}. As explained previously, the proof consists of two independent tasks: comparing the behaviours of the percolation clusters and of a system of interacting random walks, and then obtaining the scaling limit and fine estimates on such a system of random walks. Our interest mainly being statistical mechanics, we postpone all the results about interacting random walks to Section~\ref{section marches}, which is  independent of the other sections, and may be skipped by readers only interested in the percolation aspects. Section~\ref{section review oz} consists in a review of the rigorous Ornstein--Zernike results for one single subcritical cluster of FK-percolation. Section~\ref{section independent system} is devoted to the study of the scaling limit under the \textit{product measure} through a straightforward extension of the Ornstein--Zernike theory to this setting, as discussed before. Finally, Section~\ref{section RCM} is devoted to the proof of the entropic repulsion estimates, and thus of the announced result. 

\section{Ornstein--Zernike theory for a single subcritical FK-cluster}\label{section review oz}

In the remainder of the paper, we \textbf{fix} $q \geq 1$ and $0< p < p_c(q)$. Since these parameters will not change throughout the paper, we drop them from the notation and abbreviate $\phi_{p,q} := \phi$.

In this section, we review and discuss the main result of~\cite{ozrandomclustercampanino} --- the Ornstein--Zernike Theorem. Schematically, this result can be described as follows. Under the conditioned measure $\phi \left[~\cdot\vert y \in \mathcal{C}_0\right]$ where $y$ is some vertex far away from 0, the cluster of $0$ has a very particular structure. Indeed, it macroscopically looks like the geodesic from 0 to $y$. Moreover, it exhibits typical Brownian bridge fluctuations around this geodesic, and is confined in a very small tube around this Brownian bridge. The result is precisely stated in Theorem~\ref{oz theorem}.

\begin{definition}[Directed random walk]\label{def directed probability measure}
A \textit{directed probability measure} on $\Z^2$ is a probability measure on $\N^* \times \Z$. If $X_1, \dots, X_n, \dots$ are independent and identically distributed random variables sampled according to a directed probability measure on $\Z^2$, then the process $S_n = X_1+\dots +X_n$ is called a \textit{directed random walk}. We shall call a possible realization of $(S_n)$ a \textit{directed walk} on $\Z^2$.
\end{definition}

In the remainder of the paper, we will often interpret trajectories of directed walks as real-valued functions defined on $\R^+$. Indeed, let $\nu$ be a directed probability measure on $\Z^2$ and $(S_n)_{n \geq 0}$ the associated directed random walk. Since $\nu(\N^* \times \Z)=1$, for any $t\geq 0$, the trajectory of $S$ almost surely intersects the vertical line $\lbrace t \rbrace \times \R$ once. Calling this point $S(t)$ provides us with a continuous and piecewise linear function: moreover this correspondence is one-to-one. We shall often use notation as $\lbrace S \in \mathcal{A} \rbrace$, where $\mathcal{A}$ is a subset of $\mathcal{C}(\R^+, \R)$. In that case, $S$ will have to be taken as the continuous function described above. Let $(S_n)_{n \geq 0}$ be a directed random walk on $\Z^2$. If $y \in \Z$, introduce the event 
\begin{equation}
    \hit_y = \left\lbrace \exists n \geq 0, S_n = y \right\rbrace.
\end{equation}

\subsection{Diamond confinement and diamond decomposition}\label{subsubsection diamond decomposition}

We need a bit of vocabulary, in order to properly state the confinement property of a long subcritical cluster. Let $\delta > 0, x \in \Z^2$. Following~\cite{ozrandomclustercampanino}, introduce the following subsets of $\Z^2$:

\begin{itemize}
    \item The $\delta$-forward cone of apex $x$ to be the set  $\mathcal{Y}^{\delta,+}_x = x + \left\lbrace (x_1, x_2) \in \Z^2, \delta x_1 \geq |x_2|    \right\rbrace.$ 
    \item The $\delta$-backward cone of apex $x$ to be the set $\mathcal{Y}^{\delta,-}_x = x + \left\lbrace (x_1,x_2) \in \Z^2, \delta x_1 \leq -|x_2|  \right\rbrace.$ 
    \item If $x,y \in \Z^2$ are such that $x_1<y_1$, the $\delta$-diamond of apexes $x, y$ is the intersection:
\begin{equation}    
    \mathcal{D}^\delta_{x,y} = \mathcal{Y}^{\delta,+}_x \cap \mathcal{Y}^{\delta,-}_y. 
\end{equation}
If $x=0$, we abbreviate the notation by $\mathcal{D}_y^\delta$.
\end{itemize}

Let $G$ be a finite subgraph of $\Z^2$ containing the vertex $0$ (we say that $G$ is a subgraph of $\Z^2$ \textit{rooted at 0}). We say that:
\begin{itemize}
    \item $(G,v)$ is $\delta$-\textit{left-confined} if there exists $x\in V(G)$ such that $G \subset \mathcal{Y}^{\delta,-}_x$.
    \item $(G,v)$ is $\delta$-\textit{right-confined} if there exists $x\in V(G)$ such that $G \subset \mathcal{Y}^{\delta,+}_x$.
    \item $G$ is $\delta$-\textit{diamond-confined} if there exist $y \in V(G)$ such that $G \subset \mathcal{D}^\delta_{y}$. In that case, we say that $\mathcal{D}^\delta_{y}$ is the diamond containing $G$. 
\end{itemize}

Observe that in the previous definitions, if the points $x,y$ exist, they are necessarily unique. We denote the set of $\delta$-left-confined subgraphs of $\Z^2$ rooted at 0 (resp. $\delta$-right-confined subgraphs of $\Z^2$ rooted at 0, resp $\delta$-diamond-confined subgraphs of $\Z^2$ rooted at 0) by $\mathfrak{C}_L^\delta$ (resp $\mathfrak{C}_R^\delta$, resp $\mathfrak{D}^\delta$).

\begin{definition}
We now define the notion of \textit{displacement} along a left-confined, right-confined or diamond-confined subgraph of $\Z^2$. 
\begin{itemize}
    \item Let $G$ be a $\delta$-left-confined subgraph of $\Z^2$ rooted at 0. The \textit{displacement of $G$} is
    \begin{equation} 
    X^L(G) = x, 
    \end{equation}
    where $x$ is the unique vertex of $G$ such that $G \subset \mathcal{Y}_x^{-,\delta}$.
    \item Let $G$ be a $\delta$-right-confined subgraph of $\Z^2$ rooted at 0. The \textit{displacement of $G$} is 
    \begin{equation}
    X^R(G) = -x,
    \end{equation}
    where $x$ is the unique vertex of $G$ such that $G \subset \mathcal{Y}_x^{+,\delta}$.
    \item Let $G$ be a diamond-confined subgraph of $\Z^2$ rooted at 0. The \textit{displacement of $G$} is
    \begin{equation} X(G) = y, \end{equation}
    where $y$ is the only vertex of $G$ such that $G \subset \mathcal{D}^\delta_y$.
\end{itemize}
\end{definition}

In order to properly state what is a diamond decomposition of a cluster, we also need to introduce the operation of concatenation of two confined subgraphs rooted at 0. Let $G_1 \in \mathfrak{C}^\delta_L$ and $G_2 \in \mathfrak{D}^\delta$. The concatenation of $G_1$ and $G_2$, called $G_1\circ G_2$, is defined to be the subgraph
\begin{equation}
G_1\circ G_2 = G_1 \cup \left( X^L(G_1)+G_2 \right). 
\end{equation}
In the same manner, we can concatenate a $\delta$-diamond-confined rooted graph $G_1$ with a $\delta$-right-confined rooted graph by setting 
\begin{equation}
    G_1\circ G_2 = G_1\cup \left(X^R(G_2)+G_2\right).
\end{equation}
Finally observe that one can concatenate two $\delta$-diamond-confined rooted graphs by setting 
\begin{equation}
    G_1\circ G_2 = G_1\cup \left(X(G_2)+G_2\right).
\end{equation}

These definitions in hand, we can now define the \textit{diamond decompositions} of a subgraph of $\Z^2$.

\begin{definition}[Diamond decomposition of a subgraph, skeleton of a subgraph]
Let $G$ be a finite subgraph of $\Z^2$ rooted at 0. Then, to any decomposition of the type $G=G^L\circ G_1 \circ \dots \circ G_\ell \circ G^R$ with $G^L \in \mathfrak{C}_L, G^R \in \mathfrak{C}_R, G_i \in \mathfrak{D}^\delta$ for all $i \in \lbrace 1, \dots, \ell \rbrace$, can be associated the concatenation of the confining $\delta$-left cone with all the associated $\delta$-diamonds and the confining $\delta$-right cone. We call such a subset a \textit{diamond decomposition} of $G$:

\begin{equation} 
\mathcal{D}(G) = \mathcal{Y}^{\delta, -}_{X^L(G^L)} \circ \mathcal{D}^{\delta}_{X(G_1)}  \circ \dots \circ \mathcal{D}^{\delta}_{X(G_\ell)} \circ \mathcal{Y}^{\delta,+}_{X^R(G^R)}. 
\end{equation}

Let us call $x_0 = 0, x_1 = X^L(G^L)$, $x_k = X(G_{k-1})$ for $2\leq k \leq \ell+1$, and $x_{\ell+2} = X^R(G^R)$. We then define, for $0\leq n \leq \ell+2$,
\begin{equation}
    \mathcal{S}(\mathcal{D})(G)_n = \sum_{k=0}^{n} x_k.
\end{equation}
The process $\mathcal{S}(\mathcal{D})(G)_n$ is called \textit{the skeleton} of the diamond decomposition $\mathcal{D}(G).$
\end{definition}

\begin{remark}\label{remarque maximal skeleton}
    Observe that diamond decompositions of $G$ are not unique: as soon as there exists one of them with $\ell \geq 3$, merging inner diamonds allows one to create new (coarser) diamond decompositions of $G$. However, any finite subgraph rooted at 0 admits a unique \textit{maximal} diamond decomposition: we call it $\mathcal{D}^{\mathsf{max}}(G)$. The skeleton associated to this decomposition will by called $\mathcal{S}^{\mathsf{max}}(G)$ and referred to as the \textit{maximal} skeleton of $G$.
\end{remark}

\begin{remark}
    Our object of interest will be the skeleton of random diamond decompositions of subgraphs of $\Z^2$. Amongst the properties of the skeleton associated to a diamond decomposition of some rooted subgraph $G$, observe that the vertices of the skeleton of a diamond decomposition of $G$ are \textit{cone-points} of $G$, in the sense that for any $n \leq \ell+2$,
    \begin{equation}
        G \subset \mathcal{Y}_{\mathcal{S}(\mathcal{D})(G)_n}^{\delta, -} \cup \mathcal{Y}_{\mathcal{S}(\mathcal{D})(G)_n}^{\delta, +}.
    \end{equation}
    Furthermore, observe that the skeleton of a diamond decomposition of $G$ is always a finite directed walk, which motivates the terminology introduced in Definition~\ref{def directed probability measure}.
\end{remark}

\begin{remark}
The structure of the diamond decomposition is here given in the direction given by the first coordinate axis. However we see that adapting the definitions of the cones, the diamond decomposition can be defined for any direction $s \in \mathbb{S}^1$. The results of this work naturally adapt to this case, with this slight modification.
\end{remark}

\subsection{Ornstein--Zernike theory for one subcritical cluster}

We are ready to state the main result of~\cite{ozrandomclustercampanino}, which we shall refer to as the \oz Theorem. Set $\mathfrak{G}_0$ to be the set of connected subgraphs of $\Z^2$, rooted at 0.

\begin{figure}
    \centering
    \includegraphics[width = \textwidth, page = 2]{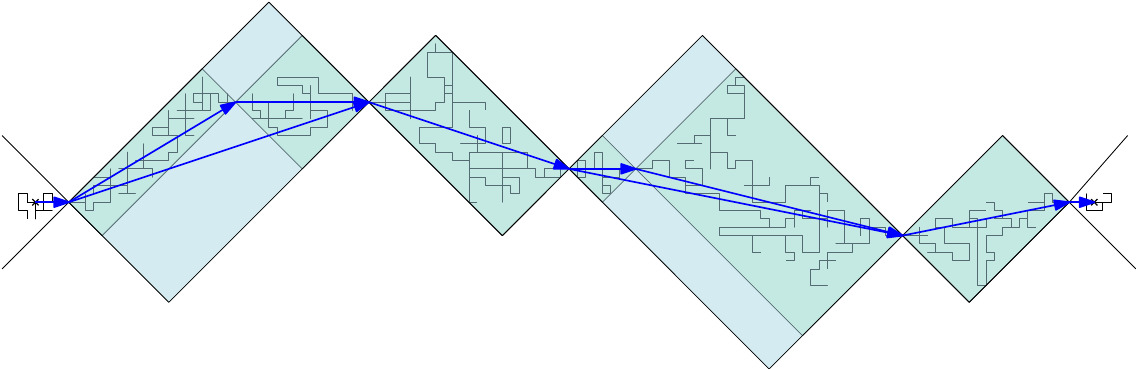}
    \caption{Two admissible diamond decompositions susceptible to appear in the Ornstein Zernike Theorem, together with their associated skeleton (in blue). Observe that the darker one is the \textit{maximal} diamond decomposition of the cluster.}
    \label{Figure ornstein zernike}
\end{figure}

\begin{theorem}[Ornstein--Zernike Theorem,~\cite{ozrandomclustercampanino}]\label{oz theorem}
There exist two constants $C,c>0$ and a positive $\delta > 0$, such that the following holds. There exist two positive finite measures $\rho_L, \rho_R$ on $\mathfrak{C}_L^\delta$ and $\mathfrak{C}_R^\delta$ respectively, and a probability measure $\mathbf{P}$ on $\mathfrak{D}^\delta$ such that for any bounded function $f: \mathfrak{G}_0\rightarrow \R$, any $y  \in \mathcal{Y}_{0}^+$,

\begin{equation}\label{equation OZ 1 cluster}
  \Big\vert \e^{\tau x_1}\phi\left[f( \mathcal{C}_0) \1_{y \in \mathcal{C}_0} \right] - \sum_{\substack{\ell \geq 0 \\ G^L\in \mathfrak{C}_L^\delta \\ G^R \in \mathfrak{C}_L^\delta \\ G_1, \dots, G_\ell \in \mathfrak{D}^\delta }}\rho_L(G^L)\rho_R(G^R)\mathbf{P}(G_1) \cdots \mathbf{P}(G_\ell) f(G) \Big\vert \leq C\Vert f\Vert_{\infty}\e^{-c\norme{y}},\end{equation}
 where the sum runs over all $G^L \in \mathfrak{C}_L^\delta, G^R \in \mathfrak{C}_R^\delta, G_1, \dots, G_\ell \in \mathfrak{D}^\delta$ satisfying the relation
 \begin{equation}
     X^L(G^L) + X(G_1)+ \dots + X(G_\ell) + X^R(G^R) = y.
 \end{equation}
 We also have written $G = G^L \circ G_1 \circ \dots \circ G_\ell \circ G^R$ in the argument of $f$. Moreover, the measures $\rho_L, \rho_R, \mathbf{P}$ have exponential tails with respect to the length of the displacement: there exist $c', C' > 0$ such that
 \begin{equation}
     \max \left\lbrace \rho_L \left[\norme{X^L(G^L)} > t  \right], \rho_R \left[ \norme{X^R(G^R)} > t  \right], \mathbf{P}\left[ \norme{X(G)} > t \right]  \right\rbrace < C'\e^{-c't}.
 \end{equation} 
\end{theorem}

In the remainder of the paper, we \textbf{fix} $\delta$ to be equal to the value given by Theorem~\ref{oz theorem}. In particular, we will not highlight the dependency anymore and we drop it from the notation.

\begin{remark}\label{Remarque definition loi nu}
For any $x \in \N^* \times \Z$, define the following three quantities:
\begin{itemize}
\item $\nu_L(x) = \sum_{G^L\in \mathfrak{C}_L, X^L(G^L)=x}\rho_L(G^L)$, 
\item $\nu_R(x) = \sum_{G^R\in \mathfrak{C}_R, X^R(G^R)=x}\rho_R(G^R),$ 
\item $\nu(x) = \sum_{G \in \mathfrak{D}, X(G) = x} \mathbf{P}\left(G\right)$.
\end{itemize}
Then, it is clear that $\nu$ is a directed probability measure on $\Z^2$, which has exponential tails. We define $\PP^\RW$ to be the directed random walk measure associated to $\nu$.
\end{remark}

These definitions allow us to formulate a second version of Theorem~\ref{oz theorem} in terms of a coupling between a percolation cluster conditioned to contain a distant point and a directed random bridge.

\begin{theorem}[Ornstein--Zernike Theorem; coupling version]\label{theoreme couplage oz 1 cluster}
Let $y \in \mathcal{Y}^+_0$. There exists a probability space $(\Omega, \mathcal{F}, \Phi_{0 \rightarrow y})$ supporting a random element $(\mathcal{C}_0, \mathcal{S})$ such that:
\begin{itemize}
    \item $\mathcal{C}_0$ has the distribution of the cluster of 0 under the measure $\phi\left[~\cdot\vert y \in \mathcal{C}_0   \right]$, ie if $C$ is a connected subgraph of $\Z^2$ containing 0,
    \begin{equation}
        \Phi_{0 \rightarrow y}\left[ \mathcal{C} = C \right] = \phi\left[\mathcal{C}_0 = C \vert y \in \mathcal{C}_0\right]
    \end{equation}
    \item $\mathcal{S}$ has the distribution of a directed random walk conditioned to hit $y$, \textit{ie} for any $\ell \geq 1$, any family $s_1, \dots, s_\ell$ of vertices of $\Z^2$,
    \begin{equation}
        \Phi_{0 \rightarrow y}\left[\mathcal{S}_1 =s_1, \dots, \mathcal{S}_\ell = s_\ell \right] \propto \nu_L(s_1)\nu_R(y-s_\ell)\prod_{k=2}^{\ell}\nu(s_k - s_{k-1}),
    \end{equation}
    where the symbol $\propto$ means that one has to normalise the latter quantity to get a proper probability measure,
    \item With probability at least $1-C\e^{-c\Vert y \Vert}$, for all $1 \leq k \leq \ell$, $\mathcal{S}_k \in \mathcal{C}_0$ and $\mathcal{S}_k$ is a renewal of $\mathcal{C}_0$. Furthermore, for any $1\leq k\leq \ell-1$, the portion of $\mathcal{C}_0$ lying between $\mathcal{S}_k$ and $\mathcal{S}_{k+1}$ (by that we mean $\mathcal{C}_0 \cap \{(x,y)\in\Z^2, (\mathcal{S}_k)_1 \leq x \leq (\mathcal{S}_{k+1})_1  \}$) is a $\delta$-diamond-confined subgraph of $\Z^2$.
\end{itemize}
\end{theorem}
\begin{proof}
Fix some $y \in \mathcal{Y}_0^+$. We define a probability distribution on the space
\begin{equation}
    \mathfrak{C}_L \times \bigcup_{l=1}^{+\infty} \left( \prod_{k=1}^\ell \mathfrak{D}\right) \times \mathfrak{C}_R
\end{equation} 
by the formula:
\begin{multline*}
    \phi^{\mathsf{Dec}}_y\left[(G^L, G_1, \dots, G_\ell, G_R)\right] \propto  \1_{X_L(G^L)+X(G_1)+\dots+X(G_\ell)+X_R(G^R) = y}\\ \times\sum_{\substack{\ell \geq 0 \\ G^L\in \mathfrak{C}_L^\delta \\ G^R \in \mathfrak{C}_L^\delta \\ G_1, \dots, G_\ell \in \mathfrak{D}^\delta }}\rho_L(G^L)\rho_R(G^R)\mathbf{P}(G_1) \cdots \mathbf{P}(G_\ell).
\end{multline*}
From~\eqref{equation OZ 1 cluster}, we get that for any percolation event $\mathcal{A}$,
\begin{equation}
\big\vert \phi[\mathcal{A} \vert y \in \mathcal{C}_0] - \phi^{\mathsf{Dec}}_y[\1_\mathcal{A}\circ \mathsf{Conc}] \big\vert \leq C\e^{-c\norme{y}},
\end{equation} 
where $\mathsf{Conc}$ denotes the concatenation operation for a sample of $\phi^{\mathsf{Dec}}_y$.
It is classical that this yields the existence of a maximal coupling between these two measures, \textit{ie} that one can construct a probability space $(\Omega, \mathcal{F}, \Phi_{0 \rightarrow y})$ supporting $(\mathcal{C}^1_0, \mathcal{C}^2_0)$ such that 
\begin{itemize}
    \item The distribution of $\mathcal{C}^1_0$ is the distribution of the cluster of 0 under $\phi\left[~\cdot\vert y \in \mathcal{C}_0\right]$,
    \item The distribution of $\mathcal{C}^2_0$ is the distribution of the concatenation $G^L\circ G^1 \circ \dots \circ G^\ell \circ G^R$ where $(G^L, G^1, \dots, G^\ell, G^R)$ are sampled according to $\phi^\mathsf{Dec}_y$,
    \item $\Phi_{0 \rightarrow y}(\mathcal{C}^1_0 \neq \mathcal{C}^2_0) \leq C\e^{-c\norme{y}}$.
\end{itemize}
Now consider the random variable $\mathcal{S}$ formed from $(G_L, G^1, \dots, G^\ell, G_R)$ by the following formula:
\begin{equation}
   \mathcal{S}_1 = X_L(G_L) ~~\text{and}~~ \mathcal{S}_k = \mathcal{S}_{k-1}+X(G^{k-1}) ~~\text{for}~~ 2\leq k \leq \ell.
\end{equation}
Then it is immediate that 
\begin{equation}
    \Phi_{(0,y)}\left[ \mathcal{S}_1 = s_1, \dots, \mathcal{S}_\ell = s_\ell \right] \propto \nu_L(s_1)\nu_R(y-s_\ell)\prod_{k=2}^{\ell}\nu(s_k-s_{k-1}).
\end{equation}
Moreover by definition, the $\mathcal{S}_k$'s are renewals of $\mathcal{C}^2_0$ and the portions of $\mathcal{C}^2_0$ lying between two consecutive $\mathcal{S}_k$'s are $\delta$-diamond-confined. Thus, $(\Omega, \mathcal{F}, \Phi_{0 \rightarrow y})$ equipped with the random element $(\mathcal{C}^1_0, \mathcal{S})$ provides us with the desired coupling.  
\end{proof}
For now, we shall only work in the extended probability space $(\Omega, \mathcal{F}, \Phi_{0 \rightarrow y})$. Thus, each percolation configuration conditioned to contain the distant point $y$ will be sampled together with a directed random walk bridge: we call this directed random bridge \textbf{the skeleton of $\mathcal{C}_0$}; the associated diamond decomposition will be called \textbf{the diamond decomposition of $\mathcal{C}_0$}. Observe that this enlarged probability space carries extra randomness than the space supporting $\phi$: indeed, to a given percolation cluster can be associated several skeletons that are randomly chosen by the measure $\Phi_{0 \rightarrow y}$(see Figure~\ref{Figure ornstein zernike}). We adopt the terminology of~\cite{ozrandomclustercampanino} by calling the points of $\mathcal{S}$ \emph{renewals} of the cluster. Observe that due to the latter discussion, all the renewals of $\mathcal{C}$ are cone-points, but the converse is not necessarily true.

In the remainder of this paper, we introduce $\PP^\RW$, the measure of the directed random walk with independent increments sampled according to $\nu$ and started from 0.
\begin{remark}
We are often going to be interested in observables of the \textit{skeleton} of a cluster sampled according to $\Phi_{0 \rightarrow y}$. In that case, Theorem~\ref{theoreme couplage oz 1 cluster} reads as follows: let $f$ be a bounded function of the set of directed random walks. Then,
    \begin{multline}\label{equation couplage squelette rw}
        \Bigg| \Phi_{0 \rightarrow y}\left[f(\mathcal{S})\right] - \sum_{x_L, x_R}\nu_L( x_L)\nu_R(y-x_R)\E^\RW\left[f( x_L \circ S \circ x_R ) \vert S \in \hit_{x_R-x_L} \right] \Bigg| \\ \leq C\Vert f \Vert_\infty \e^{-c\Vert y \Vert_2}.
    \end{multline}
In the latter writing, the notation $x_R \circ S \circ x_L$ stands for the directed walk obtained by the concatenation of $x_L$, the trajectory of $S$, and $x_R$:  $S$ is the only random object - and its law is $\PP^\RW\left[~\cdot\vert S \in \hit_{x_R-x_L}\right]$.

Finally, the unconditional version of~\eqref{equation couplage squelette rw} is the following. 
\begin{equation}\label{Equation unconditional skeleton coupling 1 cluster}
\Bigg\vert \e^{\tau y_1}\Phi_{0 \rightarrow y}\left[ f(\mathcal{S}) \right] - \sum_{x_L, x_R}\nu_L(x_L)\nu_R(y-x_R)\E^\RW\left[f(S)\1_{\hit_{x_R-x_L}}\right]  \Bigg\vert \leq C\Vert f\Vert_\infty \e^{-c\Vert y \Vert_2}.
\end{equation}
\end{remark}

 The following lemma states when looking at certain families of observables of the trajectories of directed walks, it is sufficient to study the measure $\PP^\RW$ started from 0 rather than the intricate second summand of the left-hand side of~\eqref{equation couplage squelette rw}

\begin{lemma}\label{lemme skorokhod topology}
    There exists a constant $\varsigma > 0$ such that for any $y := (y_1.y_2)\in \mathcal{Y}^+_0$, any two sequences $a_n, b_n$ of positive numbers going to infinity, any bounded function $f : \mathcal{C}([0,y_1], \R) \rightarrow \R$ continuous with respect to the topology of uniform convergence,
        \begin{multline*}
        \Big\vert \e^{-\tau b_n y_1}\Phi_{0 \rightarrow b_n y} \left[ f\left(a_n^{-1}\mathcal{S}(\lfloor b_n t\rfloor) \right)_{t \geq 0}\right] - \varsigma\E^\RW\left[f\left(a_n^{-1} S(\lfloor b_n t \rfloor )\right)_{t \geq 0}\1_{S \in \hit_{b_n y}} \right] \Big\vert \\ \goes{}{n \rightarrow \infty}{0}.
    \end{multline*}
We have used the interpretation of directed walks as real-valued functions explained above.
\end{lemma}
\begin{proof}
    Set $\varsigma = \nu_L(\Z^2)\nu_R(\Z^2)$. By~\eqref{Equation unconditional skeleton coupling 1 cluster}, it sufficient to prove that 
    \begin{multline*}
    \Big\vert \sum_{x_L, x_R \in \Z^2} \nu_L(x_L)\nu_R(x_R)\E^{\mathsf{RW}}_{x_L}\left[f\left(a_n^{-1}(x_L\circ S\circ x_R)(\lfloor b_n t  \rfloor )\right)_{t \geq 0}\1_{S \in \hit_{x_L-x_R}}\right] \\ - \varsigma\E^\RW\left[f\left(a_n^{-1} S(\lfloor b_n t \rfloor )\right)_{t \geq 0}\1_{S \in \hit_{b_n y}} \right] \Big\vert \goes{}{n \rightarrow \infty}{0}.
    \end{multline*}
The right-hand side can be dominated by
\begin{multline*}
    \sum_{x_L, x_R \in \Z^2} \nu_L(x_L)\nu_R(x_R) \E^\RW\Big[ \big\vert f\left(a_n^{-1}(x_L\circ (S+x_L)\circ x_R)(\lfloor b_n t  \rfloor )\right)_{t \geq 0}\1_{S \in \hit_{x_R - x_L}} \\ - f\left(a_n^{-1} S(\lfloor b_n t \rfloor )\right)_{t \geq 0}\1_{S \in \hit_{b_n y}} \big\vert \Big].
\end{multline*}
Now we take advantage of the exponential tails of $\nu_L$ and $\nu_R$ by splitting the sum in two parts, the first one running over $x_L,x_R \in B(0, \log(\min(a_n, b_n)))$, and the remaining one. Thanks to the exponential tails of $\nu_L$ and $\nu_R$, the remaining one can be bounded by $2\Vert f \Vert_\infty \min(a_n,b_n)^{-c}$, which indeed goes to 0. The first part is shown to go to 0 by noticing that when $x_L,x_R \in B(0, \log(\min(a_n, b_n)))$, the uniform distance between the two considered functions goes to 0. We conclude by continuity of $f$ and dominated convergence.
\end{proof}

We state two byproducts of Theorem~\ref{oz theorem}:

\begin{corollary}\label{corollaire nombre linéaire de renewals}
There exists three constants $c, C, K > 0$ such that for any $y\in \mathcal{Y}^+_0$ with $\norme{y}$ sufficiently large,
\begin{equation}
    \phi \left[ \mathcal{C}_0 \text{ has less than } K\norme{y} \text{ renewal points} \vert y \in \mathcal{C}_0 \right] < C\e^{-c\norme{y}}.
\end{equation}
\end{corollary}

\begin{corollary}\label{corollaire volume diamants}
There exists a constant $K > 0$, such that for any $y \in \mathcal{Y}^+_0$, 

\begin{equation}
    \Phi_{0 \rightarrow y}\Big[ \max_{\substack{\mathcal{D} \subset \mathcal{D}(\mathcal{C}_0) \\ \mathcal{D}\text{ diamond}}} \Vol \left(\mathcal{D}\right) > K (\log y_1)^2 \Big]  < C\norme{y}^{-c\log \norme{y}},
\end{equation}
where $\Vol$ denotes the Euclidean volume, and where the $\max$ is taken over all the diamonds appearing in the diamond decomposition of the cluster of 0 under the measure $\Phi_{(0,y)}$. 
\end{corollary}
Note that the latter bound decays faster than the inverse of any polynomial in $\norme{y}$.

\subsection{Ornstein--Zernike in a strip with boundary conditions}

We import a few facts about Ornstein--Zernike theory that will be useful later on in our analysis. They deal with the uniformity of the Ornstein--Zernike formula in the boundary conditions and are directly imported from~\cite{ozisingcampanino, greenbergioffe}. For $y= (y_1, y_2) \in \mathcal{Y}_{0}^+$, let us call $\mathsf{Strip}_y$ the strip $\mathsf{Strip}_y= \left[0, y_1 \right] \times \Z$. In the following proposition, the probability measure $\mathbf{P}$ is the same object as in Theorem~\ref{oz theorem}.

\begin{proposition}[Uniform OZ formula in a strip]\label{prop oz boundary conditions} Let $y = (y_1,y_2)\in \mathcal{Y}_{0}^+$. 
Let $ C_{\EXT}^L \ni 0$ be a finite connected subset of edges of $\mathcal{Y}^{-}_0$ and $ C_{\EXT}^R \ni y$ be a finite connected subset of edges of $\mathcal{Y}_{y}^+$. Then there exist two positive and bounded measures $\rho_L^{\EXT}, \rho_R^\EXT$ on $\mathfrak{C}_L$ and $\mathfrak{C}_R$ respectively such that for any bounded function $f : \mathfrak{G}_0\rightarrow \R$,
\begin{multline}\label{equation OZ boundary conditions}
  \bigg\vert \e^{\tau y_1}\phi\left[f( \mathcal{C}_0) \1\left(\lbrace\mathcal{C}_0 \cap \mathsf{Strip}_y^c = C^L_\EXT \sqcup C^R_\EXT \rbrace \right) \right] - \\ \sum_{\substack{\ell \geq 0 \\ G^L\in \mathfrak{C}_L \\ G^R \in \mathfrak{C}_L\\ G_1, \dots, G_\ell \in \mathfrak{D} }}\rho^\EXT_L(G^L)\rho^\EXT_R(G^R)\mathbf{P}(G_1) \cdots \mathbf{P}(G_\ell) f(G) \bigg\vert  \leq C\Vert f\Vert_{\infty}\e^{-c\norme{y}},
\end{multline}
 where the sum holds over all $G^L \in \mathfrak{C}_L, G^R \in \mathfrak{C}_R, G_1, \dots, G_k \in \mathfrak{D}$ satisfying the relation
 \begin{equation}
     X^L_0(G^L) + X(G_1)+ \dots + X(G_k) + X^R_x(G^R) = y,
 \end{equation}
 and where we have written $G = G^L \circ G_1 \circ \cdots \circ G_k \circ G^R $.
 Moreover, the measures $\rho^\EXT_L$ and $\rho^\EXT_R$ have exponential tails, uniformly in the sets $C^L_\EXT, C^R_\EXT$ satisfying the above-stated assumptions: indeed, there exist $c', C' > 0$ such that for any $t>0$,
 \begin{equation}
     \sup_{C^L_\EXT, C^R_\EXT} \max \left\lbrace \rho^\EXT_L(X(G^L) > t), \rho^\EXT_R(X(G^R) > t) \right\rbrace < C'\e^{-c't}.
 \end{equation}
\end{proposition}
Observe that in the latter formula, the event $\lbrace\mathcal{C}_0 \cap \mathsf{Strip}_y^c = C^L_\EXT \sqcup C^R_\EXT \rbrace$ implies that $y \in \mathcal{C}_0$.
\begin{remark}
As done previously, for any $y \in \mathcal{Y}_{0}^+$, we shall call 
\begin{equation}
    \nu^\EXT_L(x) = \sum_{G^L \in \mathfrak{C}^\delta_L, X^L_0(G^L) = x}\rho^\EXT_L(G^L) \text{  and  } \nu^\EXT_R(x) = \sum_{G^R \in \mathfrak{C}^\delta_R, X^R_x(G^R) = x}\rho^\EXT_R(G^R)
\end{equation}
\end{remark}
 We simply sketch the proof of the proposition, since it is a simple byproduct of the analysis of~\cite{ozrandomclustercampanino}

\begin{proof}[Proof of Proposition~\ref{prop oz boundary conditions}]
Apply the Ornstein--Zernike formula~\eqref{equation OZ 1 cluster} to the function $g(\mathcal{C}_0) = f(\mathcal{C}_0)\1\left\lbrace \mathcal{C}_0\cap \mathsf{Strip}_y^c = C^L_\EXT \sqcup C^R_\EXT \right\rbrace$. Thus, one has that $\rho_L^\EXT$ (resp. $\rho^R_\EXT$) is the restriction of $\rho_L$ (resp. $\rho_R$) to pieces of clusters compatible with $C^L_\EXT$ (resp. $C^R_\EXT$). The announced exponential decay is a byproduct of the exponential tails of $\rho_L$ and $\rho_R$.
\end{proof}

A non-trivial consequence of the latter proposition is the following estimate, appearing in~\cite[Equation (2.19)]{greenbergioffe}.  

\begin{proposition}[Ornstein--Zernike decay uniform in the boundary conditions]\label{Proposition oz uniforme boundary conditions}
There exists $\chi>0$ such that for any sets $C_\EXT^L, C_\EXT^R$ satisfying the assumptions of Proposition~\ref{prop oz boundary conditions},
\begin{equation}\label{equation oz uniform BC}
    \frac{1}{\chi} \frac{\e^{-\tau n}}{\sqrt{n}} \leq \phi\left[ 0 \overset{\mathsf{Strip}_n}{\longleftrightarrow} (n,0) \Bigg\vert 
        ~\begin{aligned} &\mathcal{C}_0 \cap \left(\Z^- \times \Z\right) = C^L_\EXT ~~\text{and} \\ &\mathcal{C}_{(n,0)} \cap \left( [n, +\infty) \times \Z\right) = C^R_\EXT \end{aligned}   \right] \leq \chi \frac{\e^{-\tau n}}{\sqrt{n}}.
\end{equation}
\end{proposition}

\section{Scaling limit for the product measure}\label{section independent system}

For our purposes, we need to develop an analog of the Ornstein--Zernike theory for $r$ non-intersecting clusters of FK-percolation. However, there is a supplementary difficulty, namely that these non-intersecting clusters are not independent, beyond the obvious interaction introduced by the conditioning. If we consider \textit{a product measure}, we can readily extend the Ornstein--Zernike Theorem to $r$ clusters sampled independently according to $\phi$. This is the goal of the present section. Even though it might seem a bit strange to consider the conditioned product measure $\phi^{\otimes r}$ instead of the real conditioned random-cluster measure, we are going to see in Section~\ref{section RCM} that because of the conditioning, these two measures behave similarly "in the bulk". This is a consequence of the spatial mixing property of the subcritical random-cluster measure combined with an a priori repulsion estimate. 

In what follows, $\phi^{\otimes r}$ will always denote the measure consisting in the product of $r$ random-cluster measures $\phi$. Moreover, $\mathcal{C}_i$ will denote $\mathcal{C}_i(\omega_i)$. In particular, if $\mathcal{A}$ is an event measurable with respect to $\left(\mathcal{C}_1, \dots, \mathcal{C}_r\right)$, we have:
\begin{equation}
\phi^{\otimes r}\left[\mathcal{A}\right] = \PP \left[\left( \mathcal{C}_1(\om_1), \dots, \mathcal{C}_r(\om_r) \right) \in \mathcal{A} \right],     
\end{equation}
 where $\om_1, \dots, \om_r$ are \textbf{independent} percolation configurations sampled according to $\phi$.

The main goal of this section is the following proposition:

\begin{proposition}\label{Proposition convergence mesure produit}
Recall the definition of the envelopes of a cluster $\Gamma^{\pm}(\mathcal{C})$ introduced in Def~\ref{definition interfaces}, and their natural parametrization. Then there exists $\sigma > 0$ such that:
\begin{equation}
 \frac{1}{\sqrt{n}} \left( \Gamma^+(\mathcal{C}_1)(nt), \dots, \Gamma^+(\mathcal{C}_r)(nt) \right)_{0 \leq t \leq 1} \goes{(d) }{n \rightarrow \infty}{\left(\sigma\bw^{(r)}_t \right)_{0 \leq t \leq 1}},
 \end{equation}
where the percolation configuration is sampled under the measure $\phi^{\otimes r} \left[ ~\cdot ~ \big\vert \con,\nonint \right]$, and the convergence occurs in the space $\mathcal{C}^r([0,1])$ equipped with the topology of uniform convergence. Moreover, almost surely, for any $1\leq i\leq r$,
\begin{equation}
    \frac{1}{\sqrt{n}}\Vert \Gamma_i^+ - \Gamma_i^- \Vert_{\infty} \goes{}{n\rightarrow +\infty}{0}.
\end{equation}
\end{proposition}

The strategy of the proof is the following: we start to state an analog of the Ornstein--Zernike Theorem in the case of a product measure in Section~\ref{subDefinition of the product measure and multidimensional version of Ornstein--Zernike theorem }, and use this coupling to compare the skeletons of a system of $r$ non-intersecting clusters with a system of $r$ non-intersecting directed random walks. However, there is a small difficulty while implementing this program: indeed, conditioning on the event that the \textit{clusters} do not intersect is not the same as conditioning on the event that the \textit{skeletons} of the clusters do not intersect. Moreover, while the latter event is very well described in terms of the Ornstein--Zernike coupling, it is not a priori clear how the former acts on the coupled walks. For that reason, we shall show first that under the conditioning on $\nonint \cap \con$, the clusters very soon get far from each other (this is the goal of Subsection~\ref{sub edge repulsion}), and thanks to this input we will be able to prove that "in the bulk" of the system, the conditioning of non-intersection for the clusters or for the skeletons of the clusters yield the same scaling limit. Thus we shall apply the invariance principle of Theorem~\ref{theoreme invariance principle for drw} to conclude in Section~\ref{sub convergence towards the brownian watermelon}.

\subsection{Multidimensional version of the Ornstein--Zernike Theorem in the product measure}\label{subDefinition of the product measure and multidimensional version of Ornstein--Zernike theorem }

We first state a $r$-dimensional version of Theorem~\ref{oz theorem} (the Ornstein--Zernike formula). Indeed, if $\om^1, \dots, \om^r$ are $r$ independent configurations of law $\phi$ and $f : \mathfrak{G}^r \rightarrow \R$ is a bounded function, then it is an easy consequence of Theorem~\ref{oz theorem} that when $n \in \N$ is sufficiently large,
\begin{multline}\label{multidim oz}
\Bigg\vert \e^{\tau r n}\phi^{\otimes r}\left[ f(\mathcal{C}_1, \dots, \mathcal{C}_r) \1_{(\om^1, \dots, \om^r) \in \con} \right]  -  \sum_{\substack{G_1^L\in \mathfrak{C}_L \\ \dots \\ G_r^L \in \mathfrak{C}_L}} \sum_{\substack{G_1^R \in \mathfrak{C}_R \\ \dots \\ G_r^R \in \mathfrak{C}_R}}\sum_{\substack{k_1 \geq 0 \\ \dots \\ k_r \geq 0}} \sum_{G_1^1, \dots, G^{k_1}_1 \in \mathfrak{D} } \dots \sum_{G_r^1, \dots, G^{k_r}_r \in \mathfrak{D}}  \\ \bigg(\prod_{i=1}^r \rho_L(G_i^L)\rho_R(G_i^R)\mathbf{P}(G^1_i) \cdots \mathbf{P}(G^{k_i}_i)\bigg) f(G_1, \dots, G_r) \Bigg\vert 
 \leq Cr\Vert f\Vert_{\infty}\e^{-cn}.
\end{multline}
where we sum over all the $G_1^L, \dots, G_r^L \in \mathfrak{C}_L$, the $G_1^R, \dots, G_r^R \in \mathfrak{C}_R$, and the families $G_1^1, \dots, G_1^{k_1}, \dots, G^1_r, \dots,G^{k_r}_r \in \mathfrak{D}$ such that for any $1 \leq i \leq r$,
\begin{equation}
    X^L_0(G_i^L)+ X(G^i_1)+\dots+X(G^{k_i}_i) +X^R_x(G_i^R) = n(y_i-x_i).
\end{equation}

The coupling stated in Theorem~\ref{theoreme couplage oz 1 cluster} is available in this context: we call it $\Phi^{\otimes r}_{(0,x)\rightarrow (n,y)}$. It simply consists in the product of the $r$ couplings $\Phi_{(0,x_j)\rightarrow (n,y_j)}$ given by Theorem~\ref{theoreme couplage oz 1 cluster}. Its main feature is that for any bounded function of the skeletons of a system of $r$ clusters, 
\begin{multline}\label{equation couplage squelette p rw} \Bigg\vert
 \Phi^{\otimes r}_{(0,x)\rightarrow (n,y)} \left[ f\left(\mathcal{S}_1, \dots, \mathcal{S}_r \right)\right] - \\ \sum_{\substack{x_L^1, \dots, x_L^r \\ x_R^1, \dots, x_R^r}}\left(\prod_{i=1}^r \nu_L(x_L^i)\nu_R(x_R^i)\right)\left(\mathbb{E}^{\RW}\right)^{\otimes r}\left[f(x_L^1\circ S^1\circ x_R^1, \dots, x_L^r\circ S^r \circ x_R^r)\1_{S\in \hit_{ x_R - x_L}} \right] \Bigg\vert  \\ < Cr\Vert f \Vert_{\infty}\e^{-cn},
\end{multline}
where we denoted by $(\E^\RW)^{\otimes r}$ the expectation under the measure of $r$ independent directed walks $(S^1, \dots, S^r)$ started from 0, and $\hit_{x_R - x_L}$ the event that each $S^i$ ever hits $x_R^i - x_L^i$. Note that the uniform Ornstein--Zernike coupling introduced in Proposition~\ref{prop oz boundary conditions} also holds in this context. Moreover, Lemma~\ref{lemme skorokhod topology} is also true in its $r$-dimensional version, so that it will be sufficient to study $(\E^\RW)^{\otimes r}$ when estimating probabilities for \textit{scaled} random walks.

Before working on the repulsion estimates as announced, we lower bound the probability of non-intersection and connection in the product measure. To this end, we will need to import a local limit theorem from Section~\ref{section marches}. We here present it in a very informal way but refer the reader to Theorem~\ref{theoreme local limit srw} for a complete statement --- the main purpose of what follows is the introduction of the function $V$. 

\begin{theoremnonumber}[Theorem~\ref{theoreme local limit srw}; very informal version]
There exists a function $V: W\cap\Z^r \rightarrow \R^+$ and a constant $C_1 > 0$ such that the probability that $r$ synchronized random walks reach the point $(n,y)$ without intersecting is of order 
\begin{equation*}
C_1\frac{V(x)V(y)}{n^{\frac{r^2}{2}}}(1+o(1)). 
\end{equation*}
Furthermore this estimate holds when $x,y$ are allowed to depend on $n$, as long as $\norme{x(n)},\norme{y(n)} = o(\sqrt{n})$ .
\end{theoremnonumber}
The desired lower bound is given by the following lemma. 
\begin{lemma}
\label{lemme minoration proba con nonint mesure produit} Let $x,y \in W\cap\Z^r$. Then, there exists $c>0$ such that 
\begin{equation}\label{equation lemme minoration proba con nonint mesure produit}
    \phi^{\otimes r}\left[ \nonint, \con \right] \geq cV(x)V(y) n^{- \frac{r^2}{2}}\e^{-\tau rn},
\end{equation}
where $V$ is the function introduced in Theorem~\ref{theoreme local limit srw}.
\end{lemma}

\begin{proof}
We use the Ornstein--Zernike coupling given by~\eqref{equation couplage squelette p rw}: indeed, up to exponential terms due to the coupling, and using the diamond confinement property,

\begin{align*}
&\phi^{\otimes r}\left[ \nonint, \con \right] \\ &= \e^{-\tau rn}\Phi^{\otimes r}_{(0,x)\rightarrow (n,y)}\left[(\mathcal{C}_1, \dots, \mathcal{C}_r) \in \nonint \right] \\
&\geq \e^{-\tau rn}\Phi^{\otimes r}_{(0,x)\rightarrow (n,y)} \left[\bigcap_{1\leq i\neq j\leq r}\lbrace\mathcal{D}(\mathcal{S}^i)\cap\mathcal{D}(\mathcal{S}^j)=\emptyset\rbrace \right] \\
&= \e^{-\tau rn}\sum_{\substack{x^1_L, \dots, x_L^r \\ x^1_R, \dots, x_R^r}}\left(\prod_{i=1}^r \nu_L(x^i_L) \nu_R(x^i_R)\right)\\  &\hspace{20pt}\times c\left(\PP^\RW\right)^{\otimes r}\left[\bigcap_{1\leq i\neq j\leq r}\lbrace \mathcal{D}\left(x^i_L \circ S^i \circ x^i_R\right)\cap \mathcal{D}\left(x^j_L \circ S^j \circ x^j_R\right) = \emptyset\rbrace , \hit_{(n,y - x)}\right].
\end{align*}
Hence the result boils down to lower bound the probability of non-intersection and connection for $r$ independent \textit{decorated} directed random walks. This is precisely the content of Lemma~\ref{lemme estimation proba non-synchronized RW}. By finiteness of the measures $\nu_L, \nu_R$, we conclude that
\begin{equation}
    \phi^{\otimes r}\left[ \nonint, \con \right] \geq cV(x)V(y) n^{-\frac{r^2}{2}}\e^{-\tau rn}.
\end{equation}
\end{proof}
Observe that the latter bound reads as follows on the coupling measure:
\begin{equation}\label{eq: lower bound proba nonint conditionnelle}
    \Phi^{\otimes r}_{(0,x)\rightarrow (n,y)}\left[ \nonint \right] \geq cV(x)V(y)n^{-\frac{r^2 }{2}}.
\end{equation}

Thanks to Proposition~\ref{prop oz boundary conditions}, the same analysis holds for deriving the analog of the latter result in a strip with boundary conditions. 
\begin{corollary}\label{corollaire majoration proba ni con mesure produit BC avec ecartement suffisant}
Let $x,y \in W$. Let $ C^L_{i,\EXT} \ni x_i,  C^R_{i,\EXT} \ni y_i$, and assume that the family $C = \left(C^L_{i,\EXT}, C^R_{i,\EXT}\right)_{1\leq i \leq r}$ satisfies the assumptions of the uniform Ornstein--Zernike coupling given by Proposition~\ref{prop oz boundary conditions}. Then, there exists a uniform constant $\chi > 0$ such that 
\begin{equation}
    \phi^{\otimes r}\left[ \con, \nonint ~\Bigg\vert \begin{aligned}
        &\mathcal{C}_{(0,x_i)} \cap (\Z^-\times\Z) = C^L_{i, \EXT},  \\
        &\mathcal{C}_{(n,y_i)} \cap ([n,+\infty]\times\Z) = C^R_{i,\EXT}, ~~~\forall 1 \leq i \leq r
    \end{aligned} \right] \geq \frac{c}{\chi}V(x)V(y)n^{-\frac{r^2}{2}}\e^{-\tau rn}. 
\end{equation}
\end{corollary}

\subsection{Edge repulsion}\label{sub edge repulsion}

The goal of this section is to prove Lemma~\ref{equation entropic repulsion independent system} which we refer to as the "edge repulsion" lemma for the independent system. Beforehand we introduce the important notion of \emph{synchronized skeleton} of a system of long clusters. 

Let $\mathcal{C} = (\mathcal{C}_1, \dots, \mathcal{C}_r)$ be sampled according to $\Phi^{\otimes r}_{(0,x)\rightarrow (n,y)}$. We say that $k \in \N$ is a \textit{synchronization time} for $\mathcal{C}$ if there exists $(s^1_k, \dots, s^r_k) \in \Z^r$ such that for any $1\leq i \leq r$, one has $(k, s^i_k) \in \mathcal{S}^i$. In other words, $k$ is a synchronization time for $\mathcal{C}$ if and only if each one of the $r$ skeletons of $\mathcal{C}$ contains a point of $x$-coordinate $k$. We define the set of synchronization times of $\mathcal{C}$ by 
\begin{equation}
    \mathsf{ST}(\mathcal{C}) = \lbrace 0 \leq k_1 < k_2 < \dots < k_l \leq n\rbrace.
\end{equation}
The \textit{synchronized skeleton} of $\mathcal{C}$, called $\check{\mathcal{S}}$ is now defined to be the process defined on $\mathsf{ST}(\mathcal{C})$, taking its values in $\Z^r$, such that for any $k \in \mathsf{ST}(\mathcal{C})$,
\begin{equation}
    \check{\mathcal{S}}_k = (\mathcal{S}^1(k), \dots, \mathcal{S}^r(k)).
\end{equation}
As previously, we extend this process as a function of $\R^+$ to $\R$ by linear interpolation.
Let us observe an important property of this process (which is the reason of its introduction) before turning to Lemma~\ref{Lemme entropic repulsion independent system}. 

\begin{figure}
     \centering
     \includegraphics[width = .45\textwidth, page = 1]{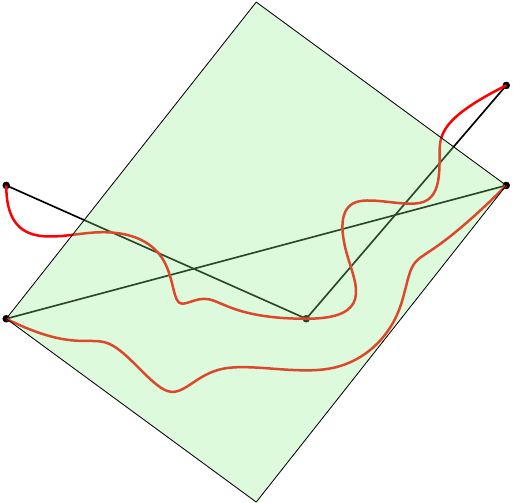}
     \caption{Illustration of the necessity of considering the process of synchronized renewals: here, the red clusters do not intersect while their associated skeletons do intersect. }
     \label{figure mauvaise intersection}
 \end{figure}

\begin{claim}
Under $\Phi^{\otimes r}_{(0,x)\rightarrow (n,y)}$, 
\begin{equation}\label{claim inclusion synchronized events}
   \left\lbrace \left(\mathcal{C}_1, \dots, \mathcal{C}_r \right) \in \nonint \right\rbrace \subseteq \left\lbrace \left( \check{\mathcal{S}}(\mathcal{C}_1), \dots, \check{\mathcal{S}}(\mathcal{C}_r)) \right) \in \mathcal{W}_{[0,n]} \right\rbrace.
\end{equation}
\end{claim}
\begin{proof}
    This is an immediate consequence of the Intermediate Value Theorem.
\end{proof}
\begin{remark}
    Observe that this inclusion would not be true when replacing $\check{\mathcal{S}}$ by $\mathcal{S}$ (see Figure~\ref{figure mauvaise intersection}).
\end{remark}

Moreover, due to the exponential tails of the increments of $\mathcal{S}$, it is clear that the increments of $\check{\mathcal{S}}$ also have exponential tails. Thus, $\check{\mathcal{S}}$ falls into the class of \textit{synchronized directed random walks}, studied in Section~\ref{subsection non-intersecting systems of synchronized directed random walks}.

Next lemma indicates that in a time less than $n-o(n)$ the clusters have been far away from at least $n^\eps$ at least once. A convenient notion for stating this result is the gap of a point $x \in W$.

\begin{definition}[Gap of a point]
    Let $x \in W$. We define its \emph{gap} to be the following quantity:
    \begin{equation}
        \Gap(x) = \min_{1\leq i \leq r-1}(x_{i+1} - x_i).
    \end{equation}
    Observe that due to the fact that $x$ lies in $W$, $\Gap(x)$ is always a positive quantity.
\end{definition}

We are ready to state the edge repulsion result for the independent system.

\begin{lemma}[Edge repulsion for the independent system]\label{Lemme entropic repulsion independent system}
There exists an $\eps > 0$ such that the following holds. Let $T_1$ and $T_2$ be the following random variables:
\begin{equation}
    T_1 = \inf \left\lbrace k \geq 0, \Gap(\check{\mathcal{S}}_k) > n^\eps \right\rbrace.
\end{equation}
and
\begin{equation}
    T_2 = \sup \left\lbrace k \geq 0, \Gap(\check{\mathcal{S}}_k)  > n^\eps \right\rbrace
\end{equation}
Then, there exist $C,c>0$ such that when $n\geq 0$ is large enough,
\begin{equation}\label{equation entropic repulsion independent system}
    \Phi^{\otimes r}_{(0,x)\rightarrow (n,y)} \left[ \left\lbrace T_1 > n^{1-\eps}\right\rbrace \cup \left\lbrace T_2 < n - n^{1-\eps} \right\rbrace \vert  \nonint \right] < \frac{1}{c}\exp(-cn^\eps).
\end{equation}
\end{lemma}

\begin{proof}
We prove that $\Phi^{\otimes r}_{(0,x)\rightarrow (n,y)}\left[ T_1 > n^{1-\eps} \vert\nonint \right] < \exp(-cn^\eps)$. By time reversal and a basic union bound, it will be sufficient to conclude. We roughly upper bound:

\begin{equation}\label{e: numerateur denominateur preuve edge rep indep}
    \Phi^{\otimes r}_{(0,x)\rightarrow (n,y)}\left[ T_1 > n^{1-\eps} \vert\nonint \right] \leq \frac{ \Phi^{\otimes r}_{(0,x)\rightarrow (n,y)}\left[ T_1 > n^{1-\eps} \right]  }{ \Phi^{\otimes r}_{(0,x)\rightarrow (n,y)}\left[ \nonint  \right]
}.\end{equation}

We are going to separately bound the numerator and the denominator. We start by the numerator of~\eqref{e: numerateur denominateur preuve edge rep indep}. Since $T_1$ is measurable with respect to the synchronized skeleton of $\mathcal{C}$ (which itself is measurable with respect to $\mathcal{S}$), we use the fact that the law of $\mathcal{S}$ under $\Phi^{\otimes r}_{(0,x)\rightarrow (n,y)}$ is that of a system of directed random bridges to write - up to an exponential correction due to the coupling:
\begin{equation}
    \Phi^{\otimes r}_{(0,x)\rightarrow (n,y)}\left[ T_1 > n^{1-\eps} \right]  = \PP^{\RW}_{(0,x)}\left[ T_1(\check{S}) > n^{1-\eps} \vert \hit_{(n,y)} \right].
\end{equation}
We then are exactly in the context of entropic repulsion for synchronized directed random bridges, and we refer to Corollary~\ref{cor: repulsion SDRbridges}, which asserts that the latter probability is upper bounded by $c^{-1}\exp(-cn^\eps)$ for some constant $c>0$.

By~\eqref{eq: lower bound proba nonint conditionnelle}, we have a polynomial lower bound on the denominator. Hence, up to slightly changing the value of $c$, we obtained
\begin{equation}
    \Phi^{\otimes r}_{(0,x)\rightarrow (n,y)}\left[ T_1 > n^{1-\eps} \vert\nonint \right] < \frac{1}{c}\exp(-cn^\eps),
\end{equation}
which achieves the proof.
\end{proof}

By the synchronized renewal time property, for any $1\leq i\leq r$, the cluster $\mathcal{C}_i$ intersects the line $\lbrace T_1 \rbrace \times \Z$ (resp. $\lbrace T_2 \rbrace \times \Z$) at a unique vertex of $\Z^2$, whose $y$-coordinate shall be called $\mathscr{X}_i$ (resp. $\mathscr{Y}_i$). Thus, $\mathscr{X}$ and $\mathscr{Y}$ are elements of $\Z^r$ satisfying $\mathscr{X}_1< \dots < \mathscr{X}_r$ (resp. $\mathscr{Y}_1< \dots < \mathscr{Y}_r$). Introduce the following \emph{edge-regularity} condition:

\begin{definition}[Edge-regularity property]\label{def: definition edge regularity property}
Let $\om \in \con\cap\nonint$ be a percolation configuration. We call $\om$ \emph{edge-regular} an abbreviate this event in $\mathsf{EdgeReg}$ if it satisfies the following properties:
\begin{itemize}
    \item $T_1 < n^{1-\eps}$ and $T_2 > n-n^{1-\eps}$,
    \item $\norme{\mathscr{X}} \leq n^{1/2 - \eps/4}$ and $\norme{\mathscr{Y}} \leq n^{1/2 - \eps/4}$
\end{itemize}
where $\eps>0$ is given by Lemma~\ref{Lemme entropic repulsion independent system}.
\end{definition}
We then prove that a percolation configuration sampled under $\Phi^{\otimes r}_{(0,x)\rightarrow (n,y)}\left[~\cdot\vert\nonint \right]$ is typically edge-regular.
\begin{lemma}\label{lem: edge regularity product measure}
    There exists a small constant $c>0$ such that 
    \begin{equation}
        \Phi^{\otimes r}_{(0,x)\rightarrow (n,y)}\left[\mathsf{EdgeReg}^c\vert\nonint\right] < \frac{1}{c}\exp(-cn^{\frac{\eps}{2}}).
    \end{equation}
\end{lemma}
\begin{proof}
    Let us work conditionally on the event $T_1 < 2n^{1-\eps}$, as it has been proved to occur with exponentially large probability in Lemma~\ref{Lemme entropic repulsion independent system}. As previously, we use the rough upper bound
    \begin{eqnarray*}
        \Phi_{(0,x)\rightarrow(n,y)}^{\otimes r}\left[ \norme{\mathscr{X}} > n^{1/2-\eps/4}\vert\nonint\right] &=& \Phi_{(0,x)\rightarrow(n,y)}^{\otimes r}\left[ \norme{\mathcal{S}(T_1)}>n^{1/2-\eps/4} \vert\nonint \right] \\
        &=& \frac{\Phi_{(0,x)\rightarrow(n,y)}^{\otimes r}\left[ \norme{\mathcal{S}(T_1)}>n^{1/2-\eps/4}\right]}{\Phi_{(0,x)\rightarrow(n,y)}^{\otimes r}\left[\nonint\right]}.
    \end{eqnarray*}
    As previously we use the lower bound~\eqref{eq: lower bound proba nonint conditionnelle} to argue that the denominator is at least polynomial, while we are going to produce a stretched-exponential upper bound on the numerator. First, observe that
    \begin{eqnarray*}
        \Phi_{(0,x)\rightarrow(n,y)}^{\otimes r}\left[ \norme{\mathscr{X}} > n^{1/2-\eps/4}\right]  &=& \PP^{\RW}_{(0,x)}\left[ \norme{S(T_1)} > n^{1/2-\eps/4} \vert\hit_{(n,y)} \right] \\ &\leq& \frac{\PP^{\RW}_{(0,x)}\left[ \norme{S(T_1)} > n^{1/2-\eps/4} \right]}{\PP^{\RW}_{(0,x)}\left[ \hit_{(n,y)}\right]}.
    \end{eqnarray*}
By Theorem~\ref{theoreme local limit srw}, the denominator is at least polynomial.
Now observe that the classical theory of large deviations for random walks allows us to produce a stretched-exponential upper bound on the numerator (remember that we work conditionally on $T_1 < n^{1-\eps}$): there exists $c>0$ such that
\begin{equation}
    \PP^{\RW}_{(0,x)}\left[\norme{S(T_1)} > n^{\frac{1-\eps}{2}+\eps/4}\right] \leq \exp(-cn^{\frac{\eps}{2}}).
\end{equation}
This proves, up to some small change in the constant $c$, that
\begin{equation}
\Phi_{(0,x)\rightarrow(n,y)}^{\otimes r}\left[ \norme{\mathscr{X}} > n^{1/2-\eps/4}\vert T_1 \leq  n^{1-\eps}, \nonint\right] \leq \exp(-cn^{\frac{\eps}{2}}).
\end{equation}
We conclude writing (the factor 2 comes from the terms in $T_2$ and $\mathscr{Y}$ that are handled by symmetry):
\begin{multline*}
    \Phi_{(0,x)\rightarrow(n,y)}^{\otimes r}\left[ \mathsf{EdgeReg}^c\vert\nonint\right] \leq \\ 2\left(2\Phi_{(0,x)\rightarrow(n,y)}^{\otimes r}\left[ T_1 > n^{1-\eps}\vert\nonint\right]+\Phi_{(0,x)\rightarrow(n,y)}^{\otimes r}\left[ \norme{\mathscr{X}} > n^{1/2-\eps/4}\vert T_1 \leq n^{1-\eps}, \nonint\right]\right).
\end{multline*}
Thus,
\begin{equation}
    \Phi_{(0,x)\rightarrow(n,y)}^{\otimes r}\left[ \mathsf{EdgeReg}^c\vert\nonint\right] \leq \frac{1}{c}\exp(-cn^{\frac{\eps}{2}}).
\end{equation}
\end{proof}

\subsection{Convergence towards the Brownian watermelon}\label{sub convergence towards the brownian watermelon}

As we shall explain here, the edge repulsion stated in Lemma~\ref{Lemme entropic repulsion independent system} is the main ingredient needed to show that the rescaled system, \textit{sampled under the product measure} and conditioned both on the mutual avoidance of the clusters and on the connection event converges in distribution towards the Brownian watermelon. 

The technique of proof will be used several times through the paper. Basically, it consists in splitting the system of clusters sampled under the measure $\phi^{\otimes r}[\cdot \vert \con, \nonint]$ in two different parts (for the sake of exposition, we explain the splitting on the first half of the cluster "near the starting point" - of course, one has to do the symmetric splitting "near the arrival point"). The first part will be given by the random time $T_1$ introduced in Lemma~\ref{Lemme entropic repulsion independent system}. At this time, the clusters are far from each other, sufficiently far for the conditioning on the non-intersection of the \textit{clusters} to be asymptotically equivalent to the conditioning on the non-intersection of the \textit{skeletons} of the clusters. This allows us to implement the Ornstein--Zernike coupling given by~\eqref{multidim oz} for the section of the clusters which is after $T_1$ (taking into account the boundary conditions enforced by the configuration outside of the strip thank to Proposition~\ref{prop oz boundary conditions}). We conclude by applying the invariance principle for directed random walks derived in Section~\ref{subsection non intersecting non-synchronized rw}.

Due to the fact that we work between the random times $T_1$ and $T_2$ we need a technical input that allows us to extend the convergence as a process on the interval $(0,1)$ to the convergence as a process defined on $[0,1]$.

\begin{lemma}\label{Lemme technique processus sto}
    Let $G_n$ be a random sequence of functions of the space $\mathcal{C}([0,1], \R^r)$ and $G$ be a continuous stochastic process of $\mathcal{C}([0,1], \R^r)$. Assume that: 
    \begin{itemize}
        \item For any $\delta>0$, for any bounded and continuous function $f: \mathcal{C}([\delta, 1-\delta], \R^r) \rightarrow \R$,
        \begin{equation}
            \E\left[f(\restriction{G_n}{[\delta, 1-\delta]})\right] \goes{}{n \rightarrow \infty}{\E\left[f(\restriction{G}{[\delta, 1-\delta]})\right]},
        \end{equation}
        \item For all $\eps > 0$
        \begin{equation}
        \lim_{t\rightarrow 0}\sup_{n\geq 0} \PP\left[ \big\vert G_n(t) - G_n(0)\big\vert > \eps\right]   = 0
        \end{equation}
        and
        \begin{equation}
        \lim_{t\rightarrow 1}\sup_{n\geq 0} \PP\left[ \big\vert G_n(t) - G_n(1)\big\vert>\eps \right]   = 0
        \end{equation}
    \end{itemize}
    Then, $G_n$ converges in distribution towards $G$ in the space $\mathcal{C}([0,1], \R^r)$.
\end{lemma}

\begin{proof}[Sketch of proof of Lemma~\ref{Lemme technique processus sto}]
The proof of Lemma~\ref{Lemme technique processus sto} relies on very classical arguments and we refer to~\cite{billingsley} for details. Observe that hypothesis $(i)$ together with the fact that the family $[\delta, 1-\delta]$ is a compact exhaustion of $(0,1)$ yields the convergence of $G_n$ towards $G$ as processes from $(0,1)$ to $\R^r$. The equicontinuity of $G_n$ at 0 and 1 (hypothesis $(ii)$) then yields the desired convergence by the Arzelà-Ascoli Theorem.
\end{proof}
This technical tool in hand, we can prove the main result of this section.

\begin{proof}[Proof of Proposition~\ref{Proposition convergence mesure produit}]
In what follows, introduce the scaled version of $\mathcal{S}$, for any $0\leq t\leq 1$:
\begin{equation}
    \mathcal{S}_n(t) = \tfrac{1}{\sqrt{n}}\mathcal{S}(nt).
\end{equation}

We are going to implement the strategy given by Lemma~\ref{Lemme technique processus sto} to show that under the measure $\Phi^{\otimes r}_{(0,x)\rightarrow (n,y)}[~\cdot\vert\nonint]$, the scaled system of skeletons $\mathcal{S}_n$ converges towards the Brownian watermelon as random functions of $\mathcal{C}([0,1], \R^r)$.

We start with the proof of point (i) (the crucial part of the proof). Fix $\delta > 0$. Fix $f^\delta : \mathcal{C}([\delta, 1-\delta], \R^r) \rightarrow \R$, continuous and bounded. Our goal is to show that there exists $\sigma > 0$, independent of $\delta$, such that
\begin{equation}\label{equation a montrer convergence FDD mesure produit}
    \Phi^{\otimes r}_{(0,x)\rightarrow (n,y)}\left[ f^\delta(\restriction{\mathcal{S}_n}{[\delta, 1-\delta]})\vert \nonint \right] \goes{}{n\rightarrow\infty}{\E\left[f^\delta(\restriction{\sigma\bw^{r}}{[\delta, 1-\delta]})\right]}.
\end{equation}
For sake of notational simplicity, the restrictions of the functions $\mathcal{S}_n$ and $\bw^{(r)}$ to the interval $[\delta, 1-\delta]$ will not be made explicit anymore.  

We first observe that, by Lemma~\ref{lem: edge regularity product measure}, and using the fact that $f^\delta$ is bounded, 

\begin{equation}
    \Phi^{\otimes r}_{(0,x)\rightarrow (n,y)}\left[ f^\delta(\mathcal{S}_n)\vert \nonint \right] = (1+o(1))\Phi^{\otimes r}_{(0,x)\rightarrow (n,y)}\left[ f^\delta(\mathcal{S}_n)\vert \nonint, \mathsf{EdgeReg} \right].
\end{equation}
Hence, it is sufficient to establish the convergence~\eqref{equation a montrer convergence FDD mesure produit} for the measure conditioned on the configuration to be edge-regular. We recall that under this conditioning, there exist $T_1$ and $T_2$ such that:
\begin{itemize}
    \item $T_1 < n^{1-\eps}, T_2 > n-n^{1-\eps}$ and $T_1$ and $T_2$ are synchronized renewal times of $\mathcal{C}$.
    \item $\Gap(\mathscr{X}), \Gap(\mathscr{Y}) > n^\eps$.
    \item $\norme{\mathscr{X}}, \norme{\mathscr{Y}} < n^{1/2 - \eps/4}$.
\end{itemize}

We chose $n$ large enough so that $n\delta > T_1$ and $n(1-\delta) < T_2$\footnote{This is the only place where we use the fact that our functions are defined on $[\delta, 1-\delta]$ and it is the reason why we follow the strategy given by Lemma~\ref{Lemme technique processus sto} instead of directly working on $[0,1]$.}.

We call $\mathsf{Strip} := [T_1, T_2]\times\Z$. We are going to use an exploration argument, by conditioning on the portion of the clusters $\mathcal{C}$ that lies outside of $\mathsf{Strip}$. To that end, for an edge-regular percolation configuration $\om \in \nonint \cap \con$, introduce the following sets of vertices:
\begin{equation}
    \EXT_i = (\mathcal{C}_i \cup \partial_{\mathsf{ext}}\mathcal{C}_i) \cap \mathsf{Strip}^c \qquad \text{and} \qquad \EXT = \bigcup_i \EXT_i.
\end{equation}

Now observe that summing over all the possible exterior edge-regular configurations yields
\begin{multline}\label{equation conditionnement preuve invariance systeme produit}
    \Phi^{\otimes r}_{(0,x)\rightarrow (n,y)}\left[ f^\delta(\mathcal{S}_n) \vert
        \nonint, \mathsf{EdgeReg} \right] = \\ \sum_{\mathsf{Ext}} \Phi^{\otimes r}_{(T_1,\mathscr{X})\rightarrow (T_2,\mathscr{Y})}\left[f^\delta(\mathcal{S}_n) \vert \nonint, \EXT = \mathsf{Ext} \right]\\\times\Phi_{(0,x) \rightarrow (n,y)}^{\otimes r}\left[\EXT = \mathsf{Ext} \vert \nonint, \mathsf{EdgeReg}\right].
\end{multline}

We would like to conclude using Theorem~\ref{theoreme invariance principle for drw}. For that, we need to understand how the measure changes when switching the conditioning from $\nonint$ to $\nidiam(\mathcal{S})$ (this is the event appearing in the statement of Theorem~\ref{theoreme invariance principle for drw} that the \textit{decorated} skeletons do not intersect).

Fix such an admissible edge-regular $\mathsf{Ext}$. We use the following important input from Section~\ref{subsection non intersecting non-synchronized rw}. By edge-regularity of $\mathsf{Ext}$, usual properties of the coupling and Lemma~\ref{lemme repulsion globale non-synchronized RW} ensure that
\begin{equation}
    \Phi^{\otimes r}_{(T_1,\mathscr{X})\rightarrow (T_2,\mathscr{Y})}\left[\inf_{t\in [T_1, T_2]} \Gap(\mathcal{S}(t)) \leq (\log n)^3 \vert \nidiam(\mathcal{S}), \EXT = \mathsf{Ext} \right] \goes{}{n \rightarrow \infty}{0}.
\end{equation}
Now under the complementary event, the diamond confinement property given by the Ornstein--Zernike coupling ensures that ($\Delta$ here means the symmetric difference) the event \\ $\nidiam(\mathcal{S}) \Delta\lbrace \mathcal{C}\in \nonint \rbrace$  can occur only if one of the diamonds appearing in the diamonds decompositions of the $\mathcal{C}_i$ has a volume larger than $(\log^3 n)^2$. This event has been shown in Corollary~\ref{corollaire volume diamants} to occur with probability going to 0. Thus, we proved that:
\begin{equation}
    \Phi^{\otimes r}_{(T_1,\mathscr{X})\rightarrow (T_2,\mathscr{Y})}\left[ \nidiam(\mathcal{S}) \Delta \left\lbrace \mathcal{C}\in \nonint \right\rbrace \vert \EXT = \mathsf{Ext}\right]   \goes{}{n\rightarrow\infty}{0}.
\end{equation}
It is an easy consequence that:
\begin{multline}\label{equation conditionnement equivalents mesure prod}
    \Big\vert  \Phi^{\otimes r}_{(T_1,\mathscr{X})\rightarrow (T_2,\mathscr{Y})}\left[ f^\delta(\mathcal{S}_n) \vert 
        \nonint,\EXT = \mathsf{Ext}  \right] - \\ \Phi_{(T_1,\mathscr{X}) \rightarrow (T_2,\mathscr{ Y})}^{\otimes r}\left[f^\delta(\mathcal{S}_n) \vert 
            \nidiam(\mathcal{S}),  \EXT = \mathsf{Ext}  \right] \Big\vert \goes{}{n \rightarrow\infty}{0}.
\end{multline}
The right-hand summand of the latter formula is measurable with respect to $\mathcal{S}$, except the conditioning on $\EXT = \mathsf{Ext}$. Thanks to the uniform Ornstein--Zernike formula stated in Proposition~\ref{Proposition oz uniforme boundary conditions} and Lemma~\ref{lemme skorokhod topology} to get rid of the boundary conditions, we obtain that - uniformly on $\mathsf{Ext}$ being edge-regular:
\begin{multline}
 \Big\vert \Phi^{\otimes r}_{(T_1,\mathscr{X})\rightarrow (T_2,\mathscr{Y})}\left[f^\delta(\mathcal{S}_n) \vert 
            \nidiam(\mathcal{S}), \EXT = \mathsf{Ext}  \right] - \\ \left(\PP^\RW_{(T_1, \mathscr{X})}\right)^{\otimes r} \left[ f^\delta(S)~ \big\vert \nidiam(S), \hit_{(T_2, \mathscr{Y})} \right] \Big\vert  \goes{}{n \rightarrow\infty}{0}
\end{multline}
The main input of Section~\ref{section marches}, namely Theorem~\ref{theoreme invariance principle for drw} then allows us to conclude that 
\begin{equation}
    \Phi^{\otimes r}_{(T_1,\mathscr{X})\rightarrow (T_2,\mathscr{Y})}\left[f^\delta(\mathcal{S}_n) \vert 
            \nidiam(\mathcal{S}), \EXT = \mathsf{Ext}  \right]  \goes{}{n\rightarrow \infty}{\E\left[f^\delta(\sigma\bw^{(r)})\right]},
\end{equation}
for some $\sigma>0$ that depends on the distribution $\mathbf{P}$ in~\eqref{equation OZ boundary conditions}, but of course not on $\delta$. This concludes the proof of point $(i)$ of Lemma~\ref{Lemme technique processus sto}.

The point $(ii)$ - the equicontinuity of the family $\mathcal{S}_n$ at 0 and 1 -  is an easy consequence of classical large deviations estimates combined with the arguments above. 

By Lemma~\ref{Lemme technique processus sto}, we thus proved that $\mathcal{S}_n$ converges in distribution towards $\bw^{(r)}$ when sampled under the distribution $ \Phi^{\otimes r}_{(0,x)\rightarrow (n,y)}\left[f^\delta(\mathcal{S}_n) \vert 
            \mathcal{S} \in \mathcal{W}_{[T_1, T_2]}, \EXT = \mathsf{Ext}  \right]$. Now, observe that the diamond confinement property and the volume estimate stated in Lemma~\ref{corollaire volume diamants} yield that
\begin{equation}
    \Phi^{\otimes r}_{(0,x)\rightarrow (n,y)}\left[ \sup_{0\leq t \leq n}\left| \Gamma^\pm(t) - \mathcal{S}(t) \right| > \log^3 n \right] \leq \exp(-c\log^2 n).
\end{equation}
This concludes the proof, by the usual observation that this decay is faster than any polynomial.
\end{proof}

\section{Brownian watermelon asymptotics for the random-cluster measure}\label{section RCM}

Now that the convergence of the rescaled clusters towards the Brownian watermelon is established in the case of the product measure, the goal of the following section is to transfer this convergence to the rescaled clusters sampled under the "real" random-cluster measure, and thus to achieve our journey towards Theorems~\ref{theoreme estimation} and~\ref{Theoreme main}. The strategy looks similar to the precedent section: indeed, we shall first prove an edge repulsion lemma in Section~\ref{sub edge rep RCM}. Then we shall prove in Section~\ref{sub global rep RCM} that the clusters remain far away from each other in the bulk. This will finally allow us to conclude that in the bulk, the conditioned random-cluster measure is close to the conditioned product measure thanks to a mixing argument, and to import the results of the precedent section to conclude the proofs in Section~\ref{Sub finale mixing et preuves}. The main difficulty and the reason why we needed to introduce and study the measure $\phi^{\otimes r}$ is that a coupling such as $\Phi_{x \rightarrow y}^{\otimes r}$ is not available in this setting. Hence, the random diamond decomposition and its associated skeleton given by the coupling $\Phi_{x \rightarrow y}^{\otimes r}$ do not exist anymore. We then work with the \emph{maximal} diamond decomposition and \emph{maximal} skeletons of the clusters (see remark~\ref{remarque maximal skeleton}). We draw the attention of the reader on the fact that this maximal skeleton \emph{does not} behave like a process with independent increments as it was the case in Section~\ref{section independent system}.

Introduce the following notation: if $\mathcal{E}$ is a set of edges of $\Z^2$ and $\eta, \om$ are two percolation configurations, we set
\begin{equation}
    \lbrace \om \overset{\mathcal{E}}{=} \eta \rbrace = \lbrace \forall e \in \mathcal{E}, \om(e) = \eta(e) \rbrace.
\end{equation}
The first easy comparison between the infinite volume and the product measures is given by the following lemma:
\begin{lemma}\label{lemme estimation fine proba con ni}
Let $\mathcal{E}$ be an arbitrary subset of $E(\Z^2)$ and $\eta$ an arbitrary percolation configuration on $\Z^2$. Then,
\begin{equation}\label{equation lemme estimation fine proba con ni}
    \phi\big[ \con, \nonint \vert \om \overset{\mathcal{E}}{=} \eta  \big] \geq \phi^{\otimes r }\big[ \con, \nonint \vert  \om_1 \overset{\mathcal{E}}{=} \dots \overset{\mathcal{E}}{=} \om_r \overset{\mathcal{E}}{=} \eta  \big].
\end{equation}
\end{lemma}

\begin{proof}

It is a simple consequence of the FKG inequality applied to $\phi[\cdot\vert \om \overset{\mathcal{E}}{=} \eta ]$, which is a random-cluster measure on the graph $\Z^2 \setminus \mathcal{E}$ with some boundary conditions imposed by the configuration $\eta$. If $C$ is a connected set of edges of $\Z^2$, introduce its edge exterior boundary by:
\begin{equation}
    \partial_{\text{ext}} C = \left\lbrace \lbrace x,y\rbrace \in E(\Z^2) \cap C^c, x \text{ is the endpoint of an edge of }C\right\rbrace.
\end{equation}
Then, we write, summing over all the possible realizations of $\mathcal{C}_1, \dots, \mathcal{C}_r$ such that $\con \cap \nonint$ occurs:

\begin{align*}
\phi \big[\nonint, \con\vert \om \overset{\mathcal{E}}{=}\eta \big] 
&= \sum_{C_1, \dots, C_r } \phi\big[\bigcap_{i=1}^r\lbrace \mathcal{C}_i = C_i\rbrace \vert \om \overset{\mathcal{E}}{=}\eta\big]\\
&= \sum_{C_1, \dots, C_r } \phi\big[ \bigcap_{i=1}^r \lbrace C_i\text{ is open },\partial_{\text{ext}}C_i \text{ is closed} \rbrace\vert \om \overset{\mathcal{E}}{=}\eta\big] \\
&=\sum_{C_1, \dots, C_r } \phi\big[ \bigcap_{i=1}^r \lbrace \partial_{\text{ext}}C_i \text{ is closed}\rbrace\vert \om \overset{\mathcal{E}}{=} \eta \big]\prod_{i=1}^r \phi^0_{C_i}\big[C_i \text{ is open } \vert ~\om \overset{\mathcal{E}\cap C_i }{=} \eta \big]\\
&\geq\sum_{C_1, \dots, C_r } \prod_{i=1}^r \phi\big[\partial_{\text{ext}} C_i \text{ is closed }\vert \om \overset{\mathcal{E}}{=} \eta \big]\phi_{C_i}\big[C_i \text{ is open } \vert \om \overset{\mathcal{E}\cap C_i}{=}\eta \big]\\
&= \sum_{C_1, \dots, C_r } \prod_{i=1}^r \phi\big[\mathcal{C}_i = C_i \vert \om \overset{\mathcal{E}}{=} \eta \big]\\
&=\phi^{\otimes r} \big[ \nonint,\con \vert \om \overset{\mathcal{E}}{=} \eta \big],
\end{align*}
where the inequality comes from the positive association property~\eqref{equation FKG} of the measure $\phi\big[\cdot\vert \om \overset{\mathcal{E}}{=}\eta \big]$.

\end{proof}

\begin{remark}\label{remarque minoration ni con arbitrary BC}
The lemma above together with Lemma~\ref{lemme minoration proba con nonint mesure produit} immediately yields that for any $x, y \in W\cap\Z^r$,
\begin{equation}\label{equation minoration proba con nonint arbitrary boundary conditions}
    \phi\left[\nonint, \con  \right] \geq cV(x)V(y)n^{-\frac{r^2}{2}}\e^{-\tau rn}.
\end{equation}
This will be of particular interest later - and is the first half of the proof of Theorem~\ref{theoreme estimation}.
\end{remark}
While the latter bound is optimal (up to a constant), we also import a rough non-optimal upper bound.
\begin{lemma}\label{lemme majoration rough upper bound con ni}
    Let $x, y \in W \cap \Z^2$. Then,
    \begin{equation}
        \phi\left[\con, \nonint \right] \leq \e^{-\tau rn}.
    \end{equation}
\end{lemma}

Before turning to the proof of Lemma~\ref{lemme majoration rough upper bound con ni}, we introduce a useful notation for the rest of the paper. When $x,y \in W\cap\Z^r$, if $1\leq i \leq r$ we write $(\con,\nonint)_{\neq i}$ for the event that $\mathcal{C}_1, \dots, \mathcal{C}_{i-1},\mathcal{C}_{i+1}, \dots, \mathcal{C}_r$ realize the connection event and are non-intersecting. Observe that whenever $1\leq i \leq r$, $(\con, \nonint) \subset (\con,\nonint)_{\neq i}$, while the opposite inclusion is obviously not true.

\begin{proof}[Proof of Lemma~\ref{lemme majoration rough upper bound con ni}]
    We proceed by induction on $r$. For $r = 1$, the statement to prove is 
    \begin{equation}
        \phi\left[(0,x) \leftrightarrow (n,y) \right] \leq \e^{-\tau rn},
    \end{equation}
    which is the consequence of a well-known subbaditivity argument. 

If the statement is established with $r$ clusters, let $x, y \in \Z^{r+1}$. Then, observe that if $\con, \nonint$ occurs, then $(\con, \nonint)_{\neq r+1}$ has to occur. Summing over all the potential realizations of $\mathcal{C}_1, \dots, \mathcal{C}_r$ under $\con, \nonint$, we get:
\begin{align*}
    \phi\left[\con, \nonint\right] &= \sum_{C_1, \dots, C_r}\phi\big[(0,x_{r+1})\overset{(C_1 \sqcup \dots \sqcup C_r)^c}{\longleftrightarrow} (n, y_{r+1}) \vert \mathcal{C}_i = C_i, \forall 1\leq i\leq r \big]\\ & \hspace{196.5pt}\times \phi\left[\mathcal{C}_i = C_i, \forall 1\leq i\leq r  \right] \\
    &= \sum_{C_1, \dots, C_r} \phi^0_{(C_1 \sqcup \dots \sqcup C_r)^c}\left[(0,x_{r+1}) \leftrightarrow (n,y_{r+1})\right]\phi\left[\mathcal{C}_i = C_i, \forall 1\leq i\leq r  \right] \\
    &\leq \phi\left[(0,x_{r+1}) \leftrightarrow (n,y_{r+1})\right]\phi\left[(\nonint, \con)_{\neq r+1}\right] \\
    &\leq \e^{-\tau n}\phi\left[(\nonint, \con)_{\neq r+1}\right],
\end{align*}
where we used~\eqref{equation smp} in the second line,~\eqref{Comparaison boundary conditions} in the third line, and the case $r=1$ in the last line. The statement follows by the induction hypothesis. 
    \end{proof}
    
We next state another consequence of these two bounds, observing that they allow us to derive a diamond confinement property for the infinite volume measure conditioned on $\con \cap \nonint$, the exact analog of Corollary~\ref{corollaire volume diamants} for the conditioned measure. We formulate it for a rather particular class of boundary conditions, in order to be able to apply it later: however the reader should think about the measure $\phi$ on the strip $\mathsf{Strip}_n$ with boundary conditions given by the trace of a subcritical cluster outside of the strip. We recall that $\mathcal{D}^\mathsf{max}(\mathcal{C})$ denotes the \textit{maximal} diamond decomposition of the cluster $\mathcal{C}$, and introduce the following events:

\begin{equation}
   \BigDiam_i = \lbrace \max_{\substack{\mathcal{D}\subset \mathcal{D}^{\mathsf{max}}(\mathcal{C}_i) \\ \mathcal{D}\text{ diamond}}} \Vol(\mathcal{D}) \geq \log^2 n \rbrace ~~~~\text{  and  }~~~~ \BigDiam = \bigcup_{i=1}^r\BigDiam_i.
\end{equation}

\begin{lemma}[Diamond confinement]\label{lemme diamond confinement mesure conditionnee} 
There exists a constant $c> 0$ such that the following occurs. Let $\Ext$ be a finite set of edges such that $\Ext \cap E(\mathsf{Strip}_n)  = \emptyset$. Then for any $n$ large enough,
\begin{equation}
    \phi^0_{\Ext^c} \Big[\BigDiam ~\big\vert \con,\nonint\Big] \leq \exp\left(-c\left(\log n\right)^{2}\right).
\end{equation} 
\end{lemma}

\begin{proof}
We write 
\begin{equation}
 \phi^0_{\Ext^c} \Big[\BigDiam ~\big\vert \con,\nonint\Big] \leq \sum_{i = 1}^r \phi^0_{\Ext^c} \Big[ \BigDiam_i ~\big\vert \con,\nonint\Big].
\end{equation}
We fix an $i \in \lbrace 1, \dots, r \rbrace$. Now we shall focus on the numerator of the latter probability, namely on estimating
the quantity $
    \phi^0_{\Ext^c} \Big[ \BigDiam_i, \con,\nonint\Big].
$
Summing over all the potential clusters $C_1, \dots, C_{i-1}, C_{i+1}, \dots, C_r$ under $\con, \nonint$, 
\begin{multline}\label{equ: equation conditionnement preuve max diamond}
\phi^0_{\Ext^c} [ \BigDiam_i, \con,\nonint] \\ = \sum_{ C_1, \dots, C_{i-1}, C_{i+1}, \dots C_r } \phi^0_{\Ext^c}[ \BigDiam_i, \con,\nonint \vert \mathcal{C}_j = C_j]\phi^0_{\Ext^c}[ \mathcal{C}_j = C_j ],
\end{multline}
where the event in the conditioning is shorthand for $\mathcal{C}_j = C_j$ for all $1\leq j\neq i\leq r$.
Fix such a system of clusters $C_1, \dots, C_{i-1}, C_{i+1}, \dots, C_r$, and call by convenience 
 $   \widetilde{\Ext} = \Ext\cup\left(\bigcup_{1\leq j \neq i \leq r}C_j\cup\partial_{\mathsf{ext}}C_j\right).
$
We then observe that - thanks to~\eqref{equation smp}, 
\begin{equation}
    \phi^0_{\Ext^c} \Big[ \BigDiam_i, \con, \nonint ~ \big\vert \mathcal{C}_j = C_j\Big] = \phi^0_{\widetilde{\Ext}^c} \Big[\BigDiam_i , (n,y_i) \in \mathcal{C}_i \Big].
\end{equation}
Moreover, since the events $\lbrace (n,y_i) \in \mathcal{C}_i \rbrace $ and $\BigDiam_i$ are increasing, we obtain:
\begin{equation}
    \phi^0_{\widetilde{\Ext}^c} \Big[ \BigDiam_i, (n,y_i) \in \mathcal{C}_i \Big] \leq 
     \phi\Big[ \BigDiam_i, (n,y_i) \in \mathcal{C}_i \Big]
\end{equation}
But we are now in the setting of Corollary~\ref{corollaire volume diamants}, which ensures that
\begin{equation}
   \phi\Big[ \BigDiam_i, (n,y_i) \in \mathcal{C}_i \Big] \leq \e^{-\tau n}n^{-c\log n}.
\end{equation}
Indeed, the diamonds appearing in the \textit{maximal} diamond decomposition are always contained in the ones appearing in the diamond decomposition given by the Ornstein--Zernike coupling.
Coming back to~\eqref{equ: equation conditionnement preuve max diamond}, we proved that 
\begin{equation}
    \phi^0_{\Ext^c} \Big[ \BigDiam_i, \con,\nonint\Big] \leq \e^{-\tau n - c\log^2 n}\phi^0_{\Ext^c}\left[ (\con, \nonint)_{\neq i} \right]
\end{equation}
 Using the rough upper bound given by Lemma~\ref{lemme majoration rough upper bound con ni}, we obtain:
\begin{equation}
      \phi^0_{\Ext} \Big[ \BigDiam_i, \con,\nonint\Big] \leq \e^{-(\tau rn + c\log^2 n) }.
\end{equation}
Now, thanks to Remark~\ref{remarque minoration ni con arbitrary BC} and the uniform Ornstein--Zernike decay~\eqref{equation oz uniform BC}, we bound the denominator:
\begin{equation}
    \phi^0_{\Ext^c}\left[\con,\nonint\right] \geq cn^{-\frac{r^2}{2}}\e^{-\tau rn}.
\end{equation}
Combining the two bounds above, we find
\begin{equation}
     \phi^0_{\Ext^c} \Big[\BigDiam_i~\big\vert \con,\nonint\Big] \leq \tfrac1c n^\frac{r^2}{2}\exp(-c (\log n)^2).
\end{equation}
Applying the union bound yields the result for an amended value of $c$.
\end{proof}

\subsection{Edge repulsion}\label{sub edge rep RCM}

The goal of this section is to prove the equivalent of Lemma~\ref{Lemme entropic repulsion independent system} for the measure $\phi[~\cdot~ \vert \con, \nonint]$. We need an alternative definition of the random times $T_1$ and $T_2$, since $\mathcal{S}$ is not available anymore. Recall the definition of the upper and lower interfaces of a cluster $\Gamma^{\pm}(t)$.

\begin{definition}
    Fix $\eps > 0$. We define the two following random variables. 
    \begin{align*}
        T'_1 &= \min\Big\lbrace t \geq 0, \min_{*, \star \in \pm}\min_{1\leq i < j \leq r}\big\vert \Gamma^\star_i(t) - \Gamma^*_j(t) \big\vert > n^\eps  \Big\rbrace
\quad  \text{  and }\\
        T'_2 &= \max\Big\lbrace t \geq 0, \min_{*, \star \in \pm}\min_{1\leq i < j \leq r}\big\vert \Gamma^\star_i(t) - \Gamma^*_j(t) \big\vert > n^\eps  \Big\rbrace.
    \end{align*}
\end{definition}
    
The analogous of Lemma~\ref{Lemme entropic repulsion independent system} is the following:

\begin{lemma}[Edge repulsion]\label{Entropic repulsion principle} There exists $\eps > 0$ and $c>0$ such that for any $n$ large enough,
\begin{equation}\label{equation entropic repulsion}
    \phi \left[ \left\lbrace T'_1 > n^{1-\eps}\right\rbrace \cup \left\lbrace T'_2 < n-n^{1-\eps}\right\rbrace \vert \con, \nonint \right] < 2\exp\left(-cn^{1-3\eps} \right).
\end{equation}
\end{lemma}

The value of $\eps> 0$ given by Lemma~\ref{Entropic repulsion principle} will be fixed in the rest of the paper.

\begin{proof}

Let $\eps> 0$, its value will be determined at the end of the proof. By symmetry, we focus on proving the following bound
\begin{equation}
    \phi\left[T'_1 > n^{1-\eps} \vert \con,\nonint \right] \leq \exp(-cn^{1-3\eps}). 
\end{equation}
As in the proof of Lemma~\ref{Lemme entropic repulsion independent system}, we will conclude by time reversal and an easy union bound. For $2 \leq i \leq r$, we define the event $\mathsf{MLCP}_i$ (meaning "many left-close points") by 
\begin{equation}
    \mathsf{MLCP}_i = \lbrace \# \lbrace k \in \lbrace 0, \dots, n^{1-\eps}\rbrace, |\Gamma^{-}_i(k) - \Gamma^+_{i-1}(k)| < n^\eps \rbrace \geq \frac{1}{r} n^{1-\eps} \rbrace. 
\end{equation}
The reason for the introduction of this event is the following inclusion (that is a simple consequence of the pigeonhole principle):
\begin{equation}
    \lbrace T'_1 > n^{1-\eps }\rbrace \subset \bigcup_{i=2}^r \mathsf{MLCP}_i.
\end{equation}
Thus, by union bound
\begin{equation}
    \phi\left[T'_1 > n^{1-\eps} \vert \con, \nonint\right] \leq \sum_{i=2}^r \phi\left[\mathsf{MLCP}_i\vert\con,\nonint\right].
\end{equation}
Fix some $i \in \lbrace 2, \dots, r\rbrace$. We upper bound $\phi\left[\mathsf{MLCP}_i\vert\con,\nonint\right]$ by separately bounding the numerator and the denominator of this fraction. We start with the numerator, and we write, conditioning over all the possible clusters $C_1, \dots, C_{i-1}, C_{i+1}, \dots, C_r$ under $\con, \nonint$:
\begin{equation}\label{eq: equation conditionnement clusters repulsion edge}
    \phi\left[\mathsf{MLCP}_i,\con,\nonint\right] =  \sum_{C_1, \dots, C_{i-1}, C_{i+1}, \dots, C_r}\phi\left[\mathsf{MLCP}_i, (n,y_i)\in\mathcal{C}_i \vert \mathcal{C}_j=C_j\right]\phi\left[\mathcal{C}_j = C_j\right],
\end{equation}
 where the event in the conditioning is shorthand for $\mathcal{C}_j = C_j$ for all $1\leq j\neq i\leq r$.
 Let us fix $C_1, \dots, C_{i-1}, C_{i+1}, \dots, C_r$ that can appear in the sum ~\eqref{eq: equation conditionnement clusters repulsion edge}.
As in the precedent proof, we define the following set of edges
\begin{equation}
    \widetilde{\Ext} = \bigcup_{1\leq j \neq i \leq r} (C_j \cup \partial _{\text{ext}}C_j).
\end{equation}
Following the previous computation and using the spatial Markov property~\eqref{equation smp}, we observe that 
\begin{equation}
    \phi\left[\mathsf{MLCP}_i, (n,y_i)\in\mathcal{C}_i \vert \mathcal{C}_j=C_j\right] = \phi^0_{\widetilde{\Ext}^c}\left[ \mathsf{MLCP}_i, (n,y_i)\in \mathcal{C}_i\right],
\end{equation}
where
\begin{multline*}
    \phi^0_{\widetilde{\Ext}^c}\left[ \mathsf{MLCP}_i,  (n,y_i)\in\mathcal{C}_i\right] =  \\ \phi^0_{\widetilde{\Ext}^c}\left[(n,y_i)\in\mathcal{C}_i, \# \lbrace k \in \lbrace 0, \dots, n^{1-\eps}\rbrace, |\Gamma^{-}_i(k) - \Gamma^+(C_{i-1})(k)| < n^\eps \rbrace \geq \tfrac{1}{r} n^{1-\eps}\right].
\end{multline*}
The event appearing on the right-hand side of the latter equation is increasing (the connection event is always increasing, and adding edges to the configuration can only push $\Gamma^-_i$  down, rendering it closer to $\Gamma^+(C_{i-1})$). Thus, by ~\eqref{Comparaison boundary conditions}, we obtain:
\begin{multline}\label{equ: equation remarque increasing event}
    \phi\left[\mathsf{MLCP}_i, (n,y_i)\in\mathcal{C}_i \vert \mathcal{C}_j=C_j\right] \leq \\ \phi\left[  (n,y_i)\in\mathcal{C}_i, \# \lbrace k \in \lbrace 0, \dots, n^{1-\eps}\rbrace, |\Gamma^{-}_i(k) - \Gamma^+(C_{i-1})(k)| < n^\eps \rbrace \geq \frac{1}{r} n^{1-\eps} \right]. 
\end{multline}
We are now in the framework of the classical one-cluster Ornstein--Zernike theory - and we are going to conclude using Lemma~\ref{Confinement lemma}. Indeed, observe that if $k \in \{0, \dots, n^{1-\eps}\}$ satisfies $|\Gamma^{-}_i(k) - \Gamma^+(C_{i-1})(k)| <n^\eps$, then one can invoke Lemma~\ref{lemme diamond confinement mesure conditionnee} applied to both $\mathcal{C}_i$ and to $C_{i-1}$ to argue that with very large probability, there exists some time $t_k$ (measurable with respect to $\mathcal{C}_i$: one can choose the next renewal after $k$) such that $t_k$ is a renewal time for $\mathcal{C}_i$ and $|\mathcal{S}(t_k) - \Gamma^+(C_{i-1})(t_k)| < 2n^\eps$. The map $k \mapsto t_k$ is not one-to-one but still due to Lemma~\ref{lemme diamond confinement mesure conditionnee} can be at most $\log^2 n$-to-1 with high probability. Thus we get that for some $\alpha < 1$, 
\begin{multline*}
    \phi\left[  (n,y_i)\in\mathcal{C}_i, \# \lbrace k \in \lbrace 0, \dots, n^{1-\eps}\rbrace, |\Gamma^{-}_i(k) - \Gamma^+(C_{i-1})(k)| < n^\eps \rbrace \geq \frac{1}{r} n^{1-\eps} \right] \\
    \leq \e^{-\tau n}\Phi_{(0,x_i)\leftrightarrow (n,y_i)}\left[\# \lbrace 0 \leq k \leq n^{1-\eps}, |\mathcal{S}^i(k) - \Gamma^+(C_{i-1})(k)|<2n^\eps \rbrace \geq \frac{\alpha}{r}n^{1-\eps} \right],
\end{multline*}
In words, the trajectory of $\mathcal{S}$ stays confined during a large amount of time close to the function $\Gamma^+(C_{i-1})$.
Since $\Gamma^+(C_{i-1})$ is the graph of a fixed function we can apply Lemma~\ref{Confinement lemma} \footnote{Lemma~\ref{Confinement lemma} is stated for unconditioned random walks, but as usual due to the Local Limit Theorem, being a bridge has a polynomial probability which is always beaten by the quantity $\exp(-cn^{1-3\eps})$} to argue that
\begin{multline*}
    \Phi_{(0,x_i)\leftrightarrow (n,y_i)}\left[\# \lbrace k \in \lbrace 0, \dots, n^{1-\eps}\rbrace, |\mathcal{S}^i(k) - \Gamma^+(C_{i-1})(k)|<n^\eps \rbrace \geq \frac{1}{r}n^{1-\eps} \right] \\  \leq \exp(-cn^{1-3\eps}),
\end{multline*}
provided that $\eps >0 $ is small enough. Coming back to ~\eqref{eq: equation conditionnement clusters repulsion edge}, we proved that 
\begin{equation}
    \phi\left[ \mathsf{MLCP^i(r)},\con,\nonint\right] \leq \exp(-cn^{1-3\eps} - \tau n)\phi\left[ (\con, \nonint)_{\neq i}\right]. 
\end{equation}
Finally using the rough bound given by Lemma ~\ref{lemme majoration rough upper bound con ni} we proved that
\begin{equation}
    \phi\left[T'_1 > n^{1-\eps}, \con, \nonint \right] \leq \exp(-(\tau rn + cn^{1-3\eps})).
\end{equation}
We conclude by using the lower bound on $\phi\left[\con, \nonint \right]$ given by ~\eqref{equation minoration proba con nonint arbitrary boundary conditions}. We obtain
\begin{equation}
     \phi\left[T'_1 > n^{1-\eps}\vert \con, \nonint \right] \leq n^{\frac{r^2}{2}}\exp(-cn^{1-3\eps}),
\end{equation}
which yields the result, up to slightly changing the value of the constant $c$.

\end{proof}

Observe that, by their definition, $T'_1$ and $T'_2$ are not necessarily synchronization times. 
We chose to define them as such so as to obtain an increasing event in~\eqref{equ: equation remarque increasing event}; otherwise the proof of Lemma~\ref{Entropic repulsion principle} would be more complicated.
Thus let us define the actual random variables $T_1$ and $T_2$ by the following:
\begin{align*}
    T_1 &= \min\lbrace k \geq T'_1, k \text{ is a synchronization time for }\mathcal{S}^{\mathsf{max}} \rbrace 
\quad \text{ and } \\    
T_2 &= \max\lbrace k \leq T'_2, k \text{ is a synchronization time for }\mathcal{S}^{\mathsf{max}} \rbrace.
\end{align*}

Similarly to Section~\ref{sub edge repulsion}, Lemma~\ref{Entropic repulsion principle} will be used to produce edge-regular configurations with large probability. We recall and modify slightly the notion of edge-regular configurations: a percolation configuration $\om \in \con,\nonint$ is called \emph{edge-regular} (also written $\om\in \mathsf{EdgeReg}$) if the following set of conditions is satisfied: 
\begin{itemize}
    \item $T_1 < n^{1-\eps}$ and $T_2 > n-n^{1-\eps}$.
    \item $\norme{\mathscr{X}}\leq n^{1/2- \eps/4}$ and $\norme{\mathscr{Y}} \leq n^{1/2 - \eps/4}$.
    \item $\Gap(\mathscr{X}) > \frac{1}{2}n^{\eps}$ and $\Gap(\mathscr{Y}) > \frac{1}{2}n^\eps$.
\end{itemize}
Then typical configurations are edge-regular under the conditioning on $\con, \nonint$.
\begin{lemma}\label{lem: edge regularity true measure}
   There exists a constant $c>0$ such that 
   \begin{equation}
       \phi\left[\mathsf{EdgeReg}^c\vert\con,\nonint\right] \leq \frac{1}{c}\exp(-cn^{\eps/2}).
   \end{equation}
\end{lemma}
We only briefly sketch the proof since it is very similar to the one of Lemma~\ref{lem: edge regularity product measure}.
\begin{proof}
By Lemma ~\ref{lemme diamond confinement mesure conditionnee}, it is easy to see that:
\begin{equation}\label{equ: modified T_1 close to T_1}
    \phi\left[ \max\lbrace |T_1-T'_1|,|T_2-T'_2| \rbrace > (\log n)^3 \vert \con,\nonint\right] \leq \e^{-c(\log n)^2}
\end{equation}
for some small constant $c>0$. Indeed, one can use Lemma ~\ref{lemme diamond confinement mesure conditionnee} together with the fact that being a cone-point is a decreasing event. This fact established, the proof is essentially the same as the proof of Lemma~\ref{lem: edge regularity product measure}: in a time smaller than $(\log n)^2$, the clusters cannot move to a polynomial distance of their starting point by a basic large deviations estimate. 
\end{proof}
\subsection{Global repulsion }\label{sub global rep RCM}

In this section, we work under the measure $\phi\left[~\cdot\vert \con, \nonint\right]$, and we want to prove that between a time $o(n)$ and $n-o(n)$, the minimal gap between the clusters diverges with $n$. The event that we are going to estimate is the following "global repulsion" event: 
\begin{align*}
    \mathsf{GlobRep} := &\lbrace T_1 < n^{1-\eps} \rbrace \cap \lbrace T_2 > n-n^{1-\eps} \rbrace 
    \cap \Big\lbrace \min_{\substack{1< i \leq r \\t \in [T_1,T_2]}} \big\vert \Gamma^-_{i}(t) - \Gamma^+_{i-1}(t) \big\vert > (\log n)^2\Big\rbrace. 
\end{align*}

The goal of this section is to prove the following statement.

\begin{proposition}[Global repulsion estimate]\label{proposition global entropic repulsion} 
    There exists $\beta > 0$, depending only on $r$, and $c>0$ such that 
    \begin{equation}\label{entropic repulsion equation}
        \phi \left[\mathsf{GlobRep} \vert \con,\nonint \right] \geq 1-cn^{-\beta}.
    \end{equation}
\end{proposition}

This lemma will be the main ingredient for the proofs of Theorems~\ref{theoreme estimation} and~\ref{Theoreme main}. 
The rest of the section is dedicated to the proof of Proposition~\ref{proposition global entropic repulsion}. We first use Lemma~\ref{lem: edge regularity true measure} to write:

\begin{eqnarray}\label{eq:globrep_decomp}
    \phi\left[\mathsf{GlobRep}^c\vert \nonint, \con\right] 
    &\leq& \phi[\mathsf{EdgeReg}^c \vert\nonint, \con] + \phi[\mathsf{GlobRep}^c \vert \nonint, \con, \mathsf{EdgeReg} ]\nonumber\\
    &\leq& \tfrac{1}{c}\exp\big(-cn^{\eps/2}\big) + \phi[\mathsf{GlobRep}^c \vert \nonint, \con, \mathsf{EdgeReg}].
\end{eqnarray}
We will focus on bounding the second term in the right-hand side of the above. We will do so by conditioning on $T_1$, $T_2$ and the shape of the clusters before $T_1$ and after $T_2$.

As previously, write $\EXT$ for the trace of the clusters and their boundaries outside of the strip $\mathsf{Strip} = [T_1,T_2] \times \mathbb Z$. Recall the notation $\mathscr{X}$ and $\mathscr{Y}$ for the vertical positions of the renewals at the times $T_1$ and $T_2$. 
Then $\phi[\cdot|\con, \nonint, \EXT]$\footnote{The conditioning on $\EXT$ contains implicitly the fact that $T_1$ and $T_2$ are renewals.} on the complement of $\EXT$ is the measure  $\phi^0_{ \EXT^c}$ conditioned on $\mathscr{X}$ being connected to $\mathscr{Y}$ by disjoint clusters contained in the respective diamonds. Write $\mathsf{Diam}_i$ for the fact that the $i$-th connection above occurs indeed in the corresponding diamond, and set $\mathsf{Diam} = \bigcap_{1 \leq i \leq r}\mathsf{Diam}_i$. Then 
\begin{align} \label{equation conditionnement trajext preuve globrep}
    \phi\big[&\mathsf{GlobRep}^c\vert\nonint, \con, \mathsf{EdgeReg}\big] \nonumber
    \\ &= \sum_{\mathsf{Ext}} \phi \left[ \mathsf{GlobRep}^c \big\vert \nonint, \con,\EXT = \mathsf{Ext} \right]\phi\left[\EXT=\mathsf{Ext} ~\big\vert\con, \nonint,\mathsf{EdgeReg}\right]\nonumber\\
    &=\sum_{\mathsf{Ext}} \frac{\phi \left[ \mathsf{GlobRep}^c,\nonint, \con|\EXT = \mathsf{Ext} \right]}
    {\phi \left[ \nonint, \con|\EXT = \mathsf{Ext} \right]}
    \phi\left[\EXT=\mathsf{Ext} ~\big\vert\con, \nonint,\mathsf{EdgeReg}\right],
\end{align}
where the sum runs over all possible edge-regular realizations $\mathsf{Ext}$ of $\EXT$.
We are going to focus on bounding the ratio above uniformly over $\Ext$. 
This ratio may be written as 
\begin{align}\label{equation mesure conditionnement trajext preuve globrep}
    \frac{\phi \left[ \mathsf{GlobRep}^c,\nonint, \con|\EXT = \mathsf{Ext} \right]}
    {\phi \left[ \nonint, \con|\EXT = \mathsf{Ext} \right]}=
    \frac{\phi_{\Ext^c}^0 \left[ \mathsf{GlobRep}^c,\nonint, \con, \mathsf{Diam}\right]}
    {\phi_{\Ext^c}^0 \left[\nonint, \con, \mathsf{Diam}\right]}.
\end{align}

\begin{figure}
\begin{center}
    \includegraphics[width = .75\textwidth, page = 7]{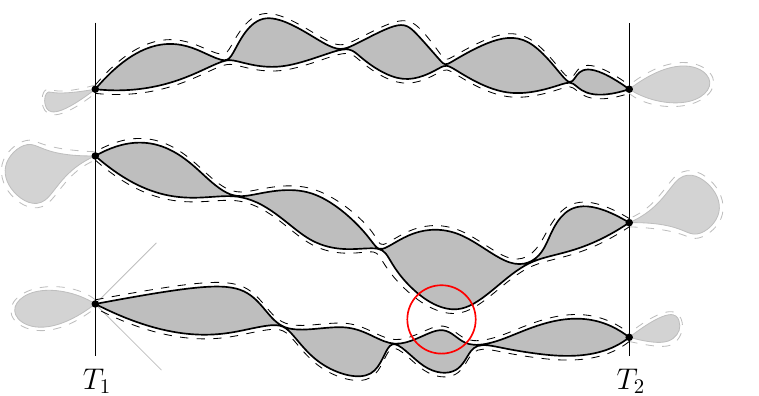}
    \caption{The region $\mathsf{Strip}$ is the vertical strip between $T_1$ and $T_2$; $\EXT$ is the configuration outside of this strip. 
    The events $\nonint$, $\con$ and $\mathsf{Diam}$ are realised by the three traversing grey clusters. 
    Here $\mathsf{GlobRep}$ is violated since the bottom two clusters come close to each-other in the marked region. }
\end{center}
\end{figure}

Lemma~\ref{lemme estimation fine proba con ni} and Corollary~\ref{corollaire majoration proba ni con mesure produit BC avec ecartement suffisant}  allow us to lower bound the denominator as
\begin{equation}\label{eq: lower bound denominator preuve globrep}
    {\phi_{\Ext^c}^0 \left[\nonint, \con, \mathsf{Diam}\right]}
    \geq \tfrac{1}{\chi} V(\mathscr{X})V(\mathscr{Y}) (T_2 - T_1)^{-\tfrac{r^2}{2}} e^{-\tau r (T_2 - T_1)},
\end{equation}
where $\chi$ and $V(\cdot)$ were described in the aforementioned lemmas. We claim the following upper bound on the numerator.

\begin{lemma}\label{lem:Ioan_numerator}
    There exist constants $\beta > 0$ and $C > 0$ such that, for any edge-regular $\mathsf{Ext}$
    \begin{equation}\label{eq:Ioan_numerator}
    \phi_{\Ext^c}^0 \left[ \mathsf{GlobRep}^c,\nonint, \con, \mathsf{Diam}\right] 
    \leq C V(\mathscr{X})V(\mathscr{Y}) (T_2 - T_1)^{-r^2/2 - \beta} e^{-\tau r (T_2 - T_1)}.
    \end{equation}
\end{lemma}

The lemma above is the main difficulty in the proof of Proposition~\ref{proposition global entropic repulsion}; 
we postpone its proof and finish that of the proposition.

\begin{proof}[Proof of Proposition~\ref{proposition global entropic repulsion}]
    Lemma~\ref{lem:Ioan_numerator} together with the estimate~\eqref{eq: lower bound denominator preuve globrep} yield that for any edge-regular $\mathsf{Ext}$, the following holds:
    \begin{equation}
        \frac{\phi_{\Ext^c}^0 \left[ \mathsf{GlobRep}^c,\nonint, \con, \mathsf{Diam}\right]}
    {\phi_{\Ext^c}^0 \left[\nonint, \con, \mathsf{Diam}\right]} \leq \tfrac{C}{\chi}(T_2-T_1)^{-\beta}.
    \end{equation}
By edge-regularity of $\Ext$, we know that $(T_2-T_1) \geq n-2n^{1-\eps}$, so that $(T_2-T_1)^{-\beta} = n^{-\beta}(1+o(1))$. Thus, inserting this estimate into~\eqref{equation conditionnement trajext preuve globrep}, we proved that
\begin{equation}
    \phi\left[ \mathsf{GlobRep}^c\vert\nonint, \con, \mathsf{EdgeReg} \right] \leq \tfrac{C}{\chi}n^{-\beta}.
\end{equation}
Proposition~\ref{proposition global entropic repulsion} is obtained by applying equations \eqref{eq:globrep_decomp}, \eqref{equation conditionnement trajext preuve globrep} and \eqref{equation mesure conditionnement trajext preuve globrep}.
\end{proof}

We now turn to the proof of Lemma~\ref{lem:Ioan_numerator}.
Fix some edge-regular $\mathsf{Ext}$.
When $\nonint, \con$ and $\mathsf{Diam}$ occur, write $\Gamma_i$ for the top-most path of the cluster of $\mathscr{X}_i$. 
This is a non-simple path of open edges contained in $\mathsf{Strip}$ (due to $\mathsf{Diam}$) connecting $\mathscr{X}_i$ to $\mathscr{Y}_i$. For any such path, write $\partial \Gamma_i$ for all the edges of $\mathsf{Strip}$ adjacent to and above $\Gamma_i$. 
Finally, write $\Gamma$ for the $r$-tuple of paths $\Gamma_1,\dots, \Gamma_r$. 

The idea of the proof of Lemma~\ref{lem:Ioan_numerator} is to bound the probability of $\lbrace \Gamma = \gamma \rbrace$ for any possible realization $\gamma$ of $\Gamma$ by the probability of a suitable event in the product measure.
Then we use the Ornstein--Zernike theory to estimate the latter probability.
The first step is contained by the following statement, which constitutes the core of the proof.

\begin{lemma}\label{lem: first step Ioan numerator} There exists $\eps > 0$ such that for any possible realization $\gamma$ of $\Gamma$,
     \begin{equation}
         \phi_{\mathsf{Ext}^c}^0 \left[ \Gamma = \gamma ,\nonint, \con, \mathsf{Diam}\right] 
      \leq (1-\exp(-n^\eps))^{-r}  (\phi_{\mathsf{Ext}^c}^0)^{\otimes r} \Big[ \bigcap_{i=1}^r\lbrace \gamma_i \text{  is open}, \partial \gamma_i \text{ is closed}\rbrace\Big ].
     \end{equation}
\end{lemma}
\begin{remark}
    The careful reader might think that this lemma is in contradiction with the lower bound given by Lemma~\ref{lemme estimation fine proba con ni}. However, observe that while the events $\mathsf{Diam}$ and $\lbrace \Gamma = \gamma \rbrace$ do imply the events $\con$ and $\nonint$ in the measure $\phi^0_{\mathsf{Ext}^c}$, it is not the case for the measure $(\phi^0_{\EXT^c})^{\otimes r}$.
\end{remark}

\begin{proof}[Proof of Lemma~\ref{lem: first step Ioan numerator}]
  Fix $\gamma$ as in the statement. The paths $\gamma_1,\dots,\gamma_r$ all cross $\mathsf{Strip}$ horizontally, are disjoint and are in increasing vertical order. The same holds for their upper boundaries $\partial \gamma_1,\dots, \partial \gamma_r$. 

    Moreover, 
    \begin{eqnarray*}
        \phi_{\mathsf{Ext}^c}^0 \left[ \Gamma = \gamma ,\nonint, \con, \mathsf{Diam}\right] 
        &=& \prod_{i=1}^r \phi^0_{\mathsf{Ext}^c} \Big[ \Gamma_i = \gamma_i, \mathsf{Diam_i} \Big\vert  \bigcap_{k=1}^{i-1} \lbrace \Gamma_k= \gamma_k, \mathsf{Diam}_k \rbrace\Big]  \\
        &\leq& \prod_{i=1}^r \phi^0_{\mathsf{Ext}^c} \Big[ \gamma_i \text{ is open}, \partial \gamma_i \text{ is closed} \big\vert  \bigcap_{k=1}^{i-1} \lbrace \Gamma_k= \gamma_k, \mathsf{Diam}_k \rbrace\Big].
\end{eqnarray*}

Our goal is now to prove that for any $i \in \lbrace 1, \dots, r \rbrace$,
\begin{multline}\label{equ: equation a prouver majoration chemin par mesure produit}
    \phi^0_{\mathsf{Ext}^c} \Big[ \gamma_i \text{ is open}, \partial \gamma_i \text{ is closed} \Big\vert  \bigcap_{k=1}^{i-1} \lbrace \Gamma_k= \gamma_k, \mathsf{Diam}_k \rbrace\Big]
     \\ 
    \leq (1-\exp(-n^\eps))^{-1}\phi^0_{\mathsf{Ext}^c}\left[ \gamma_i \text{ is open}, \partial \gamma_i \text{ is closed} \right].
\end{multline}
Indeed, assuming that the above inequality is true, we will have proved that
\begin{align*}
    \phi_{\mathsf{Ext}^c}^0 \big[ \Gamma = \gamma ,\nonint, \con, \mathsf{Diam}\big] 
    &\leq \prod_{i=1}^r (1-\exp(-n^\eps))^{-1} \phi^0_{\mathsf{Ext}}\big[ \gamma_i \text{ is open}, \partial \gamma_i \text{ is closed} \big]\\
    &= (1-\exp(-n^\eps))^{-r}(\phi_{\mathsf{Ext}^c}^0)^{\otimes r} \Big[ \bigcap_{i=1}^r\lbrace \gamma_i \text{  is open}, \partial \gamma_i \text{ is closed}\}\Big].
\end{align*}

We thus focus on~\eqref{equ: equation a prouver majoration chemin par mesure produit}.
The bound is obviously true for $i=1$. Fix next $i >1$. Write:
\begin{align}\label{equ: decomposition chemin ouvert et ouvert sachant fermé}
  &  \phi^0_{\mathsf{Ext}^c} \Big[ \gamma_i \text{ is open}, \partial \gamma_i \text{ is closed} \Big\vert  \bigcap_{k=1}^{i-1} \lbrace \Gamma_k= \gamma_k, \mathsf{Diam}_k \rbrace\Big] \\ 
   & = \phi^0_{\mathsf{Ext}^c} \Big[ \gamma_i \text{ open}\Big\vert  \bigcap_{k=1}^{i-1} \lbrace \Gamma_k= \gamma_k, \mathsf{Diam}_k \rbrace\Big] 
    \phi^0_{\mathsf{Ext}^c} \Big[ \partial \gamma_i \text{ closed} \Big\vert  \gamma_i \text{ open}, \bigcap_{k=1}^{i-1} \lbrace \Gamma_k= \gamma_k, \mathsf{Diam}_k \rbrace\Big].\nonumber
\end{align}

The first factor is easy to upper bound, as the conditioning decreases the probability  for $\gamma_i$ to be open. 
Indeed, explore the clusters of $\mathscr{X}_1, \dots, \mathscr{X}_{i-1}$ together with their boundaries, and call $\mathsf{Expl}$ the set of explored edges. Because of the conditioning on $\Diam_1, \dots, \Diam_{i-1}$ and of the disjointness of the paths of $\gamma$,  $\gamma_i$ is disjoint from the explored edges $\mathsf{Expl}\cup\Ext$. 
Furthermore, the measure induced on the complement of these clusters by this exploration procedure is $\phi^0_{(\mathsf{Expl}\cup\mathsf{Ext})^c}$. 
By~\eqref{Comparaison boundary conditions}, 
\begin{equation}
    \phi_{(\mathsf{Ext}\cup\mathsf{Expl})^c}^0\left[\gamma_i \text{ is open} \right] \leq \phi_{\mathsf{Ext}^c}^0\left[\gamma_i \text{ is open} \right],
\end{equation}
which in turn implies
\begin{equation}\label{equ: majoration premier terme proba conditionnelle chemin ouvert}
    \phi^0_{\mathsf{Ext}^c} \Big[ \gamma_i \text{ is open}\Big\vert  \bigcap_{k=1}^{i-1} \lbrace \Gamma_k= \gamma_k, \mathsf{Diam}_k \rbrace\Big] 
    \leq \phi_{\mathsf{Ext}^c}^0\left[\gamma_i \text{ is open} \right].
\end{equation}

\begin{figure}
\begin{center}
    \includegraphics[width = .48\textwidth, page = 8]{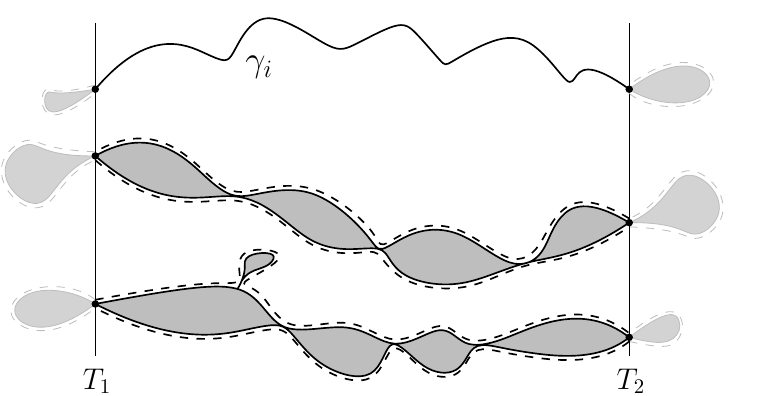}\hspace{.02\textwidth}
        \includegraphics[width = .48\textwidth, page = 9]{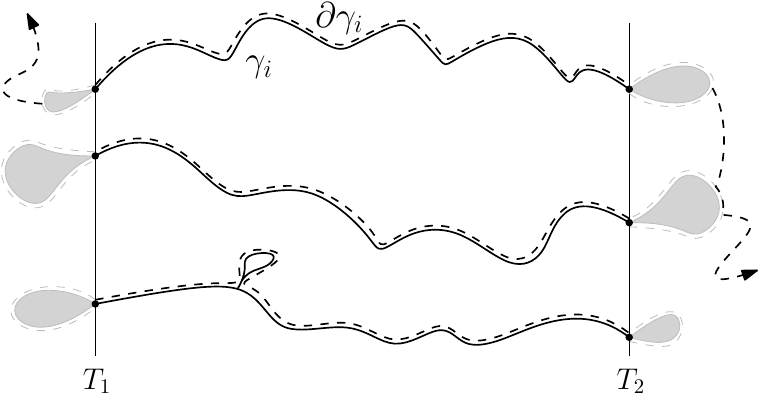}    \caption{{\em Left:} To upperbound $\phi^0_{\mathsf{Ext}^c} [ \gamma_i \text{ open}\vert  \bigcap_{k=1}^{i-1} \lbrace \Gamma_k= \gamma_k, \mathsf{Diam}_k \rbrace]$ it suffices to explore the clusters of $\Gamma_k$ for $k < i$ and observe that any such instance induces negative information on the rest of the space. 
    {\em Right:} When bounding $\phi^0_{\mathsf{Ext}^c} [ \partial \gamma_i \text{ closed} \vert  \gamma_i \text{ open}, \bigcap_{k=1}^{i-1} \lbrace \Gamma_k= \gamma_k, \mathsf{Diam}_k \rbrace]$, the 
    conditioning may increase the probability for $\partial \gamma_i$ to be closed, but not more so than the occurrence of 
    $\mathcal{H}^{\mathsf{L}} \cap \mathcal{H}^{\mathsf{R}}$. The latter events are ensured by the existence of the infinite dual paths on the left and right of the strip.}
\end{center}
\end{figure} 

We turn to the second factor in the right-hand side of \eqref{equ: decomposition chemin ouvert et ouvert sachant fermé}.
Upper bounding this term is slightly more subtle, since the boundary conditions induced by the conditioning may \emph{a priori} help $\partial \gamma_i$ to be closed. 
Introduce the following events 
\begin{align*}
    \mathcal{H}^{\mathsf{L}} &= \big\lbrace  (T_1, \mathscr{Y}_i)^*  \text{ is connected to $\infty$ by a dual open path lying in $\mathsf{Strip}^{c}$}\big\rbrace 
\quad \text{ and}\\
    \mathcal{H}^{\mathsf{R}} &= \big\lbrace  (T_2, \mathscr{Y}_i)^*  \text{ is connected to $\infty$ by a dual open path lying in $\mathsf{Strip}^{c}$}\big\rbrace.
\end{align*}
We now claim that 
\begin{align}\label{eq:HHH}
    \phi^0_{\mathsf{Ext}^c} \Big[ \partial \gamma_i \text{ closed} \Big\vert  
    \gamma_i \text{ open}, \bigcap_{k=1}^{i-1} \lbrace \Gamma_k= \gamma_k, \mathsf{Diam}_k \rbrace\Big]
    \leq \phi_{\mathsf{Ext}^c}^0\big[ \partial \gamma_i \text{ closed} \big\vert  \mathcal{H}^\mathsf{L}\cap\mathcal{H}^\mathsf{R} \cap \{\gamma_i \text{ open}\} \big].
\end{align}
Indeed, the conditioning on $ \bigcap_{k=1}^{i-1} \lbrace \Gamma_k= \gamma_k, \mathsf{Diam}_k\}$ may induce negative information, which improves the probability of $\partial \gamma_i$ to be closed. Nevertheless, this influence is weaker than that of the decreasing events $\mathcal{H}^\mathsf{L}$ and $\mathcal{H}^\mathsf{R}$. Indeed, let us assume that $\mathcal{H}^\mathsf{L}$ and $\mathcal{H}^\mathsf{R}$ are realized, and condition on $\Gamma^\mathsf{L}$ (resp. $\Gamma^\mathsf{R}$) the lowest path producing the event $\mathcal{H}^\mathsf{L}$ (resp. $\mathcal{H}^\mathsf{R}$). \footnote{Such a path is formally explorable by exploring the bottom interface of the dual cluster of $(T_1, \mathscr{Y}_i)$ in $\mathsf{Strip}^c$.}
Now, we sum over all the possible realizations of $\Gamma^\mathsf{L}$ and $\Gamma^\mathsf{R}$ as to obtain
\begin{multline}
\phi_{\mathsf{Ext}^c}^0\big[ \partial \gamma_i \text{ closed} \big\vert  \mathcal{H}^\mathsf{L}\cap\mathcal{H}^\mathsf{R} \cap \{\gamma_i \text{ open}\} \big]  =  \\\sum_{\gamma^\mathsf{L}, \gamma^{\mathsf{R}}}\phi_{\mathsf{Ext}^c}^0\big[ \partial \gamma_i \text{ closed} \big\vert \{ \Gamma^{\mathsf{L}} = \gamma^\mathsf{L} \} \cap  \{\Gamma^{\mathsf{R}} = \gamma^\mathsf{R} \} \cap \{\gamma_i \text{ open}\} \big]\\ \times \phi[\{ \Gamma^{\mathsf{L}} = \gamma^\mathsf{L} \} \cap  \{\Gamma^{\mathsf{R}} = \gamma^\mathsf{R} \} \vert \mathcal{H}^\mathsf{L}\cap\mathcal{H}^\mathsf{R} \cap \{\gamma_i \text{ open}\} ].
\end{multline}
Now we observe that the set $\partial (\Ext^i \cap \gamma^\mathsf{L}\cap\gamma^\mathsf{R} \cap \gamma^i)$ separates $\Z^2$ in two subdomains, one containing $\partial \gamma_i $, and the other one containing all the $\gamma_k$ and their top boundaries.
It is now a classical instance of ``pushing away the boundary conditions'' (or~\eqref{Comparaison boundary conditions}) that for any $\gamma^\mathsf{L}, \gamma^\mathsf{R}$ as above, since $\Gamma^\mathsf{L}$ and $\Gamma^{\mathsf{R}}$ induce free boundary conditions on $\gamma^{\mathsf{L}}$ and $\gamma^{\mathsf{R}}$, and $\{ \partial \gamma^i \text{ closed}\}$ is a decreasing event, 
 \begin{multline}
 \phi_{\mathsf{Ext}^c}^0\big[ \partial \gamma_i \text{ closed} \big\vert \{ \Gamma^{\mathsf{L}} = \gamma^\mathsf{L} \} \cap  \{\Gamma^{\mathsf{R}} = \gamma^\mathsf{R} \} \cap \{\gamma_i \text{ open}\} \big] \\ \geq \phi^0_{\mathsf{Ext}^c} \Big[ \partial \gamma_i \text{ closed} \Big\vert  
    \gamma_i \text{ open}, \bigcap_{k=1}^{i-1} \lbrace \Gamma_k= \gamma_k, \mathsf{Diam}_k \rbrace\Big],
\end{multline}
which by summing over $\gamma^{\mathsf{L}}, \gamma^\mathsf{R}$ concludes the proof of~\eqref{eq:HHH}.

The following claim will be particularly convenient to bound the right-hand side of~\eqref{eq:HHH}. 

\begin{claim}\label{claim 2}
    There exist constants $\eps > 0$ and $c>0$ such that
    \begin{equation}
    \phi^{0}_{\mathsf{Ext}^c}\big[ \mathcal{H}^\mathsf{L}\cap\mathcal{H}^\mathsf{R}\big\vert\gamma_i\text{ is open}  \big] \geq 1 - \exp\left(-cn^\eps\right).
    \end{equation}
\end{claim}

\begin{proof}[Proof of Claim~\ref{claim 2}] 

The proof is a standard argument using the properties of the subcritical regime. 
We will prove that 
\begin{align}\label{eq:HL}
    \phi^0_{\mathsf{Ext}^c} \big[ \mathcal{H}^\mathsf{L} \big\vert\gamma_i\text{ is open}  \big] 
    \geq 1 - \tfrac12 \exp\left(-cn^\eps\right),
\end{align} 
and the claim will follow by the union bound. 

The event $\mathcal{H}^\mathsf{L}$ is decreasing and the conditioning only depends on the configuration in ${\sf Strip}$. Thus~\eqref{Comparaison boundary conditions} implies 
\begin{align*}
    \phi^0_{\mathsf{Ext}^c} \big[ \mathcal{H}^\mathsf{L} \big\vert\gamma_i\text{ is open}  \big]
    \geq\phi^1_{\mathsf{Strip}^c}\big[ \mathcal{H}^\mathsf{L} \big].
\end{align*}
 For $\mathcal{H}^\mathsf{L}$  to fail, there must exist at least one index $k \geq 0$ such that $(-k, x_i)$ is connected to the vertical axis $\lbrace T_1 \rbrace\times \Z$ by a (primal) open path lying in the half-plane $(-\infty, T_1] \times \Z$. 
Thus
\begin{equation}\label{equ: equ preuve appendice finale somme}
    \phi^{1}_{\mathsf{Strip}^c}\big[(\mathcal{H}^\mathsf{L})^c\big] 
    \leq \sum_{k\geq 0}\phi^1_{\mathsf{Strip}^c} \big[ (-k, x_i) \leftrightarrow \lbrace T_1\rbrace\times \Z  \big].
\end{equation}
It is well-known that the exponential decay of the primal cluster applies also within wired boundary conditions~\cite{DuminilCopin2017SharpPT}, 
and therefore the terms in the sum above are bounded above by $e^{-c k }$ for some $c> 0$ and all $k$. 
Summing over $k\geq 0$ we find 
\begin{equation}
   \phi^{1}_{\mathsf{Strip}^c}\big[(\mathcal{H}^\mathsf{L})^c\big]  \leq C\e^{-cT_1}.
\end{equation}
By the cone confinement property $T_1 \geq \frac{n^\eps}{2\delta}$. 
This proves~\eqref{eq:HL} after altering the constants. 
\end{proof}

We are now ready to conclude. 
The claim along with~\eqref{eq:HHH} imply that 
\begin{align*}
    \phi^0_{\mathsf{Ext}^c} \Big[ \partial \gamma_i \text{ closed} \Big\vert  
    \gamma_i \text{ open}, \bigcap_{k=1}^{i-1} \lbrace \Gamma_k= \gamma_k, \mathsf{Diam}_k \rbrace\Big]
    &\leq \big(1 - \e^{-cn^\eps}\big)^{-1} \phi_{\mathsf{Ext}^c}^0\big[ \partial \gamma_i \text{ closed}\big| \gamma_i \text{ open} \big].
\end{align*}
The above, together with~\eqref{equ: majoration premier terme proba conditionnelle chemin ouvert} may be inserted into~\eqref{equ: decomposition chemin ouvert et ouvert sachant fermé} to obtain~\eqref{equ: equation a prouver majoration chemin par mesure produit}. As already mentioned, this concludes the proof of Lemma~\ref{lem: first step Ioan numerator}

\end{proof}

We turn to the second step of the proof of Lemma~\ref{lem:Ioan_numerator}. 

\begin{proof}[Proof of Lemma~\ref{lem:Ioan_numerator}]
Recall the definition of the family of paths $\Gamma$, defined when  $\nonint$, $\con$ and $\mathsf{Diam}$ occur. 
Define the set $\mathsf{ClosePath}$ of realisations $\gamma$ of $\Gamma$ for which there exists $2 \leq i\leq r$ and $T_1 \leq t \leq T_2$ such that 
\begin{align*}
    {\rm dist} \big[\gamma_{i-1} \cap (\{t\} \times \mathbb Z), \, \gamma_i \cap (\{t\} \times \mathbb Z)\big]   < (\log n)^3.
\end{align*}

We first observe that 
\begin{align}\label{equ: equation a remonter preuve lemme ioan numerator}
    \phi^0_{\mathsf{Ext}^c}&\left[ \mathsf{GlobRep}^c, \con, \nonint, \mathsf{Diam} \right]  \nonumber\\
    &\leq 
     \phi^0_{\mathsf{Ext}^c}\left[ \mathsf{GlobRep}^c,\BigDiam^c, \con, \nonint, \mathsf{Diam} \right]  + 
        \phi^0_{\Ext^c} \big[\BigDiam, \con,\nonint\big] \nonumber\\
    &\leq 
    \sum_{\gamma \in \mathsf{ClosePath}}\!\!\! \phi^0_{\Ext}\left[\Gamma= \gamma, \con, \nonint, \Diam \right] +
    \e^{-c(\log n)^{2} -\tau r (T_2 - T_1)}.
    \end{align}
Indeed, the second inequality is true term by term. For the first term, due to $\BigDiam^c$, the clusters are entirely within distance $(\log n)^2$ of the corresponding paths $\Gamma_i$. Thus, for $\mathsf{GlobRep}^c$ to occur, the paths $\Gamma_i$ need to come within distance $(\log n)^2 +2 (\log n)^2$ of each other, and in particular need to belong to $\mathsf{ClosePath}$. The bound on the second term is a direct consequence of Lemma~\ref{lemme diamond confinement mesure conditionnee} and \eqref{eq:tau}.
The second term obviously satisfies the upper bound in \eqref{eq:Ioan_numerator}, and we may focus on bounding the first term. 

By \eqref{equ: equation a remonter preuve lemme ioan numerator} and Lemma~\ref{lem: first step Ioan numerator}, 
\begin{align}
    &\phi^0_{\mathsf{Ext}^c}\big[  \mathsf{GlobRep}^c,\BigDiam^c, \con, \nonint, \mathsf{Diam} \big]   \nonumber\\
    & \leq  (1+o(1))
    \sum_{\gamma\in\mathsf{ClosePath}}(\phi^0_{\Ext^c})^{\otimes r}\Big[\bigcap_{i=1}^r\lbrace \gamma_i \text{ open}, \partial \gamma_i \text{ closed}\rbrace\Big]\nonumber\\
    &\leq
     (1+o(1))\, (\phi^0_{\Ext^c})^{\otimes r}\big[ \exists \gamma \in \mathsf{ClosePath} \text{ s.t. } \forall i\in\lbrace 1, \dots, r\rbrace, \gamma_i \text{ open}, \partial \gamma_i \text{ closed} \big].\label{eq:blabla}
\end{align}

The last upper bound is obtained by observing that when the last event is satisfied, at most one family of paths of $\mathsf{ClosePath}$ can achieve it (due to the event $\mathsf{Diam}$).
We will bound the last term using the Ornstein--Zernike coupling $\Phi^{0, \otimes r}_{\Ext^c, (T_1,\mathscr{X})\rightarrow(T_2, \mathscr{Y})}[~\cdot~]$ and the random skeleton system $\mathcal{S}$ given by this coupling. 
This argument will only be sketched as it already appeared in the proof of Proposition~\ref{Proposition convergence mesure produit}.

Under the event in the last line of \eqref{eq:blabla}, the paths $\gamma_1,\dots, \gamma_r$  contain all the renewal points of the clusters $\mathcal{C}_i$ in $\mathsf{Strip}$, and therefore the synchronised skeleton $\check{\mathcal{S}}$ is guaranteed to be non-intersecting. In addition, due to the diamond confinement property and since $\gamma \in \mathsf{ClosePath}$, 
$\inf_{t\in [T_1,T_2]} \Gap(\check{\mathcal{S}}(t)) \leq 3\log^2n + (\log n)^2$ with probability going to $1$. 
Finally we conclude that

\begin{align*}
    &(\phi^0_{\Ext^c})^{\otimes r}\big[  \exists \gamma \in \mathsf{ClosePath} \text{ such that } \forall i\in\lbrace 1, \dots, r\rbrace, \gamma_i \text{ is open}, \partial \gamma_i \text{ is closed} \big] \\ 
    &\leq (1+ o(1))\, \e^{-\tau r (T_2-T_1)}\Phi_{(T_1,\mathscr{X}) \rightarrow (T_2, \mathscr{Y})}^{\otimes r}\left[ \check{\mathcal{S}}\in \mathcal{W}_{[T_1,T_2]}, \, \Gap(\check{\mathcal{S}}(t)) \leq 4\log^2n\, \big\vert\, \EXT=\Ext \right]. 
\end{align*}

We now use once the local limit Theorem~\ref{theoreme local limit srw} and then Lemma~\ref{Lemme repulsion globale sdry} (the assumptions of the lemma are satisfied due to the edge-regularity of $\Ext$) to conclude that
\begin{align*}
    \Phi_{(T_1,\mathscr{X}) \rightarrow (T_2, \mathscr{Y})}^{\otimes r}&\left[ \check{\mathcal{S}}\in \mathcal{W}_{[T_1,T_2]}, \Gap(\check{\mathcal{S}}(t)) \leq 4\log^2n\big\vert\EXT=\Ext \right]\\ 
    &\leq CV(\mathscr{X})V(\mathscr{Y})(T_2-T_1)^{-\frac{r^2}{2}}(T_2-T_1)^{-\beta}.
\end{align*}
Putting everything together and coming back to~\eqref{equ: equation a remonter preuve lemme ioan numerator}, we obtain:
\begin{equation}
    \phi^0_{\Ext^c}\left[\mathsf{GlobRep}^c, \con, \nonint, \Diam\right] \leq C \, V(\mathscr{X})V(\mathscr{Y})\, (T_2-T_1)^{-\frac{r^2}{2}-\beta}\e^{-\tau r(T_2-T_1)},
\end{equation}
which concludes the proof.
\end{proof}

\subsection{The mixing argument and the proof of Theorems~\ref{theoreme estimation} and~\ref{Theoreme main} }\label{Sub finale mixing et preuves}

We start with the proof of Theorem~\ref{theoreme estimation}. 

\begin{proof}[Proof of Theorem~\ref{theoreme estimation}]

Recall that one bound, namely 
\begin{equation}\label{eq:TM1}
    \phi\left[\nonint, \con  \right] \geq c\,V(x)V(y)\,n^{-\frac{r^2}{2}}\e^{-\tau rn}
\end{equation}
has already been proved in Remark~\ref{remarque minoration ni con arbitrary BC}. Our goal here is to prove a matching upper bound. 

Running the argument used in the proof of Lemma~\ref{lem:Ioan_numerator}, but summing over all realisations of~$\Gamma$ rather than only those in $\mathsf{ClosePath}$ yields 
\begin{equation}\label{eq:Ioan_numerator2}
    \phi_{\Ext^c}^0 \left[ \nonint, \con \right] 
    \leq C V(\mathscr{X})V(\mathscr{Y}) (T_2 - T_1)^{-\frac{r^2}2} \e^{-\tau r (T_2 - T_1)}.
\end{equation}
The above does not match the desired bound since $T_1$ is larger than $0$ and $T_2$ smaller than $n$ by a polynomial quantity. Moreover, $V(\mathscr{X})$ and $V(\mathscr{Y})$ are also of polynomial order.
To obtain the upper bound matching \eqref{eq:TM1}, we will run the same argument with $T_1$ and $T_2$ replaced with random times of finite order.

Let $\tilde{T}_1$ (resp. $\tilde{T}_2$) be the first (resp. the last) synchronization point of the maximal skeletons of the clusters after 0 (resp. before $n$). We also call $\tilde{\mathscr{X}}$ (resp. $\tilde{\mathscr{Y}}$) the unique vector such that $(\tilde{T}_1, \tilde{\mathscr{X}}_i) \in \mathcal{C}_i$ (resp. $(\tilde{T}_2, \tilde{\mathscr{Y}}_i)\in\mathcal{C}_i$). We already argued in the proof of Lemma~\ref{lemme diamond confinement mesure conditionnee} that $\tilde{T}_1$ and $\tilde{T}_2$ have exponential tails: for any $t \geq 0$ large enough,
\begin{equation}
    \phi\left[ \max \lbrace \tilde{T}_1, \tilde{T}_2 \rbrace >t \vert \con,\nonint\right] \leq \e^{-ct}.
\end{equation}
We also already argued that 
\begin{equation}
    \phi\left[\max\lbrace\Vert\tilde{\mathscr{X}}-x\Vert, \Vert\tilde{\mathscr{Y}}-y\Vert\rbrace > t \vert \con, \nonint \right] \leq \e^{-ct},
\end{equation}
for a possibly different value of $c>0$. For the rest of this proof fix $t$ so that $e^{-ct} \leq 1/2$. Then, we upper bound:
\begin{align*}
    \phi\big[ \con, \nonint \big] &= \phi\big[\con, \nonint, \max\lbrace\Vert\tilde{\mathscr{X}}-x\Vert, \Vert\tilde{\mathscr{Y}}-y\Vert\rbrace > t\big] \\&\hspace{70pt}+ \phi\big[\con, \nonint, \max\lbrace\Vert\tilde{\mathscr{X}}-x\Vert, \Vert\tilde{\mathscr{Y}}-y\Vert\rbrace < t\big] \\ 
    &\leq \phi\big[\con, \nonint \big\vert \max\lbrace\Vert\tilde{\mathscr{X}}-x\Vert, \Vert\tilde{\mathscr{Y}}-y\Vert\rbrace < t\big] + \tfrac12\phi\big[ \con, \nonint\big].
\end{align*}
Hence, we obtain that 
\begin{equation}\label{equ: majoration con nonint via conditionnement preuve finale}
    \phi\left[\con,\nonint\right] 
    \leq 2\, \phi\big[\con, \nonint \big\vert \max\{\|\tilde{\mathscr{X}}-x\|, \|\tilde{\mathscr{Y}}-y\|\} < t\big].
\end{equation}
We focus on upper bounding $\phi[\con, \nonint \vert \max\{\|\tilde{\mathscr{X}}-x\|, \|\tilde{\mathscr{Y}}-y\|\} < t]$, and will do so using the method of Lemma~\ref{lem: first step Ioan numerator}. 
As in the proof of Lemma~\ref{lem: first step Ioan numerator}, condition on the shape of the clusters outside of $\mathsf{Strip} = [\tilde T_1,\tilde T_2] \times \mathbb Z$ to write:
\begin{multline*}
    \phi\left[\con, \nonint \vert \max\lbrace\Vert\tilde{\mathscr{X}}-x\Vert, \Vert\tilde{\mathscr{Y}}-y\Vert\rbrace < t\right] = \sum_{\Ext}\phi^0_{\Ext^c}\left[\con, \nonint, \mathsf{Dian} \right] \\ \times \phi\left[\EXT=\Ext \vert \max\lbrace\Vert\tilde{\mathscr{X}}-x \Vert, \Vert\tilde{\mathscr{Y}}-y\Vert\rbrace < t\right].
\end{multline*}
As previously, the conditioning on $\EXT$ contains the fact that $\tilde{T}_1$ and $\tilde{T}_2$ are renewals. 

Fix some $\Ext$ appearing in the sum above. Recall the definition of the top-most path $\Gamma_i$ of the cluster $\mathcal C_i$ and its upper boundary $\partial \Gamma_i$. Then 

\begin{multline}\label{eq:HH33}
    \phi^0_{\Ext^c}\big[\con, \nonint,\Diam\big]
    \leq\sum_{\gamma} \prod_{i=1}^r \phi^0_{\mathsf{Ext}^c} \Big[ \gamma_i \text{ open}\Big\vert  \bigcap_{k=1}^{i-1} \lbrace \Gamma_k= \gamma_k, \mathsf{Diam}_k \rbrace\Big] \\
  \times   \phi^0_{\mathsf{Ext}^c} \Big[ \partial \gamma_i \text{ closed} \Big\vert  \gamma_i \text{ open}, \bigcap_{k=1}^{i-1} \lbrace \Gamma_k= \gamma_k, \mathsf{Diam}_k \rbrace\Big],
\end{multline}
where the sum is over all possible realisations $\gamma$ of $\Gamma$ that induce disjoint connections between $\tilde{\mathscr{X}}$ and $\tilde{\mathscr{Y}}$.

The two terms in the right-hand side of the above may be bounded as in Lemma~\ref{lem: first step Ioan numerator} by 
\begin{align}
    \phi^0_{\mathsf{Ext}^c} \big[ \gamma_i \text{ open}\big\vert  \bigcap_{k=1}^{i-1} \lbrace \Gamma_k= \gamma_k, \mathsf{Diam}_k \rbrace\big]
    &\leq   \phi^0_{\mathsf{Ext}^c} \big[ \gamma_i \text{ open}\big] \text{ and } \nonumber  \\
     \phi^0_{\mathsf{Ext}^c} \big[ \partial \gamma_i \text{ closed} \big\vert  \gamma_i \text{ open}, \bigcap_{k=1}^{i-1} \lbrace \Gamma_k= \gamma_k, \mathsf{Diam}_k \rbrace\big]& \leq  C  \phi^0_{\mathsf{Ext}^c} \big[ \partial \gamma_i \text{ closed} \big\vert  \gamma_i \text{ open}\big]\label{eq:HHExt}
\end{align}
This is the only place where the proof differs a little from that of Lemma~\ref{lem: first step Ioan numerator}. 
Indeed, in the second bound above, we will use that 
\begin{align}\label{eq:HH3}
    \phi^0_{\mathsf{Ext}^c} \big[ \mathcal H^{\mathsf L} \cap \mathcal H^{\mathsf R} \big\vert  \gamma_i \text{ open}\big] \geq c.
\end{align}
for some universal constant $c$, where $ \mathcal H^{\mathsf L}$ and $\mathcal H^{\mathsf R}$ are defined as in Claim~\ref{claim 2}. 
The proof of \eqref{eq:HH3} is easier than that of Claim~\ref{claim 2}: it relies simply on the fact that, in the subcritical regime, the dual percolates even in a half-plane with wired boundary conditions. 

Now, injecting \eqref{eq:HHExt} back into \eqref{eq:HH33}, we find that 
\begin{align}\label{equ: equ majoration of con ni by a path event}
    \phi^0_{\Ext^c}\big[\con, \nonint, \Diam \big]
    \leq C\, (\phi^0_{\mathsf{Ext}^c})^{\otimes r} \big[\exists \gamma \text{ s.t. }  \forall i\in\lbrace 1, \dots, r\rbrace, \gamma_i \text{  open}, \partial \gamma_i \text{ closed}\big]
\end{align}

We conclude using the Ornstein--Zernike coupling for the product measure. Indeed, the event on the right-hand side of~\eqref{equ: equ majoration of con ni by a path event} implies that, in the product measure, the connection event occurs and the synchronized skeletons are non-intersecting. Thus,
\begin{equation}
    \phi^0_{\Ext^c}\big[\con, \nonint, \Diam\big] \leq C\,\e^{-\tau r (\tilde T_2-\tilde T_1)}\phi^{0, \otimes r }_{ (\tilde T_1,\tilde{\mathscr{X}})\rightarrow(\tilde  T_2,\tilde{\mathscr{Y}})}\big[  \check{\mathcal{S}}\in \mathcal{W}_{[\tilde  T_1, \tilde  T_2]} \big\vert \Ext=\Ext \big].
\end{equation}
We make use of the Local limit Theorem~\ref{theoreme local limit srw} to upper bound the right-hand side probability by $C\, V(\tilde{\mathscr{X}})V(\tilde{\mathscr{Y}})(\tilde T_2-\tilde T_1)^{-\frac{r^2}{2}}$. Using the assumption on $\Ext$, we then very roughly upper bound 
\begin{equation}
    \max\lbrace V(\tilde{\mathscr{X}}), V(\tilde{\mathscr{Y}}) \rbrace \leq  (2\max\lbrace \Vert\tilde{\mathscr{X}}-x\Vert, \Vert\tilde{\mathscr{Y}}-y\Vert \rbrace)^{\frac{r(r-1)}{2}}V(x)V(y) \leq V(x)V(y)(2t)^{\frac{r(r-1)}{2}}.
\end{equation}
Gathering everything together, we conclude that for all $\Ext$ satisfying $\max\{\|\tilde{\mathscr{X}}\|, \|\tilde{\mathscr{Y}}\|\} \leq t$ and $\max \lbrace \tilde{T}_1, n-\tilde{T}_2 \rbrace \leq t$
\begin{equation}
    \phi^0_{\Ext^c}\big[ \con,\nonint \big] \leq C\,(2t)^{\frac{r(r-1)}{2}} V(x)V(y) \e^{-\tau r(n-2t)}(n-2t)^\frac{r^2}{2}
  \leq C' V(x)V(y) \e^{-\tau r n}n^\frac{r^2}{2}. 
\end{equation}
Summing over all such $\Ext$ and using~\eqref{equ: majoration con nonint via conditionnement preuve finale} we find that
\begin{equation}
    \phi\big[\con, \nonint \big] \leq C'\, V(x)V(y)\e^{-\tau rn}n^{-\frac{r^2}{2}},
\end{equation}
where the value of $C'$ has been increased between every equation, but does not depend on $n$. This concludes the proof.

\end{proof}

We then turn to the proof of Theorem~\ref{Theoreme main}. It follows from the repulsion estimate of Proposition~\ref{proposition global entropic repulsion} and the convergence of the product system stated in Proposition~\ref{Proposition convergence mesure produit}. The observation is that when $\mathsf{GlobRep}$ occurs, then it is a consequence of the mixing property of the random-cluster measure~\eqref{weak ratio mixing} that the distribution of the system of clusters is very close to the one of an \emph{independent} system of clusters.

\begin{proof}[Proof of Theorem~\ref{Theoreme main}]

We follow the same pattern as in the proof of Proposition~\ref{Proposition convergence mesure produit}, and use the strategy given by Lemma~\ref{Lemme technique processus sto}. Fix some $\delta >0$ and arbitrary signs for the $r$ envelopes which we denote by $\pm$. As previously, define the scaled process $\Gamma^{\pm}_n(t) := \tfrac{1}{\sqrt{n}}\Gamma^\pm(nt)$, and consider a function $f : \mathcal{C}([\delta, 1-\delta], \R^r) \rightarrow \R$, continuous and bounded. As in the proof of Proposition~\ref{Proposition convergence mesure produit} we shall omit to write the restriction to the interval $[\delta, 1-\delta]$ when writing $f^\delta(\Gamma^\pm_n)$. We start by arguing that due to the boundedness of $f^\delta$ and to Lemma~\ref{lem: edge regularity true measure},
\begin{equation}
\phi\big[ f^\delta(\Gamma^\pm_n)\big\vert \nonint,\con \big] = (1+o(1))\phi\big[ f^\delta(\Gamma^\pm_n)\big\vert \nonint,\con, \mathsf{EdgeReg} \big].
\end{equation}

Next, we condition on $T_1$, $T_2$ and the edge-regular shape of the clusters outside of $\mathsf{Strip}:=\mathsf{Strip}_{[T_1,T_2]}$.
Due to the edge-regularity condition, we can chose $n$ large enough so that $n\delta>T_1$ and $n(1-\delta) < T_2$. 
We find
\begin{align*}
    \phi\big[f^\delta(\Gamma^\pm_n) \big\vert \nonint,\con, \mathsf{EdgeReg}\big]  
    =   \sum_{\Ext} &\phi\big[f^\delta\big(\Gamma^\pm_n\big) \big\vert \nonint, \con, \EXT = \mathsf{Ext} \big] \\ 
    &\times \phi\big[\EXT = \mathsf{Ext}\big\vert\nonint,\con, \mathsf{EdgeReg}\big].
\end{align*}
Fix some $\Ext$ which is edge-regular. 
We make use of Proposition~\ref{proposition global entropic repulsion} to argue that:
\begin{multline}
    \phi\left[f^\delta(\Gamma^\pm_n) \vert \nonint, \con, \EXT = \mathsf{Ext} \right] = \\ \left(1+o(1)\right)\phi\left[f^\delta(\Gamma^\pm_n) \vert \nonint, \con, \EXT = \mathsf{Ext}, \mathsf{GlobRep} \right].
\end{multline}

Now, we claim that the mixing property~\eqref{weak ratio mixing} implies that
\begin{equation}\label{equ: equation finale mixing}
    \bigg| \frac{\phi\big[f^\delta(\Gamma^\pm_n) \vert \nonint, \con, \EXT = \mathsf{Ext}, \mathsf{GlobRep} \big]}{\phi^{\otimes r}\big[f^\delta(\Gamma^\pm_n) \vert \nonint, \con, \EXT = \mathsf{Ext}, \mathsf{GlobRep} \big]} - 1\bigg| < \e^{-2(\log n)^2}.
\end{equation}
Indeed, decompose the term $\phi\big[f^\delta(\Gamma^\pm_n) \vert \nonint, \con, \EXT = \mathsf{Ext}, \mathsf{GlobRep} \big]$ as follows:
\begin{equation}
    \phi\big[f^\delta(\Gamma^\pm_n) \vert \nonint, \con, \EXT = \mathsf{Ext}, \mathsf{GlobRep} \big] = \!\!\!\sum_{C_1, \dots, C_r}\!\!\! f^\delta(\Gamma^\pm_n)\frac{\phi^0_{\Ext^c}\big[\mathcal{C}_1=C_1, \dots, \mathcal{C}_r=C_r\big]}{\phi^0_{\Ext^c}\big[\nonint, \con, \mathsf{Globrep} \big]},
\end{equation}
where the sum runs over the possible realisations $C_1, \dots, C_r$ of the clusters of $\mathscr{X}$ under the measure $\phi^0_{\Ext^c}[\cdot\vert \nonint, \con, \mathsf{GlobRep}]$. The point is that those sets are almost surely finite and have a mutual distance larger than $\delta (\log n)^3$ by the diamond confinement property. We can then apply~\eqref{weak ratio mixing} to both the numerator and the denominator of the fraction to obtain~\eqref{equ: equation finale mixing}.

The last thing to notice is that the entropic repulsion estimate~\eqref{entropic repulsion equation} also holds for the product measure: 
\begin{equation}
    \phi^{\otimes r}\left[ \mathsf{GlobRep} \vert \con_{\mathscr{X}, \mathscr{Y}}, \nonint, \EXT = \mathsf{Ext} \right] \geq 1-cn^{-\beta},
\end{equation}
because of the usual Ornstein--Zernike coupling~\eqref{equation OZ boundary conditions} and the entropic repulsion for random walks given by Lemma~\ref{lemme entropic repulsion sdrw}. 
In conclusion, we proved that:
\begin{equation}
    \phi\left[f^\delta(\Gamma^\pm_n) \vert \nonint, \con, \EXT = \mathsf{Ext} \right] = (1+o(1))\phi^{\otimes r}\left[f^\delta(\Gamma^\pm_n) \vert \nonint, \con, \EXT = \mathsf{Ext} \right].
\end{equation}
Finally, due to the assumed edge-regularity of $\Ext$ and to Proposition~\ref{Proposition convergence mesure produit}, we know that the RHS converges towards $\E\left[f^\delta(\sigma\bw^{(r)})\right].$
Hence, 
\begin{equation}
    \phi\left[ f^\delta(\Gamma^\pm_n) \vert \nonint, \con \right] \goes{}{n \rightarrow \infty}{\E\left[f^\delta(\sigma\bw^{(r)})\right]},
\end{equation}
and so is established point $(i)$ of Lemma~\ref{Lemme technique processus sto}. 
As previously the equicontinuity at 0 and 1 is an easy consequence of basic large deviations estimates. This observation achieves the proof. 
\end{proof}

\section{Local statistics of directed non-intersecting random bridges}\label{section marches}

As seen before, by the Ornstein--Zernike theory and the entropic repulsion, a system of clusters subject to the non-intersection conditioning resembles a system of non-intersecting directed random walks.  

Non-intersecting random walks, and more largely random walks in cones have a very rich combinatorial and probabilistic structure. They have been studied widely throughout the last 50 years. The seminal work is the paper of Karlin and McGregor \cite{karlin1959} which proves a determinantal formula for the probability of $r$ random walks to intersect. Their approach only applies to a very specific class of walks, and is combinatorial by nature; it lead to remarkable developments around integrable systems of walks (see \cite{GRABINER1999177,johansonn2004}).

A more probabilistic treatment has been started in~\cite{konig2005},~\cite{conditionallimittheoremsfororderedrandomwalks,randomwalksincones, Invarianceprinciplesforrandomwalksincones}. Indeed, in~\cite{konig2005} a definition of the random walk conditioned to stay in a cone was given in terms of a Doob $h$-transform by a harmonic function vanishing on the boundary of the cone, allowing the authors to obtain Local Limit Theorems and invariance principles for a much broader class of random walks.  We briefly summarize the definitions and construction of the concerned objects.

The goal of this section is then to study the properties of such systems of walks, especially their behaviour under the diffusive scaling. Let us introduce the relevant object to study.  

\begin{definition}[Directed system of random walks] \label{Definition directed random walk}
Let $r\geq 1 $ be an integer, and for $1 \leq i \leq r$, let $(\theta_n^i, X_n^i)_{n \geq 1}$ be an independent and identically distributed family of independent and identically distributed random variables on $\N^* \times \Z$. We assume that it satisfies the following properties:
\begin{itemize}
    \item Both $\theta_1^1$ and $X_1^1$ have an exponential moment.   
    \item Conditionally on $\theta_n^i$, $X^i_n$ is centered.
\end{itemize}
Call:
\begin{equation}
    \mathbf{T}^i_n = \sum_{k=1}^n \theta^i_k\qquad  \text{  and  } \qquad \mathbf{Z}^i_n = \sum_{k=1}^n X^i_k
\end{equation}
Then the system 
\begin{equation}
  \left(\mathbf{S}_n\right)_{n \geq 0} := \left((\mathbf{T}^1_n, \mathbf{Z}^1_n), \dots, (\mathbf{T}^r_n, \mathbf{Z}^r_n)\right)_{n \geq 0}
\end{equation}
is called a \textit{system of directed random walks}. For any $(k_i,x_i)_{1\leq i \leq r} \in \left(\N\times \Z\right)^r$, we write $\mathbf{P}_{(k,x)}$ for the law of the $r$-directed random walk with $\mathbf{S}^i_0 = (k_i,x_i)$ -- this is defined as above, with the addition of an initial offset. 
When all the $k_i$ are equal, which will often be the case, we make a slight abuse of notation by writing $\mathbf{P}_{(k,x)}$ with $k\in\Z, x\in\Z^r$.
\end{definition}

As observed in the precedent sections, a subcritical percolation cluster can be roughly described as the trajectory of a directed random walk \textit{decorated} with $\delta$-confined clusters of edges. For that reason, it is convenient to study directed system of non-intersecting random bridges carrying $\delta$-diamonds around their steps. We then make the following assumption:
\begin{assumption}\label{assumption 1} There exists a $\delta > 0$ such that almost surely, 
\begin{equation}\label{eq: assumption 1}
    (\theta^1_1, X^1_1) \in \mathcal{Y}^{+,\delta}_0. 
\end{equation}
\end{assumption}

Recall the definition of the diamonds from Section~\ref{section review oz}. If $(\mathbf{S}_n)_{n \geq 0}$ is a system of directed random walks, we introduce 

\begin{equation}\mathcal{D}^\delta_{i,k} = \mathcal{D}^\delta_{(\mathbf{T}^i_k, \mathbf{Z}^i_k),(\mathbf{T}_{k+1}^i, \mathbf{Z}^i_{k+1})} \qquad\text{and}\qquad \mathcal{D}(\mathbf{S}^i) := \bigcup_{k\geq 0}\mathcal{D}^\delta_{i,k}.
\end{equation}

\begin{figure}
    \centering
    \includegraphics[width = .65\textwidth, page = 4]{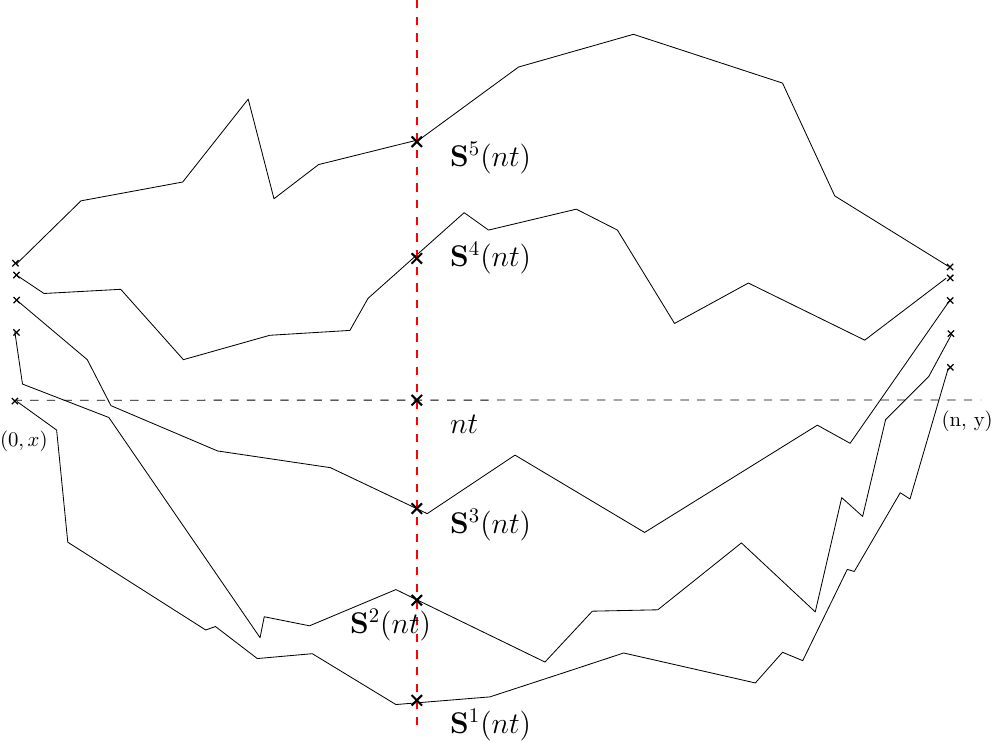}
    \caption{A depiction of a system of $r=5$ non-synchronized random walks under the event $\left\lbrace \mathbf{S} \in \mathcal{W}_n, \hit_{(n,y)}\right\rbrace$.}
    \label{figure non-synchronized random walks}
\end{figure}

We also introduce the diamond-decorated walks analogs of the events $\con$ and $\nonint$ (see Figure~\ref{figure non-synchronized random walks}).

For $y \in \R^r, n\geq 0$, the hitting event is defined by:
\begin{equation}
    \hit_{(n,y)} = \left\lbrace \exists k_1, \dots, k_r \geq 0, \forall i \in \{1,\dots,r\}, (\mathbf{T}^i_{k_i},\mathbf{S}^1_{k_1}) = (n, y_i) \right\rbrace.
\end{equation}
The non-intersection of diamond event is defined by 
\begin{equation}
\nidiam(\mathbf{S}) = \bigcap_{1\leq i \neq j \leq r} \{\mathcal{D}(\mathbf{S}^i) \cap \mathcal{D}(\mathbf{S}^j) = \emptyset\}.
\end{equation}

The goal of this section is the proof of the following result:

\begin{theorem}[Invariance principle for directed random walks]\label{theoreme invariance principle for drw}
Let $S$ be a system of $r$ directed random walks sampled according to $\mathbf{P}_{(0,x)}\left[\cdot\vert \mathbf{S} \in \hit_{(n,y)}, \nidiam(\mathbf{S}) \right]$. Then, there exists $\sigma> 0$ such that
\begin{equation}
    \left(\frac{1}{\sqrt{n}}\mathbf{S}\left(nt\right) \right)_{0\leq t \leq 1} \goes{(d)}{n \rightarrow \infty}{\left(\sigma\bw^{(r)}_{t}\right)_{0\leq t \leq 1}},
\end{equation}
where the convergence holds in the space $\mathcal{C}([0,1], \R^r)$ equipped with the topology of the uniform convergence. Moreover, $x$ and $y$ can depend on $n$ in the statement, as long as they both have norms that are $o(\sqrt{n}).$
\end{theorem}

This theorem, as well as Theorem~\ref{theoreme local limit srw}, has already been proved in the setting of \textit{regular} random walks (that is when $\theta_1^1= 1$ almost surely), and replacing the conditioning over the non-intersection of \textit{diamonds} by a conditioning of non-intersection of their \textit{spatial trajectories}. Our goal here is simply to extend this to the setting of directed random walks decorated with diamonds, and state some properties tailored to our needs. The key object to derive this statement is the \textit{embedded synchronized system of directed random walks}. We define it in the next section and derive the key input for the study, namely the Local Limit Theorem~\ref{theoreme local limit srw}.

\subsection{Synchronized directed random walks}\label{subsection non-intersecting systems of synchronized directed random walks}

As observed in Sections~\ref{section independent system} and~\ref{section RCM} our arguments are often soft enough to boil down to the study of a \emph{synchronized} system of walks, where the time reference is still random but common to every walk.

\begin{definition}[Synchronized directed random walk]\label{def synchronized system of random walks}
Let $r\geq 1$ be an integer. Let $(\theta_k, X^1_k, \dots, X^r_k)_{k\geq 0}$ be a sequence of independent and identically distributed random variables taking values in $\N^* \times \Z^r$. Moreover we assume that
\begin{itemize}
    \item Both $\theta_1$ and $X_1^1$ have an exponential moment 
    \item Conditionally on  $\theta_1$, $X_1^1, \dots, X_1^r$ are centered, independent and identically distributed. 
\end{itemize}
We call:
\begin{equation}
    T_n = \sum_{k=1}^n \theta_k \quad \text{  and  } \quad Z^i_n = \sum_{k=1}^n X^i_k
\end{equation}
Then the system 
\begin{equation}
  \left(S_n\right)_{n \geq 0} = \left(T_n, Z^1_n, \dots, Z^r_n\right)_{n \geq 0}
\end{equation}
is called a \textit{synchronized system of directed random walks}. In what follows, we see it as a random element of $\N \times \Z^{r}$, which we will refer to as $\left(S_n\right)_{n\geq 0} = \left(T_n, Z_n\right)_{n \geq 0}.$ For any $(k,x) \in \N\times\Z^r$, we will denote by $\PP_{(k,x)}$ the law of the synchronized $r$-random walk started from the point $(k,x)$, \textit{i.e.} the law of $\left( (k,x)+S_n\right)_{n \geq 0}.$ 
\end{definition}

We also assume for convenience that Assumption~\ref{assumption 1} holds. 

Introduce the following hitting event, for any $(n,y) \in \N\times\Z^r$:
\begin{equation}
    \hit_{(n,y)} = \left\lbrace \exists k \geq 0, S_k = (n,y) \right\rbrace,
\end{equation}
and the stopping time
\begin{eqnarray}
    H_{(n,y)} = \min\left\lbrace k \geq 0, S_k = (n,y) \right\rbrace.
\end{eqnarray}
Moreover, $\rho$ will denote the stopping time corresponding to the first exit of the Weyl chamber:
\begin{equation}
    \rho= \min \left\lbrace n \geq 0,  S_n \notin W \right\rbrace
\end{equation}

The key results of this section are the following:

\begin{theorem}[Local limit Theorem for synchronized, non-intersecting directed random walks]\label{theoreme local limit srw} 
Let $\left( S_n \right)_{n \geq 0}$ be a synchronized system of random walks. There exists a function $V: W\cap\Z^r \rightarrow \R^*_+$ and a constant $C_1 > 0$ such that for any pair of sequences $(x(n))_{n \geq 0}, (y(n))_{n\geq 0}$ taking values in $W$ such that $\norme{x(n)}, \norme{y(n)} = o(\sqrt{n})$, when $n \goes{}{}{\infty}$, 
\begin{equation}
    \PP_{(0,x(n))} \left[ H_{(n,y(n))}<\rho, \hit_{(n,y(n))} \right] = C_1\frac{V(x(n))V(y(n))}{n^{r^2/2}}\left(1+o(1)\right).
\end{equation}
Furthermore, the function $V$ satisfies the following set of properties:
\begin{itemize}
    \item If $x,y \in W$ are such that $\left| y_{i+1} - y_i \right| > \left|x_{i+1}-x_i\right|$ for any $1 \leq i \leq r-1$, then
    \begin{equation}
        V(y) \geq V(x).
    \end{equation}
    \item When  $\Gap(x) \rightarrow \infty$, then $\frac{V(x)}{\Delta(x)} \rightarrow 1$.
    \item There exists a positive $c>0$ such that
    \begin{equation}
            V(x) \leq c \prod_{1\leq i < j \leq r}\left|1+x_j-x_i\right|
    \end{equation}
\end{itemize}
\end{theorem}

The second local limit result is the analog of Gnedenko's Local Limit Theorem. It corresponds to~\cite[Thm. 5]{randomwalksincones}.

\begin{theorem}\label{gnedenko theorem}
Let $\left( S_n \right)_{n \geq 0}$ be a synchronized system of random walks. Then, there exists a constant $\varkappa > 0$ such that for any fixed $x \in W$,
\begin{equation}
    \sup_{y \in W} \left| n^{\frac{r(r+1)}{4}}\PP_{(0,x)}\left[ S \in \mathcal{W}_n, \hit_{(n,y)} \right] - \varkappa V(x)\Delta\left(\frac{y}{\sqrt{n}}\right)\e^{-\frac{\norme{y}^2}{2n}} \right| \goes{}{n \rightarrow \infty}{0}.   
\end{equation}
\end{theorem}

The last theorem of this section is the invariance principle stating that a synchronized system of random walks conditioned on the events $\{ H_{(n, y(n))} < \rho \}$ and $\{\hit_{(n,y(n)))}\}$ converges towards the Brownian watermelon. 

\begin{theorem}[Invariance principle for synchronized, non-intersecting random walks]\label{theoreme d'invariance sdrw}
Let $\left(S_n\right)_{n \geq 0}$ be a system of $r$ synchronized random walks. We study the trajectory of $S$ on $[0, H_{(n,y)}] $ under the measure 
\begin{equation}
\PP_{(0,x)} \left[ ~\cdot~ \vert H_{(n,y)} < \rho, \hit_{(n,y)} \right] .
\end{equation}
Let $\mathfrak{T}$ be the linear interpolation between the points $(T_1, S_1), \dots, (n,y)$, and 
$S(t)$ be the almost surely unique intersection $\mathfrak{T}\cap \left(\lbrace t \rbrace\times \R^r\right)$. Then, there exists $\sigma>0$ such that:
\begin{equation}
    \left( \frac{1}{\sqrt{n}} S(nt) \right)_{0 \leq t \leq 1} \goes{}{n\rightarrow \infty}{\left(\sigma\bw^{(r)}_t\right)_{0 \leq t \leq 1}},
\end{equation}
The convergence occurs in the space $\mathcal{C}\left([0,1], \R^r \right)$ endowed with the topology of uniform convergence. Moreover, the convergence holds when $x, y$ depend on $n$, still as long as their norm is $o(\sqrt{n})$. 
\end{theorem}

Theorems~\ref{theoreme local limit srw},~\ref{gnedenko theorem} and~\ref{theoreme d'invariance sdrw} have already been derived in the works~\cite{konig2005},~\cite{conditionallimittheoremsfororderedrandomwalks,randomwalksincones} and most importantly in~\cite{Invarianceprinciplesforrandomwalksincones} in the case of \textit{regular} random walks, meaning that $\theta_1 = 1$ almost surely. Moreover it has been explained in great detail how to adapt the proofs of these articles to the case of \textit{directed} walks in \cite{ottvelenikwachtelioffe}. For that reason, a very brief sketch of proof of these three important results is deferred to the Appendix. 

We also import fast repulsion estimates that are going to be useful later on

\begin{lemma}[Edge repulsion for synchronized random walks]\label{lemme entropic repulsion sdrw}
There exists $\eps>0$ such that the following holds. Let $\left( S^i_n\right)_{n \geq 0,  1 \leq i \leq r}$ be a synchronized system of directed random walks. Let
\begin{equation}\label{definition stopping time}
    \eta_n = \min \left\lbrace k \geq 0, \min_{1 \leq i<j \leq r} \left| S^i_k - S^j_k \right| > n^\eps  \right\rbrace.
\end{equation}
Then there exists $c>0$ such that for any $x \in \R^r$, when $n$ is sufficiently large,
\begin{equation}\label{eq:repulsion_sdrw}
    \PP_{(0,x)}\left[ \eta_n > n^{1-\eps} \right] < \tfrac{1}{c}
    \exp(-cn^\eps).
\end{equation}
\end{lemma}
\begin{proof}
This fact has been proved in~\cite[Lemma 7]{conditionallimittheoremsfororderedrandomwalks} in the case of regular random walks, with a stronger statement: indeed, the $n^{\eps}$ in the definition of $\eta_n$ is replaced by $n^{\frac{1}{2}-\eps}$ in the latter paper. We briefly explain how to derive the result in our setting. First, condition on the time increments $(\theta_k)_{k \geq 0}$. The spatial increments become a sequence of independent (though non identically distributed) random variables. However one can check that the proof of~\cite[Lemma 7]{conditionallimittheoremsfororderedrandomwalks} can be \textit{mutatis mutandi} repeated in that setting.  
\end{proof}

\begin{remark}
Since this probability in \eqref{eq:repulsion_sdrw} is stretch-exponentially small, the bound also holds - up to a change in the constant $c$ - when conditioning the synchronized system of random walks on an event of polynomial probability. 
In particular, the next corollary follows from \eqref{eq:repulsion_sdrw} and Theorem~\ref{theoreme local limit srw} (which will be proved shortly without the use of the statement below).  
\begin{corollary}\label{cor: repulsion SDRbridges}
    There exists $\eps>0$ such that the following holds. Let $\left( S^i_n\right)_{1 \leq i \leq r, n \geq 0 }$ be a synchronized system of directed random walks. Then there exists $c>0$ such that for any $x,y \in \R^r$, when $n$ is sufficiently large,
    \begin{equation}
    \PP_{(0,x)}\left[ \eta_n > n^{1-\eps} \vert \hit_{(n,y)} \right] < \tfrac{1}{c}\exp(-cn^\eps).
\end{equation}
\end{corollary}    
\end{remark}

Using the input given by the Local Limit Theorem~\ref{theoreme local limit srw}, we are now able to derive the essential bulk repulsion for non-intersecting synchronized random walks in the next two lemmas.

\begin{lemma}\label{lemme repulsion bulk}
Let $S$ be a system of directed synchronized random walks and fix $\eps> 0$. Then, for any $\delta > 0$  sufficiently small, any points $x,y \in W$ satisfying $\norme{x}, \norme{y} = o(\sqrt{n})$, there exist $\beta > 0, C>0$ such that for any $n \geq 0$ sufficiently large, 
\begin{equation}\label{equation a démontrer bulk rep synchronized rw}
    \PP_{(0,x)} \big[ \exists t \in [n ^\eps, n-n^\eps], \Gap(S(t)) \leq n^\delta \big| H_{(n,y)}<\rho , \hit_{(n,y)} \big]  \leq Cn^{-\beta}.
    \end{equation}
\end{lemma}

\begin{proof} 
First, notice that one can actually examine only integer values of $t$ in~\eqref{equation a démontrer bulk rep synchronized rw} since the minimal distance between two synchronized piecewise linear functions is achieved at a slope change time, which by definition of $S$ is an integer. Introduce the following kernel:
\begin{equation}\label{Equation definition kernel}
    q_n(x,y) = \PP_{(0,x)}\left[ S \in \mathcal{W}_n, \hit_{(n,y)} \right].
\end{equation}

By the union bound it is sufficient to prove that
\begin{equation}
    \PP_{(0,x)} \big[ \exists k \in \lbrace n^\eps, \dots, n-n^{\eps} \rbrace , |S^{i}(k) - S^{i-1}(k) | \leq n^\delta \, \big|\,  S \in \mathcal{W}_n, \hit_{(n,y)} \big] \leq Cn^{-\beta}.
\end{equation}
for any $2 \leq i <r$.
Fix such an $i$ and introduce the following subset of $W$ 
\begin{equation}
    W_{n,\delta} = \left\lbrace u \in W, \left| u_{i}- u_{i-1} \right| < n^\delta \right\rbrace.
\end{equation}
We make use of Theorems~\ref{theoreme local limit srw} and~\ref{gnedenko theorem}. 
Indeed,  choose $n$ large enough so that for $n^\eps < k < n-n^\eps$, one has that for any $u\in W_{n,\delta}$:
\begin{equation}
    \begin{cases}
        q_n(x,y) &\geq (1-\eps)V(x)V(y)n^{-\frac{r^2}{2}} \\
        q_k(x,u) &\leq 2V(x)\Delta\big(\frac{u}{\sqrt{k}}\big)k^{-\frac{r(r+1)}{4}}\e^{-\frac{\norme{u}^2}{2k}} \\
        q_{n-k}(u,y) &\leq 2V(y)\Delta\big(\frac{u}{\sqrt{n-k}}\big)(n-k)^{-\frac{r(r+1)}{4}}\e^{-\frac{\norme{u}^2}{2(n-k)}}.
    \end{cases}
\end{equation}
Then, a union bound over $k$ yields: 
\begin{eqnarray*}
&& \PP_{(0,x)} \left[ \exists k \in \lbrace n^\eps, \dots, n-n^{\eps} \rbrace , \left|S^{i+1}(k) - S^i(k) \right| \leq n^\delta \,\Big\vert\,  S \in \mathcal{W}_n, \hit_{(n,y)} \right] \\
&\leq& \sum_{k=n^\eps}^{n-n^{1-\eps}} \sum_{u\in W_{n,\delta}} \frac{q_k(x,u)q_{n-k}(u,y)}{q_n(x,y)} \\
&\leq& \tfrac{4}{1-\eps}\sum_{k=n^\eps}^{n-n^{1-\eps}} \sum_{u\in W_{n,\delta}} n^\frac{r^2}{2} (k(n-k))^{-\frac{r(r+1)}{4}}\Delta\left(\tfrac{u}{\sqrt{k}}\right)\Delta\left(\tfrac{u}{\sqrt{n-k}}\right)\e^{-\frac{\norme{u}^2}{2}\left(\frac{1}{k}+\frac{1}{n-k}\right)}.
\end{eqnarray*}

We make two observations: the first is that this sum is actually symmetric around $\frac{n}{2}$, so that it is sufficient to bound it for $k$ going from $n^\eps$ to $\frac{n}{2}$. The second is that since $u\in W_{n,\delta}$, we have:
\begin{equation}
    \Delta\left(\tfrac{u}{\sqrt{k}}\right) \leq 2\norme{u}^{\frac{r(r-1)}{2}-1}n^\delta \, k^{-\frac{r(r-1)}{4}}.
\end{equation}
Then,
\begin{multline*}
 \PP_{(0,x)} \Big[ \exists k \in \lbrace n^\eps, \dots, n-n^{\eps} \rbrace , |S^{i+1}(k) - S^i(k)| \leq n^\delta \,\Big\vert\,  \hit_{(n,y)}, \tau > H_{(n,y)} \Big] \\
\leq \tfrac{32}{1-\eps}\sum_{k=n^\eps}^{\frac{n}{2}}\left(\tfrac{n}{k(n-k)}\right)^{\frac{r^2}{2}}n^{2\delta}\underbrace{\sum_{u\in W_{n,\delta}}\norme{u}^{r(r-1)-2}\e^{-\frac{\norme{u}^2}{2k}}}_{I}.
\end{multline*}
We then evaluate the order of the sum $I$. Indeed, let us write:
\begin{equation}
    I = \sum_{\ell\geq 0}\sum_{\substack{u\in W_{n,\delta} \\ \norme{u}=\ell}} \ell^{r(r-1)-2}\e^{-\frac{\ell^2}{2k}} \lesssim \sum_{\ell \geq 0} n^\delta \ell^{r-2}\ell^{r(r-1)-2}\e^{-\frac{\ell^2}{2k}},
\end{equation}
where we have used the fact that when $\ell \rightarrow \infty$, if $B_\ell(0)$ denotes the $\norme{\cdot}$ ball of $\R^r$ centered at 0 and of radius $\ell$, then
\begin{equation}\label{equation estimation volume twisted weyl chamber}
    \left| W_{n,\delta} \cap \partial B_\ell(0) \right| \lesssim n^\delta \ell^{r-2}.
\end{equation}
We then compare the latter sum with the integral $\int_{x=0}^\infty x^{r^2-4}e^{-\frac{x^2}{2k}}\dif x,$
which after the change of variables $t = \frac{x^2}{2k}$, can be explicitly evaluated:
\begin{equation}
    \int_{x=0}^\infty x^{r^2-4}e^{-\frac{x^2}{2k}}\dif x = \sqrt{2}^{r^2-5}\Gamma\left( \tfrac{r^2-3}{2}\right)k^{\frac{r^2-3}{2}}.
\end{equation}

Inserting this into our previous computation we find
\begin{align*}
 \PP_{(0,x)} \Big[ \exists k \in \lbrace n^\eps, \dots, n-n^{\eps} \rbrace , |S^{i+1}(k) - S^i(k) | \leq n^\delta \,\Big\vert\,  S \in \mathcal{W}_n, \hit_{(n,y)} \Big] &\\
\leq C\sum_{k=n^\eps}^{n/2}\left( \tfrac{n}{k(n-k)}\right)^\frac{r^2}{2}n^{3\delta}k^{\frac{r^2-3}{2}} 
\leq Cn^{3\delta}\sum_{k=n^\eps}^\infty k^{-\frac{3}{2}} 
&\leq Cn^{3\delta}n^{-\frac{\eps}{2}}.
\end{align*}
Hence, whenever $\delta<\frac{\eps}{6}$, this probability decays polynomially, as announced.
\end{proof}

In the proofs of Sections~\ref{section independent system} and~\ref{section RCM}, we used this lemma under a slightly different form that we state now. 

\begin{lemma}\label{Lemme repulsion globale sdry}
Let $S$ be a directed system of synchronized random walks and $\eps > 0$. Let $x(n), y(n)$ two sequences of elements of $W$ such that 
\begin{equation}
    \min\lbrace\Gap(x(n)), \Gap(y(n))\rbrace \geq  n^\eps
    \quad \text{ and } \quad 
    \norme{x(n)}, \norme{y(n)} = o(\sqrt{n}).
\end{equation}
Then, for any $\delta > 0$ sufficiently small, there exist $\beta > 0$ and $C> 0$ such that for $n\geq 0$ large enough, 
\begin{equation}
    \PP_{(0,x(n))} \left[ \inf_{0 \leq t \leq n }\Gap(S(t)) \leq n^\delta~\Big\vert~S \in \mathcal{W}_n, \hit_{(n,y(n))}\right] \leq Cn^{-\beta}.
\end{equation}
\end{lemma}

\begin{proof}
All the work has been done in Lemma~\ref{lemme repulsion bulk}. Indeed, we already know that 
\begin{equation}
    \PP_{(0,x(n))} \left[ \inf_{n^\eps\leq t\leq n-n^\eps}\Gap(S(t)) \leq n^\delta~\Big\vert~S \in \mathcal{W}_n, \hit_{(n,y(n))}\right] \leq Cn^{-\beta}.
\end{equation}
It remains to control the range of indexes $k \in \lbrace 1, \dots, n^\eps\rbrace \cup \lbrace n-n^\eps, \dots, n\rbrace$ (observe that we cannot make use of the local limit theorems in this range). However it is a basic large deviations estimate: let us write it for $k \in \lbrace 0,\dots, n^\eps \rbrace$.  We roughly bound 
\begin{multline}
     \PP_{(0,x(n))} \left[ \inf_{0\leq k\leq n^\eps} \Gap(S(k)) \leq n^\delta~\Big\vert~S \in \mathcal{W}_n, \hit_{(n,y(n))}\right] \\
    \leq \frac{\PP_{(0,x(n))} \left[\inf_{0\leq k\leq n^\eps} \Gap(S(k)) \leq n^\delta\right]}{\PP_{(0,x(n))} \left[S \in \mathcal{W}_n, \hit_{(n,y(n))}\right]}.
\end{multline}
Now observe that for the event of the numerator to occur, one of the walks has to travel at a distance at least $\frac{1}{2}(n^\eps- n^\delta)$ of its starting point in a time $n^\eps$, which by large deviations occurs with stretched exponentially small probability as soon as $\delta < \eps$. Additionally, by Theorem~\ref{theoreme local limit srw}, the denominator is of order at most polynomial. Thus 
    \begin{align*}
        \PP_{(0,x(n))} \left[ \inf_{0\leq k\leq n^\eps} \Gap(S(k)) \leq n^\delta~\Big\vert~S \in \mathcal{W}_n, \hit_{(n,y(n))}\right] \leq C \e^{-cn^\eps},
    \end{align*}
    for constants $c,C >0$. The same holds for $k \in \lbrace n-n^\eps, \dots, n\rbrace$, and the union bound provides the desired result.  
\end{proof}

\subsection{Synchronized systems of random walks with random decorations}\label{Subsection synchronized systems of random walks with random decorations}

In order to prove Theorem~\ref{theoreme invariance principle for drw}, we are going to compare a system of decorated non-intersecting random bridges with a system of decorated non-intersecting synchronized random bridges. This motivates us to study the properties of such a system. Recall the definition of $\mathcal{D}(\mathbf{S}^i)$ from the precedent section. When $S$ is a synchronized system of directed walks, we simply set     
\begin{equation}
    \mathcal{D}^{\delta}_{i,k} = \mathcal{D}^\delta_{\left(T_k, S^i_k\right), \left(T_{k+1}, S^i_{k+1}\right)}.
\end{equation}
and 
\begin{equation}
    \mathcal{D}(S^i) = \bigcup_{k \geq 0} \mathcal{D}_{i,k}^\delta.
\end{equation}

The crucial result of this section is the following lemma - adapted from~\cite[Lemma 2.7]{Asymptoticsofeven-evencorrelationsintheIsingmodel}

\begin{lemma}\label{lemme decorated random walks}
Let $\delta > 0$, and $x,y \in W$. Then, there exists $c>0$ such that:
\begin{equation}
  \PP_{(0,x)}\left[S \in \hit_{(n,y)}, \nidiam(S) \right] >   c\PP_{(0,x)} \left[S \in \mathcal{W}_n, \hit_{(n,y)}\right] .
\end{equation}
\end{lemma}
\begin{proof}
We need to prove that there exists some $c<1$ such that
\begin{equation}
\PP_{(0,x)}\left[ \exists 1 \leq i< j \leq r, \mathcal{D}(S^i) \cap \mathcal{D}(S^j) \neq \emptyset \vert S \in \mathcal{W}_n, \hit_{(n,y)}\right] \leq c.
\end{equation}
By the union bound, the latter probability is lesser or equal than 
\begin{equation}
    \sum_{1 \leq i \leq r-1}\PP_{(0,x)}\left[ \mathcal{D}(S^i) \cap \mathcal{D}(S^{i+1}) \neq \emptyset \vert 
    S \in \mathcal{W}_n, \hit_{(n,y)}  \right],
\end{equation}
and we now focus on the terms of this sum. Introduce the following family of events (recall that $\theta_k = T_{k+1}-T_k$)
:
\begin{equation}
    \mathcal{L}_k = \left\lbrace \left| S^{i+1}_k - S^i_k \right| < 2\delta \theta_{k+1} \right\rbrace.
\end{equation}
Observe that due to the cone-confinement property, if $\left\lbrace \mathcal{D}(S^{i+1}) \cap \mathcal{D}(S^i) \neq \emptyset \right\rbrace $, then one of the $\mathcal{L}_k$ must occur. 
Now we call $N$ the total number of steps. There exists a constant $\mu := \E[\theta_1]^{-1}$ such that $N \in [(\mu-\eps)n, (\mu+\eps)n]$ with exponentially large probability in $n$. For sake of simplicity, we continue the computation assuming that $N = \mu n$. Formally one should sum over all the possible values of $N$ in the latter range, but it makes no difference in the proof. We even only treat the special case $\mu=1$, as a general $\mu$ would only modify the constants inside our estimates but not the dependency in $n$. Let $T>0$ be a large integer, that will be fixed later. We first argue that there exists a constant $c_1>0$ which only depends on $\delta$ such that
\begin{equation}
    \PP_{(0,x)}\left[ \bigcup_{k=1}^{ n} \mathcal{L}_k \vert S \in \mathcal{W}_n, \hit_{(n,y)} \right] \leq \e^{+c_1T}\PP_{(0,x)}\left[ \bigcup_{k=T}^{ n-T} \mathcal{L}_k\vert S \in \mathcal{W}_n, \hit_{(n,y)}\right].
\end{equation}

This is a finite-energy property, the fact that $c_1$ is uniform over $T$ comes from the cone-confinement property. Then by union bound, let us write:
\begin{eqnarray*}
  && \PP_{(0,x)}\left[\bigcup_{k=T}^{ n-T} \mathcal{L}_k\vert S \in \mathcal{W}_n, \hit_{(n,y)} \right] 
  \\&\leq& \sum_{k=T}^{ n-T} \PP_{(0,x)}\left[ \mathcal{L}_k \vert S \in \mathcal{W}_n, \hit_{(n,y)}\right]\\
  &\leq& \sum_{k=T}^{ n-T}\sum_{\ell = 1}^{ n-k} \PP_{(0,x)} \left[ \mathcal{L}_k, \theta_{k+1} = \ell \vert S \in \mathcal{W}_n, \hit_{(n,y)}\right] \\
  &\leq& \sum_{k=T}^{ n-T}\sum_{\ell = 1}^{ n-k} \PP_{(0,x)} \left[ \left| S^{i+1}_k - S^i_k \right| < 2\delta\ell, \theta_{k+1} = \ell \vert S \in \mathcal{W}_n, \hit_{(n,y)} \right] \\
  &\leq& \sum_{k=T}^{ n-T} \sum_{\ell = 1}^{ n-k}\sum_{u \in W_{\delta\ell}}\sum_{v \in W} \e^{-c_1\ell}\e^{-c_2\norme{u-v}}\frac{q_k(x,u)q_{ n-l-k}(v,y)}{q_n(x,y)},
\end{eqnarray*}
where as in the proof of Lemma~\ref{Lemme repulsion globale sdry}, we have introduced the kernel 
\begin{equation}
    q_n(x,y) = \PP_{(0,x)}\left[ S \in \mathcal{W}_n, \hit_{(n,y)} \right],
\end{equation}
and the notation
\begin{equation}
    W_{\delta\ell} = \left\lbrace u \in W, \left|u_{i+1}-u_i\right| < \delta\ell \right\rbrace.
\end{equation}
Moreover, we also used the property that both the random variables $\theta_k$ and $\check{X}_k$ have an exponential moment. We now use the same technique as in Lemma~\ref{lemme repulsion bulk} and choose $T>0$ large enough (uniformly of everything else) to upper bound the latter quantity, using Theorems~\ref{theoreme local limit srw} and~\ref{gnedenko theorem}:

\begin{multline*}
 \PP_{(0,x)}\left[\bigcup_{k=T}^{ n-T} \mathcal{L}_k\vert S \in \mathcal{W}_n, \hit_{(n,y)} \right] 
  \\ \leq \frac{2C}{1-\eps}\sum_{k=T}^{ n/2}\left(\frac{n}{k(n-k)}\right)^{\frac{r^2}{2}}\sum_{\ell=1}^{ n-k}\e^{-c_1\ell}(\delta\ell)^2\underbrace{\sum_{u \in W_{\delta\ell}}\norme{u}^{r(r-1)-2}\e^{-\frac{\norme{u}^2}{2k}}}_{I}.
\end{multline*}
As in the proof of Lemma~\ref{lemme repulsion bulk}, we now estimate the sum $I$. Here, we will use crucially the fact that we sum over $W_{\delta\ell}$ and not over $W$. We write
\begin{eqnarray*}
I &=& \sum_{s \geq 0}\sum_{\substack{u \in W_{\delta\ell} \\ \norme{u}=s}}s^{r(r-1)-2}\e^{-\frac{r^2}{2k}} \\
&\leq&C\delta\ell \sum_{s\geq 0}s^{r-2}s^{r(r-1)-2}e^{-\frac{r^2}{2k}}.
\end{eqnarray*}
We used once again the estimation~\eqref{equation estimation volume twisted weyl chamber} for the volume of the set we are summing over. As before, we compare this sum to the integral
$I_2 = \int_{x=0}^\infty x^{r^2-4}\e^{-\frac{x^2}{2k}}\dif x$,
which, after the appropriate change of variables $t = \frac{x^2}{2k}$, can be explicitly computed, yielding 
\begin{equation}
    I_2 = \sqrt{2}^{r^2-5}\Gamma(\frac{r^2-3}{2})k^{\frac{r^2-3}{2}}.
\end{equation}
Continuing our previous computation, we obtain that:
\begin{align*}
    \PP_{(0,x)}\Big[\bigcup_{k=T}^{n-T} \mathcal{L}_k\vert S \in \mathcal{W}_n, &\hit_{(n,y)}\Big] \\ &\leq \frac{2\widetilde{C}}{1-\eps}\sum_{k=T}^{n/2}\left( \frac{n}{k(n-k)} \right)^\frac{r^2}{2} k^{\frac{r^2-3}{2}}\sum_{\ell=1}^{n-k} \e^{-c_1\ell}(\delta\ell)^3 \\
    &\leq\frac{2\widetilde{C}}{1-\eps}\sum_{k=T}^{n/2} k^{-\frac{r^2}{2} + \frac{r^2}{2} - \frac{3}{2}} \\
    &\leq \frac{2\widetilde{C}}{1-\eps} T^{-\frac{1}{2}}.
\end{align*}
Chose $T>0$ large enough so that quantity is smaller than $\frac{1}{2}$.
We then showed that:
\begin{equation}
     \PP \left[ \bigcap_{k=1}^n \mathcal{L}_k^c \vert S \in \mathcal{W}_n, \hit_{(n,y)}\right] \geq \frac{1}{2}\e^{-c_1T},
\end{equation}
which conclude the proof, since $T>0$ has been chosen uniformly of $n$.
\end{proof}

\subsection{Non-intersecting systems of decorated directed random walks}\label{subsection non intersecting non-synchronized rw}

The goal of this section is to extend the result of Section~\ref{subsection non-intersecting systems of synchronized directed random walks} to the setting of \emph{non-synchronized} random walks. For that, we will interpret such a system as an embedded synchronized random walk carrying random decorations, and use the results of the precedent section. 

Before diving into the proof, we introduce the "embedded system of synchronized random walks" of a system of random walks. 

\begin{definition}[Embedded system of synchronized random walks]
Let  $\left(\mathbf{S}_n\right)= \left( \mathbf{T}_n, \mathbf{Z}_n \right)$ be a system of non-synchronized directed random walks. Introduce the random set of synchronization times:
\begin{equation}
    \mathsf{ST} = \left\lbrace \ell\geq 0, \exists k_1(\ell), \dots, k_p(\ell) \geq 0, \mathbf{T}_{k_1}^1 = \dots = \mathbf{T}_{k_p}^r = \ell \right\rbrace.
\end{equation}
Writing $\mathsf{ST} = \left\lbrace \ell_1 < \dots < \ell_r < \dots  \right\rbrace$, we define the "embedded system of synchronized random walks" to be the process:
\begin{equation}
   \left( \check{\mathbf{S}}_n\right)_{n \geq 0} = \left( \ell_n, \mathbf{Z}^1_{k_1(\ell_n)}, \dots, \mathbf{Z}^r_{k_p(\ell_n)} \right)_{n \geq 0}.
\end{equation}
Observe that in particular the trajectory of $\check{\mathbf{S}}$ is a subset of the trajectory of $\mathbf{S}$, and that by definition the system $\check{\mathbf{S}}$ is synchronized. 
\end{definition}

\begin{lemma}
The process $\check{\mathbf{S}}$ is a synchronized system of random walks (recall Definition~\ref{def synchronized system of random walks}). Moreover, for any $i \in \{1, \dots, r \}$, $\mathcal{D}(\mathbf{S}^i) \subset \mathcal{D}(\check{\mathbf{S}}^i)$.
\end{lemma}

\begin{proof}
All the statements are easy to check, the exponential tails of the length being a consequence of the Renewal Theorem of~\cite{feller15}. 
\end{proof}

\begin{lemma}\label{lemme estimation proba non-synchronized RW}
There exists a positive $c>0$ such that for any fixed $x,y \in W$,
\begin{equation}
    \mathbf{P}_{(0,x)}\left[\nidiam(\mathbf{S}), \hit_{(n,y)}\right] > c\frac{V(x)V(y)}{n^{\frac{r^2}{2}}}.
\end{equation}
\end{lemma}
\begin{proof}
Observe that 
\begin{equation}
    \mathbf{P}_{(0,x)}\left[\nidiam(\mathbf{S}), \hit_{(n,y)}\right] \geq \PP_{(0,x)}\left[\nidiam(\check{\mathbf{S}}), \hit_{(n,y)}\right],
\end{equation}
so that we focus on lower bounding the right-hand side.
We are in the setting of Lemma~\ref{lemme decorated random walks}, allowing us to write:
\begin{eqnarray*}
\PP_{(0,x)}\left[\nidiam(\check{\mathbf{S}}), \hit_{(n,y)}\right]
&\geq& c\PP_{(0,x)}\left[\check{\mathbf{S}} \in  \mathcal{W}_n, \hit_{(n,y)}\right] \\
&\geq& cV(x)V(y)n^{-\frac{r^2}{2}},
\end{eqnarray*}
where the first inequality comes from Lemma~\ref{lemme decorated random walks} and the second one comes from the fact that $\check{S}$ has the distribution of a synchronized system of random walks, so that Theorem~\ref{theoreme local limit srw} applies. 
\end{proof}

\begin{remark}
The exact same technique of proof can be used to show an analog of Lemma~\ref{Lemme repulsion globale sdry} for non-synchronized random walks. Indeed the probability of two non-synchronized random walks coming close one from each other can be upper bounded by the probability of two decorated synchronized random walks coming close one from each other. Making use of Lemmas~\ref{Lemme repulsion globale sdry} and~\ref{lemme estimation proba non-synchronized RW} we obtain:
\begin{lemma}\label{lemme repulsion globale non-synchronized RW}
Let $\mathbf{S}$ be a system of non-synchronized random walks. Let $x, y$ two sequences of elements of $W$ such that 
\begin{equation}
    \Gap(x), \Gap(y) \geq  n^\eps
\end{equation}
and
\begin{equation}
    \norme{x}, \norme{y} = o(\sqrt{n}).
\end{equation}
Then for any $\delta > 0$ sufficiently small,, there exists $\beta > 0$, $C> 0$ such that for $n\geq 0$ large enough, 
\begin{equation}
    \mathbf{P}_{(0,x)} \left[ \inf_{1\leq k\leq n} \Gap(\mathbf{S}_k) \leq n^\delta~\big\vert ~\mathbf{S} \in \hit_{(n,y)}, \nidiam(\mathbf{S})\right] \leq Cn^{-\beta}.
\end{equation}
\end{lemma}

\end{remark}

The next step in our way to the proof of Theorem~\ref{theoreme invariance principle for drw} is then to show a fast repulsion estimate near the starting and ending points stated in Lemma~\ref{lemme entropic repulsion sdrw} in the setting of non-synchronized systems of non-intersecting bridges. Let $\eps > 0$. As in Sections~\ref{section independent system} and~\ref{section RCM} we introduce the following times:
\begin{equation}
    T_1(\mathbf{S}) = \min_{k \geq 0} \left\lbrace k \geq 0, \Gap(\mathbf{S}_k) > n^\eps  \right\rbrace 
\end{equation}
and
\begin{equation}    
    T_2(\mathbf{S}) = \max_{k \geq 0} \left\lbrace k \geq 0, \Gap(\mathbf{S}_k) > n^\eps  \right\rbrace .
\end{equation}

\begin{lemma}\label{entropic repulsion non-synchronized RB}

There exists $\eps >0$ sufficiently small such that there exists a positive constant $c>0$ such that:
\begin{equation}
    \mathbf{P}_{(0,x)} \left[ T_1(\mathbf{S}) > n^{1-\eps}, T_2(\mathbf{S}) < n-n^{1-\eps} \vert \mathbf{S} \in \hit_{(n,y)}, \nidiam(\mathbf{S}) \right] < \frac{1}{c}\exp\left(-cn^\eps\right),
\end{equation}
where $\mathbf{P}_{(0,x)}$ is the distribution of a non-synchronized system of directed random walks. 
\end{lemma}
\begin{proof}
Recall that $\check{\mathbf{S}}$ denotes the synchronized system of random walks embedded in $\mathbf{S}$. Then, 
\begin{align*}
 \mathbf{P}_{(0,x)} \big[ T_1(\mathbf{S})& > n^{1-\eps},  T_2(\mathbf{S}) < n-n^{1-\eps} \vert  \mathbf{S} \in \hit_{(n,y)}, \nidiam(\mathbf{S})\big] \\
&\leq \mathbf{P}_{(0,x)} \left[ T_1(\check{\mathbf{S}}) > n^{1-\eps}, T_2(\check{\mathbf{S}}) < n-n^{1-\eps} \vert  \mathbf{S} \in \hit_{(n,y)}, \nidiam(\mathbf{S})\right] \\
&\leq 2 \frac{\mathbf{P}_{(0,x)}\left[T_1(\check{\mathbf{S}}) > n^{1-\eps}  \right]}{\mathbf{P}_{(0,x)}\left[\mathbf{S} \in \hit_{(n,y)}, \nidiam(\mathbf{S})\right]} \\
&\leq \frac{2}{c}\exp(-cn^\eps)n^{\frac{r^2}{2}}(V(x)V(y))^{-1},
\end{align*}
which proves the lemma for another constant $c'<c$ provided that $n$ is large enough. The last inequality comes from Lemma~\ref{lemme entropic repulsion sdrw} for upper bounding the numerator, and from Lemma~\ref{lemme estimation proba non-synchronized RW} for lower bounding the denominator. 
\end{proof}
We are now ready to prove Theorem~\ref{theoreme invariance principle for drw}. The technique is very similar to the proof of Theorem~\ref{Theoreme main}. Indeed, we shall wait for a  sublinear time that the walks attain a gap of order $n^\eps$. After this time, we know that - looking at the process as a system of synchronized decorated random walks - the diamonds are very likely not to intersect so that the convergence of the synchronized embedded system towards the Brownian watermelon can be transmitted to the whole system.  
\begin{proof}[Proof of Theorem~\ref{theoreme invariance principle for drw}]
Let $\mathbf{S}$ be sampled according to the measure
\begin{equation}
\mathbf{P}_{(0,x)}\left[ ~.~\big\vert \mathbf{S} \in \mathsf{NonIntDiam}(\mathbf{S}), \hit_{(n,y)}\right].
\end{equation}
Again, since we are going to work between the random times $T_1$ and $T_2$, we need to implement the strategy given by Lemma~\ref{Lemme technique processus sto}. Let $\delta > 0$ and $f^\delta : \mathcal{C}([\delta, 1-\delta], \R^r) \rightarrow \R$, continuous and bounded. Introduce $\mathbf{S}_n(t)$ the scaled version of $\mathbf{S}$:
\begin{equation}
    \mathbf{S}_n(t)= \frac{1}{\sqrt{n}}\mathbf{S}(nt).
\end{equation}
Our goal is to show that (we keep implicit the restrictions of $\mathbf{S}_n$ and $\bw^{(r)}$ to the interval $[\delta, 1-\delta]$):
\begin{equation}
    \E\left[f^\delta(\mathbf{S}_n) \vert \mathbf{S} \in \nidiam(\mathbf{S}), \hit_{(n,y)} \right] \goes{}{n \rightarrow \infty}{\E\left[f^\delta(\sigma\bw^{(r)})\right]}.
\end{equation}

We claim that - thanks to Lemma~\ref{entropic repulsion non-synchronized RB} and the usual deviation argument for random walks - with probability $1+o(1)$, there exist $T_1>0$ and $T_2< n$ two random times such that $T_1$ and $T_2$ are synchronization times for $\mathbf{S}$, and such that 
\begin{equation}\label{conditions t1 t2 marches}
    T_1 < 2n^{1-\eps} \text{ and  } T_2 > n - 2n^{1-\eps},
\end{equation}
and 
\begin{equation}\label{condition ecartement marches aleatoires}
    \begin{cases}
    \norme{\mathbf{S}(T_1)}, \norme{\mathbf{S}(T_2)} = o(\sqrt{n}), \\
        \min\lbrace \Gap(\mathbf{S}(T_1)), \Gap(\mathbf{S}(T_2)) \rbrace > \frac{1}{2}n^\eps.
    \end{cases} 
\end{equation}
In the rest of the proof, we then condition on the values of $T_1, T_2, \mathbf{S}(T_1)$ and $\mathbf{S}(T_2)$ satisfying~\eqref{conditions t1 t2 marches} and~\eqref{condition ecartement marches aleatoires}. Moreover for sake of simplicity in the proof let us call $u = (T_1, \mathbf{S}(T_1))$ and $v=(T_2,\mathbf{S}(T_2))$. As soon as $n$ is large enough so that $n\delta>T_1$ and $n(1-\delta) < T_2$, the Markov property for random walks ensures that:

\begin{equation}
    \E\Big[f^\delta(\mathbf{S}_n \vert u,v, \mathbf{S} \in \nidiam(\mathbf{S}), \hit_{(n,y)} \Big]  = \E_{u}\Big[ f^\delta(\mathbf{S}_n \vert \mathbf{S} \in \nidiam(\mathbf{S}), \hit_{v} \Big],
\end{equation}
where $\E_u$ denotes the expectation under the measure $\mathbf{P}_{u}$.

Let us consider $\check{\mathbf{S}}$ to be the synchronized system embedded into $\mathbf{S}$, and $\check{\mathbf{S}}(t)$ be its linear interpolation. By standard estimates on the max of a linear number of independent random variables with exponential tails, one gets: 
\begin{align*}
    \mathbf{P}_{(0,x)} \Big[ \sup_{0\leq t \leq n}\Big| \mathbf{S}(t) - \check{\mathbf{S}}(t) \Big| & > \log^2 n \vert \mathbf{S}\in \nidiam(\mathbf{S}), \hit_{(n,y)}  \Big] \\ &\leq \frac{\mathbf{P}_{(0,x)} \left[\sup_{0\leq t \leq n}\left| \mathbf{S}(t) - \check{\mathbf{S}}(t) \right| > \log^2 n\right]}{\mathbf{P}_{(0,x)}\left[\mathbf{S} \in \nidiam(\mathbf{S}), \hit_{(n,y)}\right]} \\
    &\leq \frac{1}{c}\exp\left(-c(\log^2 n)\right)(V(x)V(y))^{-1}n^{\frac{r^2}{2}},
\end{align*}
where we used Lemma~\ref{lemme estimation proba non-synchronized RW} for the last step. We now work under the event that
\begin{equation}
    \sup_{0\leq t \leq n}\left| \mathbf{S}(t) - \check{\mathbf{S}}(t) \right| > \log^2 n.
\end{equation}
Hence, for our purpose it is sufficient to show that:
\begin{equation}
    \E_{u}\big[ f^\delta(\check{\mathbf{S}}_n(t)) \vert \mathbf{S} \in \hit_{v}, \nidiam(\mathbf{S}) \big]  \goes{}{n\rightarrow \infty}{\E\left[f^\delta(\sigma\bw^{(r)})\right]}.
\end{equation}
The next step is to replace the conditioning over $\mathbf{S}$ belonging to the non-intersection of diamonds and connection event by a conditioning over $\check{\mathbf{S}}$ belonging to the non-intersection and connection event. Indeed, assuming that we managed to show that this change of conditioning was justified, the result would follow by Theorem~\ref{theoreme d'invariance sdrw}.
Our target estimate is then:
\begin{equation}
    \mathbf{P}_u\left[ \left\lbrace \nidiam(\mathbf{S}) \right\rbrace \Delta \left\lbrace \check{\mathbf{S}}\in \mathcal{W}_{T_2-T_1} \right\rbrace \vert \check{\mathbf{S}}\in \mathcal{W}_{T_2-T_1}, \hit_v \right]  \goes{}{n \rightarrow \infty}{0}.
\end{equation}
Observe that because we work under the event $\lbrace\sup_{0\leq t \leq n}\left| \mathbf{S}(t) - \check{\mathbf{S}}(t) \right| > \log^2 n \rbrace$, then
\begin{multline}
  \mathbf{P}_u\left[ \left\lbrace \mathbf{S} \in \nidiam(\mathbf{S}) \right\rbrace \Delta \left\lbrace \check{\mathbf{S}}\in \mathcal{W}_{T_2-T_1} \right\rbrace \right] \\ \leq \mathbf{P}_u \big[ \inf_{t \in [T_1,T_2]} \Gap(\check{\mathbf{S}}(t)) < 4\log^3 n \vert \check{\mathbf{S}}\in\mathcal{W}_{[T_1, T_2]}, \hit_v\big]. 
 \end{multline}
 By Lemma~\ref{Lemme repulsion globale sdry} we know that this probability decays to 0 at least polynomially fast. Thus, we proved that 
 \begin{equation}\label{equivalence des conditionnements preuve marches}
     \Big|\E_{u}\Big[ f^\delta(\check{\mathbf{S}}_n(t)) \vert \nidiam(\mathbf{S}), \hit_v \Big] -  \E_{u}\big[ f^\delta(\check{\mathbf{S}}_n(t)) \vert \check{\mathbf{S}}\in\mathcal{W}_{[T_1, T_2]}, \hit_v \big] \Big|  \goes{}{n\rightarrow \infty}{0}.
 \end{equation}
 Now because of~\eqref{conditions t1 t2 marches} and~\eqref{condition ecartement marches aleatoires}, Theorem~\ref{theoreme d'invariance sdrw} applies and we get that
 \begin{equation}
     \E_{u}\big[ f^\delta(\check{\mathbf{S}}_n(t)) \vert \check{\mathbf{S}}\in\mathcal{W}_{[T_1, T_2]}, \hit_v \big]  \goes{}{n\rightarrow \infty}{\E\left[f^\delta(\sigma\bw^{(r)})\right]}
 \end{equation}
 for some $\sigma> 0$. This concludes the proof of the theorem: as previously condition $(i)$ of Lemma~\ref{Lemme technique processus sto} is a simple consequence of basic large deviations estimates. 
 \end{proof}

\section*{Acknowledgements}
We warmly thank Ioan Manolescu and Sébastien Ott for very useful and instructive discussions. We thank Romain Panis, Ulrik Thinggaard Hansen and Maran Mohanarangan for a careful reading of an early draft of this paper. The author was supported by the Swiss National Science Foundation grant n° 182237.

\bibliographystyle{amsalpha}
\bibliography{biblio_multiple_interfaces.bib}
\section*{Appendix}
\subsection*{Non-confinement in small tubes for a single non-degenerate directed random walk}
\begin{lemma}[Non-confinement of single non-degenerate directed random walk]\label{Confinement lemma}
    There exists $\eps_0> 0$ such that that the following holds. Fix $\eps<\eps_0$. Let $(S_n)_{n\geq 0}$ be a non-degenerate directed random walk, and remember that $S(t)$ denotes its linear interpolation. Let $f: \R^+ \rightarrow \R$ be any function.    
    Then, for any $\alpha \in (0, 1]$, there exists $c>0$ such that for any $x \in \R$,
    \begin{equation}\label{eq:RWconfinement }
        \PP_{(0,x)}\left[ \# \left\lbrace k \in \lbrace 0, \dots, \lfloor n^{1-\eps} \rfloor \rbrace, \left| S(k) - f(k) \right| < n^\eps \right\rbrace > \alpha n^{1-\eps} \right] < \e^{-cn^{1-3\eps}}.
    \end{equation}
\end{lemma}

Limiting the times considered in~\eqref{eq:RWconfinement } to $k\leq \lfloor n^{1-\eps}\rfloor$ rather than the more natural choice $k \leq n$ is done only for coherence with the uses of this statement in other parts of the paper. 

\begin{proof}
In the following proof we reason up to integer rounding for the time indexes (so that $n^\eps$ might be used instead of $\lfloor n^\eps \rfloor$).  
We cut up the interval $\lbrace 0,\dots, n^{1-\eps} \rbrace$ in intervals of alternating lengths 
$C \tfrac12 \alpha n^{2\eps}$ and $C (1 - \tfrac12 \alpha)n^{2\eps}$ (where $C$ is some fixed constant to be determined).
Call these buffer and main intervals. The buffer intervals occupy a proportion $\alpha/2$ of the whole walk, so
\begin{multline*}
\PP_{(0,x)}\left[ \# \left\lbrace k \in \lbrace 0, \dots, n^{1-\eps}\rbrace, \left| S(k) - f(k) \right| < n^\eps \right\rbrace > \alpha n^{1-\eps} \right] \\
\leq 
 \PP_{(0,x)}\left[ \# \{ k \in \text{main intervals }, \left| S(k) - f(k) \right| < n^\eps \}  > \tfrac12\alpha n^{1-\eps} \right].
\end{multline*}
Call the indices $k$ considered above  ``close points''. Call a main interval \textit{bad} if it has a proportion of close points larger than $\alpha/4$. Then, for the above to be realized, one needs a proportion of at least $\alpha/4$ bad main intervals (the good main intervals account for at most $\tfrac14\alpha n^{1-\eps}$ close points). 
\\
Condition now on the trajectory in each of main interval. The only randomness comes from the starting positions of these main intervals, which are dictated by the buffer intervals.

One can then check that due to the pigeonhole principle, for each main interval, there are at most $\tfrac4\alpha  \times n^{\eps}$ starting positions that render them bad. 
Thus, just because of the buffer interval preceding each main interval, due to the Central limit Theorem, the probability of a main interval to be bad may be rendered small (smaller than any given constant, by choosing $C$ large enough). Thus choose $C$ so that 
\begin{align*}
    \PP[\text{ one main interval is bad } | \text{ all the RW except the preceding buffer interval} ] \leq \alpha/8.
\end{align*}
In total we have have $\frac1C n^{1- 3\eps}$ pairs of buffer and main intervals. Each main interval has a probability at least $1 - \alpha/8$ to be good, independently of all other. Thus, the probability of having a proportion $\alpha/4$ of bad intervals is a large deviation estimate, and thus has a probability of order $e^{- c n^{1- 3\eps}}$ for some constant $c$ that depends on $C$ (itself depending on $\alpha$).
\end{proof}

\subsection*{A brief sketch of proof of Theorems~\ref{theoreme local limit srw},~\ref{gnedenko theorem} and~\ref{theoreme d'invariance sdrw}}

As explained previously, Theorems~\ref{theoreme local limit srw},~\ref{gnedenko theorem} and~\ref{theoreme d'invariance sdrw} have already been proved in the case of regular random walks in \cite{conditionallimittheoremsfororderedrandomwalks, randomwalksincones,Invarianceprinciplesforrandomwalksincones}, under a weaker moment assumption and the assumption that the coordinates of the walk are exchangeable --- which suits to our setting. It has already been explained in~\cite{ottvelenikwachtelioffe} how to transfer Local limit Theorems proved for regular random walks to the case of directed random walks, and the same method applies \textit{mutatis mutandi} to our setting. 
Indeed, the proof consists in conditioning on the number of steps of the walk called $N$, an considering three different cases. Indeed, large deviation estimates allow to rule out the case $N \notin [(\mu \pm \eps)n]$, with $\mu := \E[\theta_1]^{-1}$. Then, the contribution of the indexes $N \in [(\mu - \eps)n, \mu n - A\sqrt n ] \cup [\mu n + A\sqrt n, (\mu + \eps)n]$ is shown to be of order $f(A)n^{-\frac{r^2}2}V(x)V(y)$, with $f(A) \rightarrow 0$ when $A \rightarrow \infty$. As explained in~\cite{ottvelenikwachtelioffe}, the important idea is to perform an exponential tilt of the random walk by the length of its time increments, and to analyze this new measure by standard random walks estimates. The proof finally reduces to the case where $N$ lies in the interval $[\mu n \pm A\sqrt n]$, where $A$ is a large constant. The proofs of~\cite{randomwalksincones} can then be mimicked. The discussion of~\cite[Proof of Thm 5.1]{ottvelenikwachtelioffe} --- in particular the observation that the harmonic function $V$ does not depend on the time reference --- shows that the result is uniform in starting points satisfying $\Vert x \Vert, \Vert y \Vert = o(\sqrt n)$.  
It then may be considered as folklore that Theorems~\ref{theoreme local limit srw},~\ref{gnedenko theorem} and~\ref{theoreme d'invariance sdrw} do hold in the case of synchronized directed random walks.

\end{document}